\documentclass[microtype]{gtpart}     
\gtart  
\title{Intersection theory of punctured pseudoholomorphic curves}

\author{Richard Siefring}
\givenname{Richard}
\surname{Siefring}
\address{Department of Mathematics \\
    Michigan State University \\
    East Lansing, MI 48824}
\address{Max Planck Institute for Mathematics in the Sciences \\
     Inselstra{\ss}e 22 \\
     D-04103 Leipzig \\
     Germany}
\email{siefring@mis.mpg.de}
\urladdr{http://personal-homepages.mis.mpg.de/siefring/}


\arxivreference{0907.0470}
\arxivpassword{n3uf4}

\volumenumber{}
\issuenumber{}
\publicationyear{}
\papernumber{}
\startpage{}
\endpage{}
\doi{}
\MR{}
\Zbl{}
\received{}
\revised{}
\accepted{}
\published{}
\publishedonline{}
\proposed{}
\seconded{}
\corresponding{}
\editor{}
\version{}

\newtheorem{theorem}{Theorem}[section]   
\newtheorem{lemma}{Lemma}[section]          
\newtheorem{corollary}{Corollary}[section]
\newtheorem{proposition}{Proposition}[section]
\theoremstyle{definition}

\newtheorem{remark}{Remark}[section]

\numberwithin{equation}{section}

\makeatletter
  \let\c@lemma=\c@theorem
\makeatother

\makeatletter
  \let\c@corollary=\c@theorem
\makeatother

\makeatletter
  \let\c@proposition=\c@theorem
\makeatother

\makeatletter
  \let\c@definition=\c@theorem
\makeatother

\makeatletter
  \let\c@remark=\c@theorem
\makeatother

\newcommand{\N}{\mathbb{N}}
\newcommand{\ep}{\varepsilon}
\newcommand{\dbar}{\bar\partial}
\newcommand{\tl}{\tilde}
\newcommand{\wt}{\widetilde}
\newcommand{\ds}{\partial_s}
\newcommand{\dt}{\partial_t}
\newcommand{\X}{X_{\Ha}}
\newcommand{\A}{\mathbf{A}}
\newcommand{\wrel}{\wind_{rel}}
\newcommand{\wtwo}{\wind_{2}}
\newcommand{\winfty}{\wind_{\infty}}
\newcommand{\M}{\mathcal{M}}
\newcommand{\Ha}{\mathcal{H}}
\newcommand{\xh}{\xi^{\Ha}}
\renewcommand{\J}{\mathcal{J}}
\newcommand{\Pe}{\mathcal{P}}
\newcommand{\PT}{\widetilde\Pe^{0}_{\tau}}
\newcommand{\PMH}{\widetilde\Pe^{0}(M,\Ha)}
\newcommand{\PTMH}{\PT(M,\Ha)}
\newcommand{\cover}{\bar\sigma}
\newcommand{\DD}{\mathbb{D}}
\newcommand{\Done}{\Delta_{1}}
\newcommand{\Dtwo}{\Delta_{2}}
\newcommand{\Dtot}{\Delta_{tot}}
\newcommand{\MH}{(M, \Ha)}
\newcommand{\MW}{\M(W, \barJ, \Ha^{+}, \Ha^{-})}
\newcommand{\MM}{\M(M, \Ha, J)}
\newcommand{\CW}{C^{\infty}(W, \Ha^{+}, \Ha^{-})}
\newcommand{\CM}{C^{\infty}(M, \Ha)}
\newcommand{\Mhar}{\M_{\Delta}(M, \Ha, J)}

\newcommand{\wind}{\operatorname{wind}}
\newcommand{\cov}{\operatorname{cov}}
\newcommand{\Hom}{\operatorname{Hom}}
\newcommand{\sing}{\operatorname{sing}}
\newcommand{\ind}{\operatorname{ind}}
\newcommand{\inum}{\operatorname{int}}
\newcommand{\rot}{\operatorname{rot}}
\newcommand{\rin}{i} 
\newcommand{\ain}{i_{\infty}}
\newcommand{\lne}{\sigma^{-}_{max}}

\def\gin#1#2{[#1]*[#2]}

\def\c#1{[#1,\infty)\times S^{1}}
\def\ip#1{\langle#1\rangle}
\def\bp#1{\left(#1\right)}
\def\bbr#1{\left[#1\right]}
\def\br#1{\left\{#1\right\}}
\def\fl#1{\lfloor#1\rfloor}
\def\ceil#1{\lceil#1\rceil}
\def\Mx{\overline\Lambda}
\def\Mxx#1{\Mx\bp{#1}}

\newcommand{\tildeJ}{\wt J}
\newcommand{\barJ}{\overline J}

\newcounter{ham}\renewcommand{\theham}{(H\arabic{ham})}
\newcounter{ham2}
\newenvironment{hamlist}
 {
 \begin{list}{\theham}{\usecounter{ham}\setcounter{ham}{\value{ham2}}}
 }
{\setcounter{ham2}{\value{ham}}\end{list}}

\begin{document}
\begin{abstract}    
We study the intersection theory of punctured pseudoholomorphic curves in $4$-dimensional
symplectic cobordisms.  Using the asymptotic results from
\cite{sief08}, we first study the local intersection properties of
such curves at the punctures.  We then
use this to develop topological controls on the intersection number of two curves.
We also prove an adjunction formula which gives a topological condition that will guarantee
a curve in a given homotopy class is embedded, extending previous work of Hutchings from \cite{hut02}.

We then turn our attention to curves in the symplectization $\R\times M$ of a $3$-manifold
$M$ admitting a stable Hamiltonian structure.
We investigate controls on intersections of the projections of curves to the $3$-manifold, and
we present conditions that will guarantee the projection of a curve to the $3$-manifold is an embedding.

Finally we consider an application concerning pseudoholomorphic curves in
manifolds admitting a certain class of holomorphic open book decomposition, and
an application concerning the existence of generalized pseudoholomorphic curves, as introduced
by Hofer in \cite{h2000}.
\end{abstract}

\maketitle


\tableofcontents

\section{Introduction}
In this paper we will study the intersection theory of punctured pseudoholomorphic curves
which arise in symplectic field theory  \cite{egh}.
Positivity of intersections for pseudoholomorphic curves has been an important tool in applications of
pseudoholomorphic curves to $4$-dimensional symplectic topology.  First stated by Gromov in
\cite{gr}, rigorous proofs were subsequently provided by McDuff \cite{mcd94},
and Micallef and White \cite{mw}.  Put simply, positivity of intersections states
that isolated intersections between two curves contribute positively to their intersection number, and that
singular points contribute positively to the self-intersection number of a single curve.
These local statements translate into useful global
topological statements for closed pseudoholomorphic curves.
If $u:(\Sigma, j)\to (M, J)$ and $v:(\Sigma', j')\to (M, J)$ are pseudoholomorphic maps with
closed connected domains and nonidentical images, then the intersection number $[u]\cdot[v]$
is always nonnegative and $[u]\cdot[v]=0$ if and only if $u$ and $v$ do not intersect.
Moreover, if $u:(\Sigma, j)\to (M, J)$ is a simple, closed, pseudoholomorphic map, then
$u$ satisfies the inequality
\begin{equation}\label{e:adj-closed}
[u]\cdot[u]-\ip{c_{1}(TM, J), [u]}+\chi(\Sigma)\ge 0
\end{equation}
and equality occurs if and only if $u$ is an embedding.

While the local results on the intersections of pseudoholomorphic curves apply to punctured curves
appearing in symplectic field theory, coming up with generalizations of the global results is subtle
due to the fact that the intersection number is no longer homotopy invariant when the domains of the
curves are noncompact.  Indeed, in a pair of smooth homotopies of curves,
intersections can escape or appear at the ends.
One way to deal with this issue is to perturb one of the curves being considered near the ends, and compute
the intersection number between one curve and the perturbation of the second curve.
Given a sufficiently precise description of the asymptotic behavior, it is then possible to
compute this intersection number in terms of the intersection number of the original curves and
behavior near the punctures.
This idea was studied by Kriener in \cite{kr}, where the self-intersection
number of a single embedded half-cylinder asymptotic to a multiply covered orbit is considered.
These ideas were further pursued by Hutchings in \cite{hut02} and \cite{hut09},
where he proves an index inequality for curves in symplectizations.
This index inequality, important for the foundations of embedded contact homology,
gives a topological criterion that will guarantee a curve is embedded and has
asymptotic behavior which satisfies a technical ``admissibility'' condition.

The goal of the present paper is to further develop some of the techniques used in \cite{hut02},
and study algebraic controls on intersections and embeddedness of punctured pseudoholomorphic curves.
We first give a complete study of the local ``asymptotic intersection theory,''
made possible by the asymptotic descriptions of curves in \cite{sief08}.  With this in hand, we introduce the
notion of the ``generalized intersection number'' of two smooth proper maps from punctured Riemann surfaces
which are asymptotic at the punctures to cylinders over periodic orbits.  From the results we prove about
asymptotic intersection theory, it will follow that the generalized intersection number of two distinct
pseudoholomorphic maps
is always nonnegative, and is equal to zero if and only if the curves do not intersect \emph{and} the
curves do not have any tangencies at infinity,
where the notion of tangency at infinity can be made precise in terms of the
asymptotic description from \cite{sief08}.
Moreover, we state a generalization of the adjunction formula
\eqref{e:adj-closed} in terms of the generalized intersection number.

Having addressed these things we specialize to the case of a cylindrical cobordism
$\R\times M$
equipped with an $\R$-invariant almost
complex structure.  Here we relate the number of intersections of two curves with the number of intersections
of each curve with the asymptotic limits of the other,
and the winding of the curves around their asymptotic limits.
This allows us to state a set of necessary and sufficient conditions that will guarantee that the projection of the
two curves to the $3$-manifold $M$ do not intersect.
These same techniques, used with the adjunction formula for punctured curves, allow us to
give conditions that will guarantee the projection of a curve to the 3-manifold is embedded.
These results and ideas are useful in the study of finite energy foliations, as initiated by Hofer, Wysocki, and Zehnder in \cite{hwz:fol}.
As an
application of these results, we prove a result about the contact homology of a
$3$-manifold admitting a special class of holomorphic open book decompositions.

Finally, we consider so-called generalized pseudoholomorphic curves as introduced by Hofer in \cite{h2000}.
We show that the generalized intersection product can be used to develop topological obstructions to
the existence of these curves.

\subsubsection*{Acknowledgements}
This paper has its roots in work I began as a graduate student, and I would
like to thank my advisor, Helmut Hofer, for his encouragement and support.

During the writing of this paper, I have had many in-depth conversations about this material
which have helped to clarify my thinking and have shaped the exposition here.
In particular I would like to thank Barney Bramham, Michael Hutchings, Al Momin, 
Eric Schoenfeld, and Chris Wendl for helpful conversations and for their interest in this work.

\section{Background and main results}
\subsection{Hamiltonian structures}\label{ss:ham}
In this section we describe a structure defined on
$3$-manifolds that is central to all we do here.
For further background and examples, see Section 2 of \cite{behwz} or Example 4.2 in \cite{ekp06}.

Let $M$ be a compact oriented $3$-manifold equipped with a
pair $\Ha=(\lambda, \omega)$ where $\lambda$ is a $1$-form and
$\omega$ is a
$2$-form on $M$.
Assume that
\begin{hamlist}
\item\label{h:vform} $\lambda\wedge\omega$ is a volume form on $M$.
\end{hamlist}
Then $\omega$ must be rank $2$ everywhere, and thus defines a line bundle
$\ell_{\omega}\subset TM$ by
\[
\ell_{\omega}=\cup_{p\in M}(p,\ker\omega_{p})
\]
where $\ker\omega_{p}$ is the kernel of the linear map
$\omega_{p}:TM_{p}\to T^{*}M_{p}$ defined by
$v\mapsto \omega_{p}(v, \cdot)$.
The condition \ref{h:vform} implies that $\lambda$ is non-zero on $\ell_{\omega}$, that
the hyperplane distribution defined by
\[
\xh:=\ker\lambda
\]
is everywhere transverse to
$\ell_{\omega}$, and that $\omega$ is nondegenerate on $\xh$.  If we define a vector field $\X$
to be the unique section of $\ell_{\omega}$ satisfying $\lambda(\X)=1$, we see that condition
\ref{h:vform} implies that the pair $(\lambda, \omega)$ determines a splitting
\begin{equation}\label{e:splitting}
TM=(\ell_{\omega},\X)\oplus(\xh, \omega)
\end{equation}
of the tangent space of $M$ into a framed line bundle $(\ell_{\omega}, \X)$ and a symplectic
$2$-plane bundle $(\xh, \omega)$.

If we further require that
\begin{hamlist}
\item\label{h:closed} $\omega$ is closed
\end{hamlist}
then we can conclude that any section $v\in\Gamma(\ell_{\omega})$ of $\ell_{\omega}$,
in particular $\X$,
satisfies
\[
L_{v}\omega=i_{v}d\omega+d(i_{v}\omega)=0
\]
so the flow of any section of $\ell_{\omega}$ preserves $\omega$.
Finally, if we require that
\begin{hamlist}
\item\label{h:kernel} $d\lambda$ vanishes on $\ell_{\omega}$,
\end{hamlist}
then we find that
\[
L_{\X}\lambda=i_{\X}d\lambda+d(i_{\X}\lambda)=0
\]
so $\lambda$ is preserved by the flow of $\X$, and hence the splitting
\eqref{e:splitting} is also preserved by the flow of $\X$.
Following \cite{ekp06}, we will refer to a pair $\Ha=(\lambda, \omega)$
satisfying \ref{h:vform}--\ref{h:kernel}  as a
{\em stable Hamiltonian structure} on $M$.
We will refer to the vector field $\X$ arising from such a structure
as the {\em Reeb vector field} associated to $\Ha$, 
and we will refer the the hyperplane bundle $\xh$ as the
{\em hamiltonian hyperplane field} associated to $\Ha$.

In what follows the dynamics of the vector field $\X$ will play an important role,
and the periodic orbits of $\X$ will be of particular interest.  For our purposes, it will be convenient to
think of periodic orbits as maps parametrized by $S^{1}\approx\R/\Z$ equipped with the basepoint
$0\in\R/\Z$.  More precisely, for $\tau> 0$ we define the set $\PTMH$ of
{\em simple $\tau$-periodic orbits} of $\X$ to be the set of all $\gamma\in C^{\infty}(S^{1}, M)$,
such that $\gamma$ is an embedding and $\gamma$ satisfies the equation
\[
d\gamma(t)\dt=\tau\X(\gamma(t))
\]
for all $t\in S^{1}$.
We will denote the set of all simple periodic orbits by
\[
\PMH:=\cup_{\tau> 0}\PT(M,\Ha).
\]
We note that each set $\PT(M,\Ha)$ is invariant under the
$S^{1}$-action on $C^{\infty}(S^{1}, M)$ defined by
$c*\gamma(t)=\gamma(t+c)$ for $c\in\R/\Z$ and $\gamma\in C^{\infty}(S^{1},M)$.
We define the space of {\em unparametrized simple periodic orbits}
\[
\Pe^{0}\MH=\PMH/S^{1}
\]
to be the space of $S^{1}$-orbits in $\PMH$.
For some of the results we present (particularly the asymptotic results in
\fullref{ss:asymp}) the choice of basepoint $0\in\R/\Z$ will be important for precise statements.
Otherwise, we will
generally want to think of two simple orbits as the same if they belong to the same
class in $\Pe^{0}\MH=\PMH/S^{1}$.  In what follows,
we will use the same notation for an orbit $\gamma\in\PMH$ and its equivalence class in
$\Pe^{0}\MH$, and if we write
$\gamma_{2}=\gamma_{1}$ for two simple periodic orbits, the ``$=$'' should be interpreted
modulo the $S^{1}$-action on $\PMH$ unless otherwise stated.

Given a $\gamma\in\PMH$, we define for any $m\in\Z\setminus\br{0}$
a map $\gamma^{m}\in C^{\infty}(S^{1}, M)$ by
\[
\gamma^{m}=\gamma\circ p_{m}
\]
where $p_{m}:S^{1}\to S^{1}$ is the $m$-fold covering map defined by
$p_{m}([t]_{\R/\Z})=[mt]_{\R/\Z}$.
We refer to $\gamma^{m}$ as a {\em multiply covered periodic orbit} with multiplicity $m$, and we
denote the set of all periodic orbits (simple and multiply covered) by
\[
\Pe(M,\Ha).
\]

Let $\psi:\R\times M\to M$ be the flow of the $\X$, i.e.\
\[
\dot\psi_{t}(x)=\X(\psi_{t}(x))
\]
for all $(t,x)\in\R\times M$.  If $\gamma\in\PMH$ is a
simple $T$-periodic orbit, then 
it follows from
the fact that $L_{\X}\lambda=0$ and $L_{\X}\omega=0$ observed earlier
that $d\psi_{mT}|_{\xh_{\gamma(0)}}\in Sp(\xh_{\gamma(0)}, \omega_{\gamma(0)})$
for all $m\in\Z$.
We say that the periodic orbit $\gamma^{m}$ is {\em nondegenerate} if
$d\psi_{mT}|_{\xh_{\gamma(0)}}$ does not have $1$ in its spectrum.
We will say a Hamiltonian structure $\Ha=(\lambda, \omega)$ is nondegenerate if all
periodic orbits of the corresponding vector field $\X$ are nondegenerate.
If $\gamma^{m}$ is nondegenerate we say that
\begin{itemize}
\item $\gamma^{m}$ is {\em hyperbolic} if $d\psi_{mT}|_{\xh_{\gamma(0)}}$ has real eigenvalues, and that
\item  $\gamma^{m}$ is {\em elliptic} if $d\psi_{mT}|_{\xh_{\gamma(0)}}$ has complex eigenvalues.
\end{itemize}
We will furthermore say that
\begin{itemize}
\item $\gamma^{m}$ is {\em even} if it is hyperbolic and $d\psi_{mT}|_{\xh_{\gamma(0)}}$
 has positive eigenvalues, and
\item $\gamma^{m}$ is {\em odd} if it is either elliptic, or if
it is hyperbolic and $d\psi_{mT}|_{\xh_{\gamma(0)}}$ has negative eigenvalues.
\end{itemize}
The designation of a periodic orbit as even or odd will correspond to the parity of the Conley--Zehnder
index of that orbit (see \fullref{ss:cz}).

\subsection{Almost complex cobordisms}
Let $(M, \Ha)$ be a manifold equipped with a stable Hamiltonian structure
$\Ha=(\lambda, \omega)$.  We would like to define a preferred class of almost complex structures
on $\R\times M$ which interact in a specific way with the Hamiltonian structure.  First recall
that given any symplectic vector bundle $(E,\omega)$, a complex structure $J$ on $E$ is said to be
compatible with $\omega$ if the bilinear form defined by
\[
g_{J}(\cdot, \cdot)=\omega(\cdot, J\cdot)
\]
is a metric on $E$.  It is a well know fact that the space of all such $J$ is nonempty and contractible
in the $C^{\infty}$ topology (see e.g.\ Proposition 5 and discussion following in Section 1.3 of \cite{hz}).

Recalling now that $(\xh,\omega)$ is a
symplectic vector bundle, we define the set $\J(M, \Ha)$ to be the set of complex structures on $\xh$
which are compatible with $\omega$.  Given a $J\in\J(M,\Ha)$, we can extend it to an
$\R$-invariant almost complex structure $\tildeJ$ on $\R\times M$ by requiring
\[
\tildeJ\partial_{a}=\X \hskip.25in\text{ and }\hskip.25in \tildeJ|_{\xh}=J
\]
where $a$ is the parameter along $\R$.  We will refer to the almost complex structure $\tildeJ$ on
$\R\times M$ defined in this way as the \emph{standard cylindrical almost complex structure}
associated to $J\in\J(M, \Ha)$.

Let $W$ be a $4$-manifold without boundary.
We define a \emph{positive Hamiltonian structured} end to be data
$\wt E^{+}=(E^{+}, \Phi^{+}, M^{+}, \Ha^{+})$
where  $E^{+}\subset W$ is an open subset of $W$,
$M^{+}$ is a closed (possible disconnected) manifold equipped with a stable Hamiltonian structure
$\Ha^{+}$,
and $\Phi^{+}:E^{+}\to\R^{+}\times M^{+}$ is a diffeomorphism.
Similarly we will define a negative Hamiltonian structured end to be data
$\wt E^{-}=(E^{-}, \Phi^{-}, M^{-}, \Ha^{-})$ where everything is as before except that
$\Phi^{-}$ is now a diffeomorphism mapping $E^{-}$ to $\R^{-}\times M^{-}$.
A $4$-manifold $W$ equipped with (possibly empty) cylindrical ends
$\wt E^{\pm}=(E^{\pm}, \Phi^{\pm}, M^{\pm}, \Ha^{\pm})$
will be called a \emph{cobordism of Hamiltonian structures}
or a \emph{manifold with Hamiltonian structured ends} if
$W\setminus E^{+}\cup E^{-}$ is a compact manifold with (possibly empty) boundary.
If a manifold $W$ with Hamiltonian structured ends
$\wt E^{\pm}=(E^{\pm}, \Phi^{\pm}, M^{\pm}, \Ha^{\pm})$ is equipped with an almost complex structure
$\barJ$, we say that $\barJ$ is compatible with the ends $\wt E^{\pm}$ if
it is conjugated to a standard cylindrical almost complex structure on the ends,
that is, if $d\Phi^{\pm}\circ\barJ=\tildeJ^{\pm}\circ d\Phi^{+}$
for some $J^{\pm}\in\J(M^{\pm}, \Ha^{\pm})$.
We will refer to $(\R\times M, \barJ)$ as a \emph{cylindrical cobordism} if
$\barJ=\tildeJ$ is the standard cylindrical
almost complex structure associated to $J\in\J(M,\Ha)$ for some Hamiltonian structure $\Ha$ on $M$.
In practice, we will usually suppress the set $E^{\pm}$ and diffeomorphism $\Phi^{\pm}$
and just refer to the data $(\R^{\pm}\times M^{\pm}, \Ha^{\pm})$ as a cylindrical end,
or $(\R^{\pm}\times M^{\pm}, \Ha^{\pm}, J^{\pm})$ when we wish to specify the almost complex structure on
the ends.

Now, let $(W_{1}, \barJ_{1})$ and $(W_{2}, \barJ_{2})$ be
cobordisms  of Hamiltonians structures equipped with compatible almost complex structures,
and assume that $W_{1}$ is equipped with negative cylindrical end
$(\R^{-}\times M_{1}, \Ha_{1}, J_{1})$
and that $W_{2}$ is equipped with positive cylindrical end
$(\R^{+}\times M_{2}, \Ha_{2}, J_{2})$.  We say that
$(W_{1}, \barJ_{1})$ can be stacked on $(W_{2}, \barJ_{2})$ if
$M_{1}=M_{2}=M$, $\Ha_{1}=\Ha_{2}$, and $J_{1}=J_{2}$.
In this case, we define the concatenation $W_{1}\odot W_{2}$
of $W_{1}$ and $W_{2}$ to be the
$C^{0}$-manifold obtained by compactifying the negative end
of $W_{1}$ with $\br{-\infty}\times M$, compactifying the positive end of
$W_{2}$ with $\br{+\infty}\times M$
and making the obvious identification of
$\br{-\infty}\times M$ with $\br{+\infty}\times M$.
This operation can obviously be generalized to
an arbitrarily long list of cobordisms
$\br{(W_{1}, \barJ_{1}), (W_{2}, \barJ_{2}), \dots, (W_{N}, \barJ_{N})}$
provided of course that the negative end of $W_{i}$ matches with the positive end of
$W_{i+1}$.

We remark that in the setting of symplectic field theory,
the cobordisms considered are usually equipped
with a symplectic form, and
it is assumed that the almost complex structure is compatible with the symplectic
form.
In some cases, these conditions allow one to obtain topological controls on energy, and this control
in turn is important for compactness theorems necessary to define the algebraic structure of the
theory.  However, here we are primarily concerned with intersection-theory-related algebraic
invariants which only require an almost complex structure with the correct asymptotic behavior.
Therefore, we will not make any assumptions about the behavior of the almost complex structure away from
the cylindrical ends.

\subsection{Asymptotically cylindrical maps and pseudoholomorphic curves}
Let $(M, \Ha=(\lambda,\omega))$ be a closed manifold equipped with a stable Hamiltonian structure,
and let $\gamma\in\PMH$ be a simple $\tau$-periodic orbit of $\X$.  
For some $m\in\Z\setminus\br{0}$, consider the map
\[
\tl \gamma^m:\R\times S^1\to\R\times M
\]
defined by
\[
\tl\gamma^m(s,t)=(m\tau s, \gamma^m(t))\in\R\times M.
\]
We will refer to such a map as a \emph{cylinder over the periodic orbit $\gamma^m$},
or simply as an \emph{orbit cylinder}.
The reader should note that for any $J\in\J\MH$, orbit cylinders have $\tildeJ$-invariant tangent spaces.

In this paper, the main objects we will study are maps from a punctured Riemann surface to
a cobordism of Hamiltonian structures
with the maps asymptotic at the punctures to orbit cylinders.
More precisely, we consider a quadruple
$(\Sigma, j, \Gamma, \tl u)$ where
$(\Sigma, j)$ is a closed Riemann surface, 
$\Gamma\subset\Sigma$ is a finite set,
and
\[
\tl u=(a,u):\Sigma\setminus\Gamma\to\R\times M
\]
is a smooth map.
We say that \emph{$\tl u$ is asymptotically cylindrical over $\gamma^m$ at $z_*\in\Gamma$}
if there exists a holomorphic embedding
\[
\phi:[R,\infty)\times S^{1}\subset\C/i\Z \to\Sigma\setminus\Gamma
\]
satisfying $\lim_{s\to\infty}\phi(s,t)=z_*$ so that the maps
\[
\tl v_c:[R,\infty)\times S^1\to\R\times M
\]
defined by
\[
\tl v_c:=\bp{a(\phi(s+c,t))-m\tau c, u(\phi(s+c,t))}
\]
satisfy
\[
\lim_{c\to\infty}\tl v_c=\tl\gamma^m|_{[R,\infty)\times S^1} \text{ in } C^1([R,\infty)\times S^1, \R\times M).
\]
The map $\tl u$ is said to be a \emph{smooth asymptotically cylindrical map}
if $\tl u$ is asymptotically cylindrical at each $z\in\Gamma$ over some periodic orbit
$\gamma_z^{m_z}$.
In this case we will say that $\gamma_z^{m_z}$ is the asymptotic limit
of $\tl u$ at $z$.
Note that if $m_z>0$ (resp.\ $<0$) then the $\R$-component of $\tl u$ approaches
$+\infty$ (resp. $-\infty$) near $z$.
If $m_z>0$, we will refer to $z$ as a \emph{positive puncture} of $\tl u$, and similarly
if $m_z<0$ we will refer to $z$ as a \emph{negative puncture} of $\tl u$.
When  convenient we will write
$\Gamma=\Gamma_+\cup\Gamma_-$
to indicate how $\Gamma$ decomposes as positive punctures, $\Gamma_+$ and
negative punctures $\Gamma_-$.

We will denote the space of smooth asymptotically cylindrical maps
in $\R\times M$
from a genus $g$ surface with $n$ punctures by
\[
C^\infty_{g, n}\MH,
\]
and we will let
\[
C^\infty\MH:=\cup_{g>0}\cup_{n\ge 0}C^\infty_{g, n}\MH
\]
denote the space of all smooth asymptotically cylindrical maps in $\R\times M$.

In a similar manner, we can define asymptotically cylindrical maps in a
$4$-manifold $W$ with Hamiltonian structured ends $E^\pm=(\R^{\pm}\times M^\pm, \Ha^\pm)$.
Let $\tl u:\Sigma\setminus\Gamma\to W$ be a smooth map, and
assume that each $z\in\Gamma$ has an open neighborhood $U_z\subset\Sigma$
so that the image $\tl u(U_z\setminus\br{z})$ of the punctured neighborhood lies entirely
within one of the cylindrical ends $E^+$ or $E^-$.  Then we can think of
the map
$\tl u|_{U_z\setminus\br{z}}$
as a map to $\R\times M^\pm$, and define what it means for $\tl u$ to be asymptotically
cylindrical over a periodic orbit as we did above.  The map $\tl u$ is then said to be asymptotically
cylindrical if it is asymptotically cylindrical at each $z\in\Gamma$ over some
periodic orbit $\gamma_z^{m_z}$.
For asymptotically cylindrical maps in a cobordism, the punctures approaching periodic
orbits $\gamma_z^{m_z}$ in the positive cylindrical end $\R^+\times M^+$ are always positive punctures
(i.e. have $m_z>0$), and similarly those punctures at which $\tl u$ approaches a periodic orbit
$\gamma_z^{m_z}$
in the negative end $\R^-\times M^-$ have $m_z<0$.
We will denote the space of genus-$g$, $n$-punctured asymptotically cylindrical
maps in $W$ by
\[
C^\infty_{g,n}(W, \Ha^+, \Ha^-)
\]
and the space of all asymptotically cylindrical maps in $W$ by
\[
C^\infty(W, \Ha^+, \Ha^-).
\]

Now, let $(W_{1}, \barJ_{1})$ and $(W_{2}, \barJ_{2})$ be almost complex cobordisms, and
assume that $W_{1}$ is equipped with negative cylindrical end
$(\R^{-}\times M, \Ha, J)$
and that $W_{2}$ is equipped with positive cylindrical end
$(\R^{+}\times M, \Ha, J)$ so that we can form the concatenation $W_1\odot W_2$.
For $i\in\br{1, 2}$, let
\[
\tl u_i:\Sigma_i\setminus(\Gamma_{i,+}\cup\Gamma_{i,-})\to W_i
\]
be asymptotically cylindrical maps.
Assume that there exists a bijection
\[
i:\Gamma_{2,+}\to\Gamma_{1,-}
\]
so that if $\tl u_2$ as asymptotic at $z\in\Gamma_{2,+}$ to
$\gamma^{m}$, then $\tl u_1$ is asymptotic at $i(z)\in\Gamma_{1,-}$ to
$\gamma^{-m}$.
Then we can form the concatenated map
\[
\tl u_1\odot \tl u_2:\Sigma_1\odot\Sigma_2\setminus(\Gamma_1^+\cup\Gamma_2^-)
\to W_1\odot W_2
\]
where $\Sigma_1\odot\Sigma_2$ is the topological surface formed by circle compactifying
$\Sigma_1$ at its negative punctures, circle compactifying $\Sigma_2$ at its positive punctures, and then
identifying 
each negative circle with its corresponding (under the bijection $i$) positive circle in a way that
makes $\tl u_1\odot\tl u_2$ a continuous map.  Note that when the asymptotic data
contains multiply covered orbits, 
$\Sigma_1\odot\Sigma_2$ is only well-defined up to
Dehn twists unless further choices (namely so-called asymptotic markers)
are made.  The specifics won't be important here, so we won't address this issue any further.
We will refer to a map $\tl u_{1}\odot \tl u_{2}$ constructed in this way as a
\emph{smooth asymptotically cylindrical building} in $W_{1}\odot W_{2}$.

Consider an asymptotically cylindrical map
$(\Sigma, j, \Gamma, \tl u)\in C^{\infty}(W, \Ha^{+}, \Ha^{-})$
and assume that $W$ is equipped with an almost complex structure $\barJ$
compatible with the cylindrical ends.
If the map $\tl u$ satisfies the equation 
\begin{equation}\label{e:j-hol}
d\tl u\circ j=\barJ\circ d\tl u
\end{equation}
we say that $(\Sigma, j, \Gamma, \tl u)$ is an \emph{asymptotically cylindrical pseudoholomorphic map}.
We define an equivalence relation on punctured pseudoholomorphic maps by
saying that $(\Sigma, j, \Gamma, \tl u)$ is equivalent to $(\Sigma', j', \Gamma', \tl u')$
if there exists a biholomorphic map
\[
\phi:\Sigma'\to\Sigma
\]
so that $\Gamma=\phi(\Gamma')$ and $\tl u\circ\phi=\tl u'$.  An equivalence class of maps
$[\Sigma, j, \Gamma, u]$ will be referred to as an \emph{asymptotically cylindrical pseudoholomorphic curve}.
We will use the notations
\[
\M(W, \barJ, \Ha^{+}, \Ha^{-}) \text{ and } \M(M, \Ha, J)
\]
to denote the set of asymptotically cylindrical pseudoholomorphic curves in $(W, \barJ)$
or $(\R\times M, \tildeJ)$ respectively,
and similarly we will use
\[
\M_{g, n}(W, \barJ, \Ha^{+}, \Ha^{-}) \text{ and } \M_{g, n}(M, \Ha, J).
\]
if we wish to specify the genus and number of punctures.

N.B.\ To simplify our language, we will henceforth simply use the term
``pseudoholomorphic curve/map'' when we are referring to asymptotically cylindrical
pseudoholomorphic curves/maps.
{\bf
The reader should always assume that all pseudoholomorphic curves are asymptotically cylindrical unless stated otherwise.
}

As with smooth asymptotically cylindrical maps, we can concatenate pairs (or finite lists) of asymptotically
cylindrical pseudoholomorphic maps into cobordisms with matching ends,
providing the asymptotic data match appropriately.
Following \cite{behwz}, we will call such maps \emph{pseudoholomorphic buildings}.

\subsection{Main results}\label{ss:intro-main}
Throughout this section, we will let $(M, \Ha=(\lambda, \omega))$
denote a $3$-manifold equipped with a \emph{nondegenerate} stable Hamiltonian structure,
and we will let $(W, \barJ)$ denote an almost complex $4$-manifold equipped with cylindrical
ends $E^\pm=(\R^\pm\times M^\pm, \Ha^\pm, J^\pm)$ where the
stable Hamiltonian structures $\Ha^\pm$ are assumed to be nondegenerate.

The main results of this paper are concerned with the intersection properties of
pseudoholomorphic curves in $4$-manifolds with Hamiltonian structured cylindrical ends.
Due to the noncompactness of the domains of the curves we consider,
the problem of understanding intersection behavior
becomes particularly subtle in the case that
the curves have
multiple ends which
approach coverings of the same orbit, or --
when considering self-intersection problems -- in the case that an end approaches a multiple cover
of an orbit.
Indeed, while all intersections (or self-intersections and singularities) must be 
isolated by the results of Micallef and White in \cite{mw}, it is not immediately obvious
why the algebraic intersection number of two curves -- computed by summing local intersection indices -- 
must be finite since it is conceivable that a pair of curves could have a sequence of
intersections approaching the punctures.
Similarly, it is not clear that a single curve couldn't have a sequence of self-intersections
approaching a puncture at which the curve is asymptotic to a multiple cover of an orbit
or 
approaching a pair of punctures at which the curve is asymptotic to coverings of the same orbit.
While we will see that the asymptotic results of \cite{sief08} imply
that the
intersection number of two curves or the self-intersection index of a single curve must in fact be finite
(see \fullref{c:finite-isect} and \fullref{c:finite-sing} below),
an additional unavoidable complication is that these quantities
may not be stable under homotopies since
intersections could run in or out of the punctures at shared asymptotic
limits.
We thus seek to determine to what degree these quantities
are topologically controlled
by the homotopy classes of the maps in
$C^\infty(M, \Ha)$ or $C^\infty(W, \Ha^+, \Ha^-)$ (depending on the target manifold).

As a first step towards
finding intersection-related  topological invariants
of a pair of maps with common asymptotic limits, 
we perturb one of the maps in a prescribed direction near the
ends.
This idea was studied locally
for a single embedded pseudoholomorphic half-cylinder in
\cite{kr} and further pursued in \cite{hut02}, \cite{sief05}, and \cite{hut09}.
More precisely, let $\Phi$ denote a choice of trivialization of
the plane-field $\xh$ along every simple periodic orbit.
Then given two maps $\tl u$, $\tl v\in C^\infty(W, \barJ, \Ha^+, \Ha^-)$, we define
the \emph{relative intersection number} $i^\Phi(\tl u, \tl v)$ of $\tl u$ and $\tl v$ by
\[
i^\Phi(\tl u, \tl v)=\inum(\tl u,\tl v_\Phi)
\]
where $\tl v_\Phi$ is the map obtained by perturbing $\tl v$ near the punctures in a direction determined by
the trivialization $\Phi$, and where ``$\inum$'' denotes the algebraic intersection number computed by
perturbing $\tl u$ and $\tl v_\Phi$ 
on compact subsets of their domain so that they are transverse, and
then counting intersections with sign.
It is well known that the resulting count is independent of the compactly supported perturbation
of $\tl u$ and $\tl v_{\Phi}$ since $\tl u$ and $\tl v_{\Phi}$ are disjoint outside of a compact set.
Thus the relative intersection number computed in this way depends on the homotopy classes of
$\tl u$ and $\tl v$ in $C^\infty(W, \barJ, \Ha^+, \Ha^-)$, and the homotopy class of the trivialization $\Phi$.

While the relative intersection number gives a topological invariant of a pair of asymptotically
cylindrical maps,
its relationship to the algebraic intersection number of
two (unperturbed) pseudoholomorphic curves is not clear, and in particular it's
not immediately clear whether or how it can be used to find a bound on the algebraic intersection number
of two pseudoholomorphic curves.
We will see below that through careful consideration of
the relative asymptotic behavior of the curves from
\cite{sief08},
such a bound depending only on the homotopy
classes of the maps can be 
obtained
by adding an appropriate quantity to
the relative intersection number which balances the
trivialization dependence.  Motivated by the local analysis in \fullref{s:local}, we define
the \emph{generalized intersection number} $\gin{\tl u}{\tl v}$ of two asymptotically cylindrical maps
$(\Sigma, j, \Gamma,  \tl u)$, $(\Sigma', j', \Gamma',  \tl v)\in C^\infty(W, \Ha^+, \Ha^-)$ as follows:
assuming that $\tl u$ is asymptotic at $z\in\Gamma$ to a cylinder over $\gamma_{z}^{m_{z}}$
and similarly that $\tl v$ is asymptotic at $w\in\Gamma'$ to a cylinder over $\gamma_{w}^{m_{w}}$
we define $\gin{\tl u}{\tl v}$ by
\begin{equation}\label{e:gin-def-intro}
\gin{\tl u}{\tl v}=i^\Phi(\tl u, \tl v)
+\sum_{\substack{(z, w)\in\Gamma\times\Gamma' \\ \gamma_z=\gamma_w \\ m_z m_w>0}}
m_{z}m_{w}\max\br{\tfrac{\fl{\mu^\Phi(\gamma_z^{m_z})/2}}{|m_z|}, \tfrac{\fl{\mu^\Phi(\gamma_w^{m_w})/2}}{|m_w|}}
\end{equation}
where $\mu^\Phi(\gamma^{m_z})$ denotes the Conley--Zehnder index of the periodic orbit $\gamma_z^{m_z}$
(see \cite{hwz:prop2} and \fullref{s:local} below),
and where $\fl{\cdot}$ denotes the greatest integer function.  Note that the sum here is taken over
all pairs of punctures with the same sign
where the maps in question are asymptotically cylindrical over coverings of the same underlying
simple (unparametrized) periodic orbit.

The following theorem summarizes the main properties of the generalized intersection number.

\begin{theorem}[Properties of the generalized intersection number]\label{t:gin-prop}
Let $W$, $W_{1}$ and $W_{2}$ be $4$-manifolds with Hamiltonian structured
cylindrical ends, and assume we can form the concatenation $W_{1}\odot W_{2}$.  Then:
\begin{enumerate}
\item  If $(\Sigma, j, \Gamma, \tl u)$ and $(\Sigma', j', \Gamma', \tl v)\in\CW$
are asymptotically cylindrical maps
then the generalized intersection
number $\gin{\tl u}{\tl v}$ depends only on the homotopy classes of $\tl u$ and $\tl v$
in $\CW$.

\item For any $(\Sigma, j, \Gamma, \tl u)$ and $(\Sigma', j', \Gamma', \tl v)\in\CW$
\[
\gin{\tl u}{\tl v}=\gin{\tl v}{\tl u}.
\]

\item If $(\Sigma, j, \Gamma, \tl u)$, $(\Sigma', j', \Gamma', \tl v)$,
$(\Sigma'', j'', \Gamma'', \tl w)\in\CW$
then
\[
\gin{\tl u+\tl v}{\tl w}=\gin{\tl u}{\tl w}+\gin{\tl v}{\tl w}
\]
where ``$+$'' on the left hand side denotes the disjoint union of the maps $\tl u$ and $\tl v$.

\item
If $\tl u_{1}\odot \tl u_{2}$ and $\tl v_{1}\odot \tl v_{2}$ are asymptotically cylindrical buildings in
$W_{1}\odot W_{2}$ then
\[
\gin{\tl u_{1}\odot \tl u_{2}}{\tl v_{1}\odot \tl v_{2}}\ge \gin{\tl u_{1}}{\tl v_{1}}+\gin{\tl u_{2}}{\tl v_{2}}.
\]
Moreover, strict inequality occurs if and only if there is a periodic orbit $\gamma$ so
that $\tl u_{1}$ has a negative puncture asymptotic to $\gamma^m$, $\tl v_1$ has a negative puncture asymptotic
to $\gamma^n$, and both $\gamma^m$ and $\gamma^n$ are odd orbits.
\end{enumerate}
\end{theorem}

In item $(4)$ above, the possibility of strict inequality has to do with the fact that at an odd
orbit (with a fixed multiplicity) the eigenvectors of the asymptotic operator controlling
the direction of approach of negative pseudoholomorphic
ends must have strictly greater winding (computed relative
to the direction of the Reeb flow) than those controlling the direction of approach of
positive pseudoholomorphic ends.  The relevant details from \cite{hwz:prop1}, \cite{hwz:prop2} are reviewed in
\fullref{ss:asymptotics} below.  In some applications it is convenient to
modify the generalized intersection number
to include information about shared odd orbits so that
one has a product which is level-wise additive, i.e.\ always satisfies equality in item $(4)$ above.
This approach is taken by Momin in \cite{mom10-2} to study
a variation on contact homology where only curves contained in
the complement of a prescribed collection of elliptic orbits are included in the differential.

One of the motivations for defining the generalized intersection number as we did is the next theorem,
which generalizes the fact that two closed curves without common components have a nonnegative homological
intersection number, and that the intersection number vanishes only if the two curves do not intersect.
The
\emph{total asymptotic intersection index}
$\delta_\infty(\tl u, \tl v)$ mentioned in the theorem
is defined in \fullref{s:local} below.  For the moment
the reader should know that it is a nonnegative quantity
defined for pairs of pseudoholomorphic curves having no common components
that can be thought of as a
measure of the degree of tangency at infinity between the two curves.
Here, we say the
two maps $\tl u$ and $\tl v$
have no common components if there is no component of the domain of $\tl u$ which has
image identical to that of a component of the domain of $\tl v$.

\begin{theorem}\label{t:gin-positivity-intro}
Let $(W^{4}, \barJ)$ be an almost complex cobordism with cylindrical ends
$(\R^\pm\times M^\pm, J^\pm, \Ha^\pm)$, and let
$[\Sigma, j, \Gamma, \tl u]$, $[\Sigma', j', \Gamma', \tl v]\in\M(W, \barJ, \Ha^+, \Ha^-)$ be pseudoholomorphic curves in $W$ with no common components.
Then
\begin{equation}\label{e:positivity}
 \gin{\tl u}{\tl v}=\inum(\tl u, \tl v)+\delta_\infty(\tl u, \tl v),
\end{equation}
where $\inum(\tl u, \tl v)$ is the algebraic intersection number of $\tl u$ and $\tl v$, and
$\delta_{\infty}(\tl u, \tl v)$ is the asymptotic intersection index of $\tl u$ and $\tl v$.
In particular
\[
\gin{\tl u}{\tl v}\ge\inum(\tl u, \tl v)\ge 0,
\]
and
\[
\gin{\tl u}{\tl v}=0
\]
if and only if $\tl u$ and $\tl v$ don't intersect, and the total asymptotic intersection index vanishes,
i.e.\
$\delta_\infty(\tl u, \tl v)=0$.
\end{theorem}

As an immediate corollary to this, we can conclude that if a homotopy class 
of maps with connected domains in
$C^\infty(W, \Ha^+, \Ha^-)$ has a negative generalized self-intersection number, then that homotopy class
can contain at most one pseudoholomorphic curve for a given compatible
almost complex structure.

We note that in contrast to the case of closed pseudoholomorphic curves,
positivity of the generalized intersection number of two given curves does
not guarantee that the two curves intersect.  This is because the algebraic intersection number
of two punctured curves is not
a homotopy-invariant quantity, and intersections between two pseudoholomorphic curves can disappear
out the punctures.  The degree to which intersections can escape or appear at the punctures
is however topologically controlled: even though both quantities on the right hand side
of \eqref{e:positivity} can vary under homotopies of curves, their sum is topologically determined,
and this fact coupled with the nonnegativity of those quantities allows a topological bound
on the total count of intersections and ``tangencies at infinity.'' 
One might consider this fungibility of intersections and ``tangencies at infinity''
as motivation to think of tangencies at infinity, or
``asymptotic intersections,'' as being somehow equivalent to actual intersections.
This viewpoint is further supported by the fact that
considering weighted Fredholm theory arguments as in \cite{hwz:prop3, dr},
one would expect that the space of pairs of curves $(\tl u, \tl v)$
with $\delta_\infty(\tl u,\tl v)>0$ should have positive codimension in the universal moduli space of
pseudoholomorphic curves.

We next state a generalization for punctured curves of the adjunction formula
\eqref{e:adj-closed}.  We first establish some notation and terminology.
Let $(\Sigma, j, \Gamma, \tl u)\in C^\infty(W, \Ha^+, \Ha^-)$ be an asymptotically cylindrical map,
and assume that at $z\in\Gamma$, $u$ is asymptotic to the periodic orbit
$\gamma_z^{m_z}$.  A choice of complex trivialization of the hamiltonian plane-field distribution
$(\xi^{\Ha^\pm}, J^\pm)$ along each $\gamma_z$ induces a complex trivialization of
the pull-back bundle $(\tl u^*TW, \barJ)$ since on the cylindrical ends we have
\[
(TW, \barJ)
\approx (\R\oplus\R X_{\Ha^{\pm}}\oplus\xi^{\Ha^\pm}, \tildeJ^{\pm})
\approx\C X_{\Ha^{\pm}}\oplus(\xi^{\Ha^\pm}, J^{\pm}).
\]
Given such a choice of trivialization, $\Phi$, we can define the relative Chern number
$c_1^\Phi(\tl u^*TW)$ which is the obstruction to extending over $\Sigma\setminus\Gamma$
the trivialization of
$(\tl u^*TW, \barJ)\to\Sigma\setminus\Gamma$
that has been chosen at the ends
(see \cite{hut02,sief05} or \fullref{ss:chern} below for a precise definition).
We define the \emph{total Conley--Zehnder index} $\mu(\tl u)$ of the map $\tl u$ by
\[
\mu(\tl u):=2c_1^\Phi(\tl u^*TW)+\sum_{z\in\Gamma}\mu^\Phi(\gamma_z^{m_z}).
\]
It follows from change of trivialization formulas given below
for the relative first Chern number and the
Conley--Zehnder index that the total Conley--Zehnder index does not depend on the choice
of trivialization.

Next, assuming $\gamma$ is a simple periodic orbit of a Reeb vector field
we define the \emph{spectral covering number} of the periodic orbit $\gamma^m$ by
\[
\cover(\gamma^m)=\gcd(m, \fl{\mu^\Phi(\gamma^m)/2}).
\]
We observe that this quantity does not depend on a choice of trivialization since
changing the trivialization changes $\fl{\mu^\Phi(\gamma^m)/2}$ by an integer multiple of $m$
(see comments following \fullref{l:alpha-triv}).
With $u$ as in the previous paragraph, we then define the
\emph{total spectral covering number} $\cover(u)$ of $u$ by
\[
\cover(u)=\sum_{z\in\Gamma}\cover(\gamma_z^{m_z}).
\]
We observe that the total spectral covering number of a map depends only
on its asymptotic limits, and not on the map itself.

We now state a generalization of \eqref{e:adj-closed} involving the generalized
intersection number.  The quantity $\delta(\tl u)$ in equation \eqref{e:adj-main-intro},
which we call the self-intersection index,
is a nonnegative integer-valued quantity which records information about
double points and singular points, and is zero if and only if $u$ is an embedding
(see \cite{mw} and \fullref{sss:adjunction} below).
The \emph{asymptotic self-intersection index} $\delta_\infty(\tl u)$ of $\tl u$ is defined in
\fullref{ss:asympt-int} below.  For the moment, the reader should know that it is a nonnegative,
integer-valued quantity, defined for a simple pseudoholomorphic curve that
can be thought of as a measure of the degree of self-tangency at infinity.
Here, we say a pseudoholomorphic curve is simple if it does not factor
through a branched cover.

\begin{theorem}\label{t:adj-main-intro}
Let $[\Sigma, j, \Gamma, \tl u]\in\M(W, \barJ, \Ha^+, \Ha^-)$
be a connected pseudoholomorphic curve,
and assume that $\tl u$ is simple.  Then
\begin{equation}\label{e:adj-main-intro}
\gin{\tl u}{\tl u}-\tfrac{1}{2}\mu(\tl u)+\tfrac{1}{2}\#\Gamma_{\text{odd}}+\chi(\Sigma)-\cover(\tl u)=2[\delta(\tl u)+\delta_\infty(\tl u)],
\end{equation}
where $\#\Gamma_{\text{odd}}$ is the number of punctures of $\tl u$ with odd Conley--Zehnder indices.
In particular, 
\[
\gin{\tl u}{\tl u}-\tfrac{1}{2}\mu(\tl u)+\tfrac{1}{2}\#\Gamma_{\text{odd}}+\chi(\Sigma)-\cover(\tl u)\ge 0
\]
and equality occurs if and only if $\tl u$ is an embedding and the asymptotic self-intersection index
of $\tl u$ vanishes, i.e.\ $\delta_\infty(\tl u)=0$.
\end{theorem}

Note that in contrast to the case of closed pseudoholomorphic curves and \eqref{e:adj-closed},
positivity of the left-hand side of equation \eqref{e:adj-main-intro}
does not guarantee a curve is not embedded.  This is because the self-intersection index
$\delta(\tl u)$ is not,
in general, topologically determined for punctured curves
since double points can escape out the ends in families of curves.
However, much like in the discussion following \fullref{t:gin-positivity-intro} above, this theorem
shows that the \emph{sum} of the 
self-intersection index and the asymptotic self-intersection index
\emph{is} topologically determined, and that double points escaping out the ends are traded for
self-tangencies at infinity.  Thus,
we might think of having a positive asymptotic self-intersection index as being somehow equivalent to having
double points or singular points,
a point of view which is again supported by the fact that the Fredholm theory
suggests that curves with $\delta_\infty(\tl u)>0$ should form strata of positive codimension in
the universal moduli space.

Results closely related to
\fullref{t:adj-main-intro} have been proved by
Hutchings in \cite{hut02} and \cite{hut09}.
Due to the specific application being pursued in those papers,
no theorem equivalent to \fullref{t:adj-main-intro} appears.
Equation (18) in Remark 3.2 of \cite{hut02} identifies embedded curves, but not in terms
of quantities that are topologically determined.
Still, the appropriate analogue of that result to the present context is the key step in proving
\fullref{t:adj-main-intro}.
On the other hand, the index inequality
in Theorem 1.7 of \cite{hut02} gives a purely topological condition that guarantees
a curve will be embedded,
but only identifies those curves
whose asymptotic data satisfy an additional ``admissibility'' condition.
\fullref{t:adj-main-intro} above, can thus be thought of as being intermediate between these two results.

\subsubsection{Cylindrical cobordisms}
In this section we specialize to the case of a cylindrical cobordism $(\R\times M, \tildeJ)$ equipped with
an $\R$-invariant almost complex structure $\tildeJ$ arising from a compatible $J\in\J\MH$ for
some nondegenerate stable Hamiltonian structure $\Ha=(\lambda, \omega)$ on $M$.

An important (and well-known) observation is that space of curves $\M(M, \Ha, J)$
is equipped with an $\R$-action, defined by translating the $\R$-coordinate of the given curve.
More precisely, if $[\Sigma, j, \Gamma, \tl u]$ is a pseudoholomorphic curve, and
we write $\tl u=(a,u)\in\R\times M$, then
$[\Sigma, j, \Gamma, \tl u_c]$ where $\tl u_c(z):=(a(z)+c, u(z))$ is also a pseudoholomorphic curve.
Furthermore, in the cylindrical case, 
the way the almost complex structure is defined allows
the pseudoholomorphic curve equation \eqref{e:j-hol} 
for the map $(\Sigma, j, \Gamma, (a, u))$ to be rewritten,
\begin{equation}\label{e:j-hol-cylindrical}
\begin{gathered}
u^{*}\lambda\circ j=da \\
\pi_{\xh}\circ du\circ j=J\circ\pi_{\xh}\circ du
\end{gathered}
\end{equation}
where $\pi_{\xh}:TM=\R\X\oplus\xh\to\xh$ is the projection of $TM$ onto $\xh$
determined by the splitting \eqref{e:splitting}.
Thus, the $M$-component, $u$, of a curve in $\R\times M$ determines the $\R$-component, $a$, up
to a constant.

As an immediate consequence of the existence of the $\R$-action,
\fullref{t:gin-positivity-intro}, and the homotopy invariance of the generalized intersection number,
we can conclude that if a connected curve $[\Sigma, j, \Gamma, \tl u]$ has no component whose image lies
in an orbit cylinder, then
\[
\gin{\tl u}{\tl u}=\gin{\tl u}{\tl u_{c}}\ge 0
\]
so such a curve must have a nonnegative self-intersection number.

Since, in the cylindrical case, the projection of the curve to the $3$-manifold contains
all of the information (up to an $\R$-shift), it is not surprising that generalized intersection number of two
curves in $\R\times M$ can be computed in terms of $3$-dimensional invariants associated to the curves;
namely we can state a formula in terms of
the intersections between one of the curves with periodic orbits which are asymptotic limits
of the other curve, and how the ends of the curves wind around the periodic orbits
(see \fullref{t:gin-cylindrical} below).  Moreover,
since the formula we give decomposes into a collection of nonnegative terms,
this computation then allows one to deduce necessary and
sufficient conditions for the generalized intersection number to vanish, which we state in
\fullref{c:gin-zero}.  

Part of the reason we seek to find conditions to characterize when the generalized
intersection number of two curves
$[\Sigma, j, \Gamma, \tl u=(a, u)]$ and $[\Sigma', j', \Gamma', \tl v=(b, v)]\in\M(M, \Ha, J)$
is zero (or nonzero) is that $\gin{\tl u}{\tl v}=0$ implies that the projections $u$, $v$, of the maps to $M$
don't intersect, provided that these maps have no components projecting to identical images in $M$.
Indeed, if we can find a $z\in\Sigma$ and $w\in\Sigma'$ so that
$u(z)=v(w)$ then we can find a $c$ so that $\tl u(z)=\tl v_c(w)$, which in turn implies
that $\inum(\tl u, \tl v_c)>0$ and hence $\gin{\tl u}{\tl v}=\gin{\tl u}{\tl v_c}\ge \inum(\tl u, \tl v_{c})>0$
by \fullref{t:gin-positivity-intro}.
Thus the vanishing of the generalized intersection number gives a sufficient condition for the projection
of the curves to the three-manifold to not intersect.  
The vanishing of the generalized intersection number is not a necessary condition for
the projected curves to not intersect since
it is possible for the algebraic intersection number
of $\tl u$ and $\tl v_c$ to be zero for all $c\in\R$, but still have $\gin{\tl u}{\tl v}>0$,
since the asymptotic intersection index could be nonzero.
However, in the $\R$-invariant setting, the asymptotic intersection number changes in a predictable manner
when $\R$-shifting one of the two curves, and we are able to establish a set
of necessary and sufficient conditions for the two projected curves to not intersect.

Before stating the relevant theorem we discuss some
of the more immediate
necessary conditions for the projected curves
$u$ and $v$ to not intersect.
We first recall some facts from \cite{hwz:prop1} and \cite{hwz:prop2}.
Namely, if a curve $[\Sigma, j, \Gamma, \tl u=(a,u)]$ does not have any components with image lying
in an orbit cylinder, then there is a neighborhood of each puncture on which the map does
not intersect any of its asymptotic limits.  Thus choosing a loop in one of these neighborhoods which
winds around the puncture once in a clockwise 
	direction,\footnote{This convention agrees with that used in \cite{hwz:prop2} for positive punctures,
	but is opposite that used in \cite{hwz:prop2} for negative punctures.
	What this means geometrically is that we compute windings by traversing the orbit 
	in a direction determined
	by the orientation of the boundary of the $S^1$-compactified punctured surface, while in
	\cite{hwz:prop2} the convention is to compute windings
	by traversing the orbit in the direction determined by the
	Reeb vector field.
	We use this convention
	because it simplifies the statements of most of our results.}
and choosing a trivialization $\Phi$ of $\xh$
along the orbit, we get a well defined winding number $\winfty^\Phi(\tl u; z)$ for each $z\in\Gamma$.
Moreover, it is shown in \cite{hwz:prop2} (and reviewed in \fullref{l:wind-infinity} below) that
\[
\winfty^\Phi(\tl u; z)\le \fl{\mu^\Phi(\gamma_z^{m_z})/2}
\]
where we assume $\tl u$ to be asymptotic to $\gamma_z^{m_z}$ at $z\in\Gamma$.
If we consider a small torus $T^{2}_{\gamma_{z}}$
bounding a tubular neighborhood of 
$\gamma_{z}$
and use the framing $\Phi$ to identify $H_{1}(T^{2}_{\gamma_{z}})$ with $\Z\oplus\Z$, the intersection 
of the projected map $u$ with $T^{2}_{\gamma_{z}}$ will be a curve in the homology
class $(m_{z}, \wind^{\Phi}(\tl u; z))$.  If $\tl v$ also has an end approaching $\gamma_{z}$ at
$w\in\Gamma'$, the homology class of
intersection of $v$ with $T^{2}_{\gamma_{z}}$ will then be given
by $(m_{w}, \wind^{\Phi}(\tl v; w))$.  Thus, if the projected maps $u$ and $v$ are to be disjoint,
we must have that the intersection number
\[
(m_{z}, \wind^{\Phi}(\tl u; z))\cdot (m_{w}, \wind^{\Phi}(\tl v; w))=m_{z}\wind^{\Phi}(\tl v; w)-m_{w}\wind^{\Phi}(\tl u; z)
\]
is zero or equivalently
\begin{equation}\label{e:torus-inum}
\tfrac{\winfty^{\Phi}(\tl u; z)}{m_z} = \tfrac{\winfty^{\Phi}(\tl v; w)}{m_w}.
\end{equation}
Therefore, if the projections $u$ and $v$ don't intersect, it is necessary that
\eqref{e:torus-inum} holds at any pair of punctures $(z, w)\in\Gamma\times \Gamma'$ at which 
$u$ and $v$ are asymptotic to coverings of a common orbit.

Next suppose that $u$ intersects one of the asymptotic limits of $v$.
Then one can use the asymptotic results of
\cite{hwz:prop1} (or see \fullref{t:asymp} below)
with the fact that the orbit is a projection to $M$ of a pseudoholomorphic curve in $\R\times M$
to argue that
$u$ must intersect $v$.
Thus if $u$ and $v$ don't intersect, it must be the case that $u$ intersects none of the asymptotic limits
of $v$, and vice versa.

As the following theorem shows, the necessary conditions we've just stated 
for the projections of two curves to
$M$ to not intersect are also sufficient, and in fact somewhat weaker conditions are sufficient.

\begin{theorem}\label{t:no-isect-intro}
Let $[\Sigma, j, \Gamma, \tl u=(a, u)]$ and $[\Sigma', j', \Gamma', \tl v=(b, v)]\in\M(M, \Ha, J)$
be pseudoholomorphic curves,
and assume that no component of $\tl u$ or $\tl v$ lies in in orbit cylinder,
and that the projected curves $u$ and $v$ do not have identical image on any component of their domains.
Then the following are equivalent:
\begin{enumerate}
\item The projected curves $u$ and $v$ do not intersect.

\item All of the following hold:
     \begin{enumerate}
     \item The map $u$ does not intersect any of the positive asymptotic limits of $v$.
     \item The map $v$ does not intersect any of the negative asymptotic limits of $u$.
     \item If $\gamma$ is a periodic orbit so that at $z\in\Gamma$, $\tl u$ is asymptotic
     to $\gamma^{m_z}$ and at $w\in\Gamma'$, $\tl v$ is asymptotic to $\gamma^{m_w}$, then:
          \begin{enumerate}
          \item If $m_z$ and $m_w$  have the same sign then
          \[
          \tfrac{\winfty(\tl u; z)}{m_z}\ge\tfrac{\winfty(\tl v; w)}{m_w}.
          \]
          \item If $m_z<0$ and $m_w>0$ then
          \[
          \tfrac{\winfty^\Phi(\tl u; z)}{m_z}=\tfrac{\fl{\mu^\Phi(\gamma^{m_z})/2}}{m_z}
          =\tfrac{\fl{\mu^\Phi(\gamma^{m_w})/2}}{m_w}=\tfrac{\winfty^\Phi(\tl v; w)}{m_w}
          \]
          (this is only possible if $\gamma^{m_z}$ and $\gamma^{m_w}$ are both even orbits).
    \end{enumerate}
     \end{enumerate}

\item All of the following hold:
     \begin{enumerate}
     \item The map $u$ does not intersect any of the asymptotic limits of $v$.
     \item The map $v$ does not intersect any of the asymptotic limits of $u$.
     \item If $\gamma$ is a periodic orbit so that at $z\in\Gamma$, $\tl u$ is asymptotic
     to $\gamma^{m_z}$ and at $w\in\Gamma'$, $\tl v$ is asymptotic to $\gamma^{m_w}$, then
          \[
          \tfrac{\winfty(\tl u; z)}{m_z} = \tfrac{\winfty(\tl v; w)}{m_w}.
          \]
     \end{enumerate}
\end{enumerate}
\end{theorem}

When two curves, $[\Sigma, j, \Gamma, \tl u]$ and $[\Sigma', j', \Gamma', \tl v]\in\MM$,
each have as an asymptotic limit a cover of the same even periodic orbit, it is sometimes possible to
conclude that the generalized intersection number $\gin{\tl u}{\tl v}$ is positive.
Let $\gamma$ denote either a simple, even orbit or the double cover of an odd, hyperbolic orbit.
In this case it is possible to use the asymptotic description from \fullref{t:asymp} below
to define a notion of two curves approaching
a cover of $\gamma$ \emph{in the same direction}.
We delay the precise definition to \fullref{ss:direction-of-approach}
because it is somewhat technical.  The following theorem then says that two curves
approaching a cover of $\gamma$ in the same direction is a sufficient condition
for a positive generalized intersection number.

\begin{theorem}\label{t:dir-of-approach-intro}
Let $\gamma$ be a periodic orbit satisfying the above assumptions, and
let $[\Sigma, j, \Gamma, \tl u=(a, u)]$ and $[\Sigma', j', \Gamma', \tl v=(b, v)]\in\MM$
be connected pseudoholomorphic curves.
Assume that at punctures $z\in\Gamma$ and $w\in\Gamma'$,
$\tl u$ and $\tl v$ approach a cover of $\gamma$ in the same direction, 
and that
there do not exist neighborhoods $U$ of $z$ and $V$ of $w$ so that
$u(U\setminus\br{z})=v(V\setminus\br{w})$.
Then
\[
\gin{\tl u}{\tl v}>0.
\]
\end{theorem}

Just as it is of use to know when the generalized intersection number of a pair
of curves vanishes, it is of use to know when the generalized self-intersection
number of a single curve $[\Sigma, j, \Gamma, \tl u=(a,u)]\in\M(M, \Ha, J)$ vanishes.
Combining  the adjunction formula \eqref{e:adj-main-intro} with
the previously referenced formula \eqref{e:gin-cylindrical}
for the generalized intersection number in terms of data associated to the projected curves,
we get a collection of conditions equivalent to the vanishing of the
generalized self-intersection number $\gin{\tl u}{\tl u}$, provided the map $\tl u$ is simple and
doesn't have any components contained in an orbit cylinder (\fullref{c:self-gin-zero} below).
From this result we will see that under the stated assumptions, $\gin{\tl u}{\tl u}=0$ implies that
$\tl u$ is an embedding which projects to an immersion in $M$ transverse to the flow
of $\X$ which doesn't intersect any of its asymptotic limits.

Combining this with the above discussion, we see that if $\tl u$ is furthermore connected, then
$\gin{\tl u}{\tl u}=0$ implies that $\tl u$ doesn't intersect any of its $\R$-translates, and thus that
the projected curve $u$ is an embedding (since it is an injective immersion that never intersects any of its asymptotic limits).
We can therefore conclude that if the projected map $u$ is not an embedding then
$\gin{\tl u}{\tl u}>0$.
As with the discussion above concerning intersections of projections of curves to $M$, 
the converse is not true, since it could be the
case that
the generalized intersection number
$\gin{\tl u}{\tl u_c}$
is positive,
but still that the intersection number $\inum(\tl u, \tl u_c)$ 
is zero for all $c\in\R\setminus\br{0}$.

However, again as above,
the asymptotic intersection index of $\tl u$ and $\tl u_c$ changes in a predictable manner
as $c$ varies, and we can use this fact to come up with the following criteria for
the map $u$ to be an embedding.
As with \fullref{t:no-isect-intro}, this result can be viewed as saying that certain subsets of
the ``obvious'' necessary conditions for $u$ to be an embedding are also sufficient.
In condition $(4)$ below, the
relative asymptotic intersection number
\[
\ain^\Phi([\tl u; z], [\tl u; w])
\]
is a count of how many intersections appear between a neighborhood of $z\in\Gamma$ and a neighborhood of $w\in\Gamma$
when $\tl u$ is perturbed near $w$ in a direction determined by a trivialization $\Phi$ of the
asymptotic limit of $\tl u$ at $w$.  More detail on this quantity can be found in
\fullref{ss:asympt-int} below.

\begin{theorem}\label{t:embedded-projection-intro}
Let $[\Sigma, j, \Gamma, \tl u=(a,u)]\in\M(M, \Ha, J)$
be a connected, simple pseudoholomorphic curve, and assume
that $\tl u$ does not have image contained in an orbit cylinder.
Then the following are equivalent:
\begin{enumerate}
\item The projected map $u:\Sigma\setminus\Gamma\to M$ is an embedding.

\item The algebraic intersection number $\inum(\tl u, \tl u_{c})$ between $\tl u$ and an $\R$-translate
     $\tl u_{c}=(a+c, u)$ is zero for all $c\in\R\setminus\br{0}$.

\item All of the following hold:
     \begin{enumerate}
          \item $u$ does not intersect any of its asymptotic limits.
          \item If $\gamma$ is a periodic orbit so that $u$ is asymptotic at $z\in\Gamma$ to $\gamma^{m_z}$
          and $u$ is asymptotic at $w\in\Gamma$ to $\gamma^{m_w}$, then
          \[
          \tfrac{\winfty(\tl u; z)}{m_z}=\tfrac{\winfty(\tl u; w)}{m_w}.
          \]
     \end{enumerate}

\item All of the following hold:
     \begin{enumerate}
          \item The map $\tl u$ is an embedding.
          \item The projected map $u$ is an immersion which is everywhere transverse to $\X$
          \item For each $z\in \Gamma$, we have
          \[
          \gcd(m_z, \winfty(\tl u; z))=1.
          \]
          \item If $\gamma$ is a periodic orbit so that $u$ is asymptotic at $z$ to $\gamma^{m_z}$ and
          $u$ is asymptotic at $w\ne z$ to $\gamma^{m_w}$ with $m_{z}m_{w}>0$, then the
          relative asymptotic intersection
          number of the ends $[\tl u; z]$ and $[\tl u; w]$ satisfies
          \[
          \ain^\Phi([\tl u; z], [\tl u; w])=-m_{z}m_{w}\max\br{\tfrac{\winfty^\Phi(\tl u; z)}{m_z}, \tfrac{\winfty^\Phi(\tl u; w)}{m_w}}.
          \]
     \end{enumerate}
\end{enumerate}
\end{theorem}
We note that the conditions of this theorem simplify somewhat if every asymptotic limit of
$\tl u$ is geometrically distinct.  In particular, condition $(3)$ implies that
a simple curve with geometrically distinct asymptotic limits
projects to an embedding in $M$
if and only if the projection doesn't intersect any of its asymptotic limits.
For planes, this generalizes Theorem 1.1 from \cite{hwz:prop2}.

As an application of some of these ideas, we can prove a result about the contact homology of a manifold
admitting a holomorphic open book decomposition satisfying some additional assumptions.
Recall that an open book decomposition of a $3$-manifold $M$ is a 
pair $(L, \pi)$ where $L\subset M$ is a link, and
$\pi:M\setminus L\to S^1$ is a fibration, the fibers $\pi^{-1}(\theta)$ of which are
embedded surfaces bounded by $L$.
In this case the link $L$ is referred to as the \emph{binding} of the open book decomposition, and the fibers
of the fibration are referred to as \emph{pages} of the open book decomposition.

Now consider $M$ equipped with a stable Hamiltonian structure $\Ha$ and compatible
complex multiplication $J\in\J\MH$.  We say that $(M, \Ha, J)$ admits
a \emph{holomorphic open book decomposition}, if $M$ admits an open book decomposition
$(L, \pi)$ so that the link $L$ is comprised of elliptic periodic orbits of the vector field $\X$,
and the pages are projections of $\tildeJ$-holomorphic curves in $\R\times M$.
We say a holomorphic open book decomposition is \emph{stable} if the
pseudoholomorphic curves which project to pages are genus $0$,
have only positive simple punctures, and have Fredholm index
$2$.

\begin{theorem}\label{t:openbook-intro}
Assume that $(M,\Ha, J)$ admits a stable, holomorphic open book decomposition.
Let $[\Sigma, j, \Gamma, \tl u=(a, u)]\in\MM$ be a connected
pseudoholomorphic curve, and assume that the image of $u$ is not a page of the open book decomposition,
and that $\tl u$ does not have image contained in an orbit cylinder.
Then at least one of the following is true:
\begin{enumerate}
\item At least one of the positive punctures of $\tl u$ limits to an orbit
that is not a binding of the open book
decomposition.
\item At least one of the positive punctures of $\tl u$ limits to a multiple cover of a binding orbit of
the open book.
\end{enumerate}
\end{theorem}

This result in particular immediately implies the following corollary.
See \cite{egh} for the definition of contact homology and cylindrical contact homology.

\begin{corollary}
Let $(M, \Ha, J)$ admit a stable, holomorphic open book decomposition $(L, \pi)$,
let $\gamma\subset L$ be a binding orbit, and assume that the
contact homology (resp.\ cylindrical contact homology) of the triple $(M, \Ha, J)$ is well-defined.
Then $\gamma$ is a cycle
in the contact homology (resp.\ cylindrical contact homology) of $(M, \Ha, J)$.
\end{corollary}

\subsubsection{Generalized pseudoholomorphic curves}
Finally we present a result concerning the intersection properties of
so-called generalized pseudoholomorphic curves.  These curves are introduced
in \cite{h2000} as a possible way to extend the theory of finite-energy foliations (see e.g.\ \cite{hwz:fol})
to include curves with nontrivial genus.

The generalized pseudoholomorphic curve equation for
an asymptotically cylindrical map $(\Sigma, j, \Gamma, (a, u))\in C^{\infty}(M, \Ha)$
is obtained by twisting the first of the two equations in $\eqref{e:j-hol-cylindrical}$
by a harmonic form on $\Sigma$.
More precisely, we consider quintuplets
$(\Sigma, j, \Gamma, \tl u=(a,u), \nu)$ where
the quadruple $(\Sigma, j, \Gamma, (a, u))$
is a smooth asymptotically cylindrical map, (i.e. belongs to $C^{\infty}(M, \Ha)$),
and where
$\nu$ is a $1$-form on the unpunctured surface $\Sigma$.
Such a quintuple is called a \emph{generalized pseudoholomorphic map}
in $(M, \Ha, J)$
if it satisfies
\begin{equation}\label{e:gen-hol}
\left\{
\begin{gathered}
u^{*}\lambda\circ j=da+\nu \\
\pi_{\xh}\circ du\circ j=J\circ\pi_{\xh}\circ du \\
d\nu=0=d(\nu\circ j).
\end{gathered}
\right.
\end{equation}
A \emph{generalized pseudoholomorphic curve} is
the equivalence class $[\Sigma, j, \Gamma, \tl u=(a,u), \nu]$  of the quintuple
$(\Sigma, j, \Gamma, \tl u=(a,u), \nu)$ under the
the equivalence relation of holomorphic reparametrization of the domain.

As with pseudoholomorphic maps, a finite energy condition guarantees that the maps
are asymptotically cylindrical.  We will not give the definition of energy here, but
will instead just assume that we are dealing with asymptotically cylindrical maps.
We will denote the set of asymptotically cylindrical generalized pseudoholomorphic curves
in $(M, \Ha, J)$ by $\M_{\Delta}(M, \Ha, J)$

In contrast to genuinely pseudoholomorphic curves,
generalized pseudoholomorphic curves
do not satisfy local positivity of intersections.  Indeed,
it is possible to construct a pair of local solutions
$\tl u=(a, u)$ and $\tl v=(b, v)$ to \eqref{e:gen-hol}
which have an isolated intersection of negative index
(see \fullref{a:gen-hol}).
In light of this, one might expect that a pair of generalized pseudoholomorphic
curves could have arbitrary generalized intersection  number, but
as the following theorem shows, this is not the case.

\begin{theorem}\label{t:gen-hol-pos-intro}
Let $[\Sigma , j, \Gamma, \tl u, \nu]$, $[\Sigma', j', \Gamma', \tl v , \nu']\in\M_{\Delta}(M, \Ha, J)$
be generalized pseudoholomorphic curves, and assume that no
component of $\tl u$ or $\tl v$ is contained in an orbit cylinder.  Then
\[
\gin{\tl u}{\tl v}\ge 0.
\]
\end{theorem}

This result puts topological restrictions on what
homotopy classes in $C^{\infty}(M, \Ha)$ can contain generalized pseudoholomorphic maps.
Indeed we have the following immediate corollaries of \fullref{t:gen-hol-pos-intro}.

\begin{corollary}
Let $(\Sigma, j, \Gamma, \tl u)\in C^{\infty}(M, \Ha)$
be a smooth asymptotically cylindrical map with no component homotopic to an orbit cylinder.
If
\[
\gin{\tl u}{\tl u}<0
\]
then there are no generalized pseudoholomorphic curves in the same relative homotopy class
as $\tl u$.
\end{corollary}

\begin{corollary}
Let $(\Sigma, j, \Gamma, \tl u)\in C^{\infty}(M, \Ha)$
be a smooth asymptotically cylindrical map with no component homotopic to an orbit cylinder.
Assume there exists
a generalized pseudoholomorphic  curve
$[\Sigma, j, \Gamma, \tl v, \nu]\in\M_{\Delta}(M, \Ha, J)$
so that
\[
\gin{\tl u}{\tl v}<0.
\]
Then there are no generalized pseudoholomorphic curves in the same relative homotopy class
as $\tl u$.
\end{corollary}

\subsection{Overview}\label{ss:overview}
The remainder of the paper is devoted to proving the results above concerning
the properties and applications of
generalized intersection number.  We review relevant technical
background
material where appropriate.

As mentioned above, some difficulties arise in trying to understand the intersection
and embedding properties
of punctured pseudoholomorphic curves when multiple ends approach the same orbit,
or, when concerned with embedding controls, when an end approaches a multiple cover
of an orbit.  
First, it is not a priori clear that the algebraic intersection number
is finite since the domains are not compact.
Second, even if the algebraic intersection number is finite,
in families of curves, intersections (or self-intersections)
can appear or disappear at the ends and thus, the algebraic intersection number
is not
topologically determined.

A precise understanding of the asymptotic behavior of punctured pseudoholomorphic curves is
the key both 
to establishing finiteness of the algebraic intersection number and
to understanding to what degree intersections or self-intersections can appear or disappear at the
ends.
The relevant asymptotic results are proved in \cite{sief08} and are reviewed in
\fullref{ss:asymp} after reviewing facts about asymptotic operators and the Conley--Zehnder
index in \fullref{ss:cz}.  With the appropriate asymptotic results in hand,
we study some local asymptotic winding invariants in \fullref{ss:asymp-winding-nos}
which are convenient for encoding  intersection related invariants we develop.
Then in \fullref{ss:asympt-int} we take up the study of local asymptotic intersection
invariants.
In particular we
define and establish the main properties
of the asymptotic intersection index $\delta_{\infty}([\tl u; z], [\tl v; w])$ of two ends
and the asymptotic self-intersection index $\delta_{\infty}([\tl u; z])$ of a single end.
For pseudoholomorphic curves, these quantities can, in light of the asymptotic results of
\cite{sief08}, be thought of as measures respectively of tangency or self-tangency at infinity,
and, in this case, these quantities give the upper bound on the count of intersections or self-intersections
that could appear that those ends.  

The theorems stated above in \fullref{ss:intro-main} are then proved in \fullref{s:global}
and \fullref{s:cylindrical} along
with some complementary results.
We establish the basic properties of the relative and generalized intersection numbers
in \fullref{ss:gen-int} building up to the proof of \fullref{t:gin-positivity-intro}.
\fullref{ss:adjunction} is then devoted to proving
\fullref{t:adj-main-intro} and some related results.
The proofs of these theorems rest on the fact that, in the cases we consider,
the generalized intersection number can be understood completely by combining
facts about the local behavior of pseudoholomorphic curves
from \cite{mw} with the asymptotic intersection
invariants developed in \fullref{s:local}.
In \fullref{s:cylindrical}, we then specialize to the case of a cylindrical cobordism.
The key idea here is that the homotopy invariance of the generalized intersection number
coupled with the existence of the $\R$-action allows the generalized intersection number to
be computed in terms of quantities associated to the projection of the curve to the
three-manifold.

Finally we close with two appendices.  In \fullref{a:windpi} we
adapt a proof of Hofer, Wysocki, and Zehnder from \cite{hwz:prop2} to show that the projection to
$\xh$ along $\X$ 
of the differential of a connected curve
in $\MM$ 
either vanishes identically, or has a finite number of zeroes of
positive order.
In \fullref{a:gen-hol} we show that local solutions to the
generalized pseudoholomorphic curve equations \eqref{e:gen-hol}
can exhibit intersection behavior not found in genuinely pseudoholomorphic maps.

\section{Local intersection numbers at punctures}\label{s:local}

\subsection{Asymptotics and asymptotic winding numbers}\label{ss:asymptotics}
\subsubsection{Asymptotic operators and the Conley--Zehnder index}\label{ss:cz}
Throughout this section, we will assume that $(M, \Ha)$ is a $3$-manifold equipped with a
stable Hamiltonian structure $\Ha=(\lambda,\omega)$, and that $J\in\J(M, \Ha)$ is a compatible
complex multiplication on $\xh$.

We will associate to any periodic orbit a differential operator related to the linearized flow.
Let $\gamma\in\Pe(M,\Ha)$ be a $\tau$-periodic orbit (possibly multiply covered) and
let $h$ be a vector field along $\gamma$, that is $h:S^{1}\to TM$ is a smooth function satisfying
$h(t)\in T_{\gamma(t)}M$ for all $t\in S^{1}$.
Since $h$ is defined along a flow line of $\X$ we can define the Lie derivative $L_{\X}h$ of $h$ by
\[
L_{\X}h(t)=\left.\tfrac{d}{ds}\right|_{s=0}d\psi_{-s}(\gamma(t+s/\tau))h(t+s/\tau).
\]
Since the flow $\psi_{t}$ of $\X$ preserves the splitting \eqref{e:splitting},
so must $L_{\X}$, and we can conclude that if
$h(t)\in\xh_{\gamma(t)}$ for all $t\in S^{1}$,
then $L_{\X}h(t)\in\xh_{\gamma(t)}$ for all $t\in S^{1}$.
Moreover, if $\nabla$ is a symmetric connection on $TM$, we can use
$d\gamma(t)\dt=\tau\cdot\X(\gamma(t))$ to write
\[
\tau\cdot L_{\X}h=L_{\tau\cdot\X}h=\nabla_{\tau\cdot\X}h-\nabla_{h}(\tau\cdot\X)=\nabla_{t}h-\tau\nabla_{h}\X,
\]
and therefore the differential operator $\nabla_{t}\cdot-\tau\nabla_{\cdot}\X$ maps sections of
$\gamma^{*}\xh$ to
sections of $\gamma^{*}\xh$ and is independent of the choice of symmetric connection.

Choosing some $J\in\J(M,\Ha)$, we associate to each $\tau$-periodic orbit
$\gamma\in\Pe_{\tau}(M, \Ha)$
a differential operator
$\A_{\gamma,J}:C^{\infty}(\gamma^{*}\xh)\to C^{\infty}(\gamma^{*}\xh)$
acting on the space of smooth sections of $\xh$ along $\gamma$ defined by
\[
\A_{\gamma,J}\eta=-J(\nabla_{t}\eta-\tau\nabla_{\eta}\X).
\]
We note that the discussion of the previous paragraph implies that $\A_{\gamma,J}$ does in fact map the
space of sections of $\xh$ along $\gamma$ to itself, and that $\A_{\gamma,J}$ is independent of
symmetric connection $\nabla$ used to define it.
Define an inner product on $C^{\infty}(\gamma^{*}\xh)$ by
\[
\ip{h,k}_{J}=\int_{S^{1}}\omega_{\gamma(t)}\big(h(t),J(\gamma(t))k(t)\big)\, dt.
\]
Recalling that $L_{\X}\omega=0$, we have for any $h$, $k\in C^{\infty}(\gamma^{*}\xh)$
that
\[
\tfrac{\partial}{\partial t}\omega_{\gamma(t)}(h(t), k(t))
=\omega_{\gamma(t)}\bp{\tau(L_{\X}h)(t), k(t)}+\omega_{\gamma(t)}\bp{h(t), \tau(L_{\X}k)(t)}
\]
Noting that the compatibility of $J$ with $\omega|_{\xh}$
implies that $\omega(J\cdot, J\cdot)=\omega(\cdot, \cdot)$ on $\xh\times\xh$,
we integrate this to give
\[
\ip{h, \A_{\gamma, J}k}_{J}=\ip{ \A_{\gamma, J}h,k}_{J}
\]
Therefore $\A_{\gamma,J}$ is formally self-adjoint, and $\A_{\gamma,J}$
induces a self-adjoint operator
\[
\A_{\gamma,J}:D(\A_{\gamma,J})=H^{1}(\gamma^{*}\xh)\subset L^{2}(\gamma^{*}\xh)\to L^{2}(\gamma^{*}\xh).
\]
We will refer to $\A_{\gamma, J}$ as the {\em asymptotic operator} associated to the orbit $\gamma$.

The kernel of $\A_{\gamma, J}$ is closely related to the degeneracy of the orbit $\gamma$.
Indeed, if $\gamma$ is degenerate then there is a nonzero $v_{0}\in\xh_{\gamma(0)}$ so that
$d\psi_{\tau}(\gamma(0))v_{0}=v_{0}$.  Defining $v\in C^{\infty}(\gamma^{*}\xh)$
by $v(t)=d\psi_{\tau t}(\gamma(0))v_{0}$, we have that $L_{\X}v\equiv 0$, so
$v\in\ker\A_{\gamma, J}$.  Conversely, if $h\in C^{\infty}(\gamma^{*}\xh)$ satisfies
$h\in\ker\A_{\gamma, J}$, then we can conclude that $L_{\X}h\equiv 0$ so we must have
that $h(t)=d\psi_{\tau t}(\gamma(0))h(0)$, and in particular $h(0)=d\psi_{\tau}(\gamma(0))h(0)$.
Therefore $\gamma$ is nondegenerate if and only if $\A_{\gamma, J}$ has trivial kernel.

It will be important to understand the behavior  --  particularly the winding  --
of eigenvectors of asymptotic operators associated to periodic orbits in trivializations of $\xh$.
We establish our conventions for dealing with trivializations here.
Given a simple periodic orbit $\gamma\in\Pe^{0}(M, \Ha)$ we 
will use the term \emph{unitary trivialization} of $(\gamma^{*}\xh, \omega, J)$
to refer to a trivialization
\[
\Phi_{\gamma}:S^{1}\times\R^{2}\to\xh|_{\gamma}
\]
of $\gamma^{*}\xh$
satisfying
\begin{gather*}
\omega(\Phi_{\gamma}\cdot, \Phi_{\gamma}\cdot)=dx\wedge dy \\
\Phi_{\gamma}\circ J_{0}=J\circ\Phi_{\gamma}
\end{gather*}
where $J_{0}$ is the standard complex multiplication used to identify $\R^{2}$ with $\C$.
Having chosen a unitary trivialization $\Phi_{\gamma}$
of $(\gamma^{*}\xh, \omega, J)$, we get an induced trivialization
for $\xh$ along the multiply covered orbit $((\gamma^{m})^{*}\xh, \omega, J)$
by pulling back $\Phi_{\gamma}$ via the $m$-fold covering map
$[t]_{\R/\Z}\to[mt]_{\R/\Z}$.
When dealing with trivializations of $\xh$ along
an orbit $\gamma^{m}$, we will always assume that the trivialization arises from 
a choice of trivialization along the
underlying simply covered
orbit $\gamma\in\PMH$.  Moreover, we will generally
use the same notation to indicate the trivialization $\Phi$ of $\xh$ along $\gamma$,
and that induced on $\xh$ along $\gamma^{m}$.

In a unitary trivialization of the $(\gamma^{*}\xh, \omega, J)$
along a simple periodic orbit $\gamma:S^{1}\to M$,
the asymptotic operator $\A_{\gamma^{m},J}$ takes the form
\[
(\A_{\gamma^{m},J} h)(t)=-i\tfrac{d}{dt}h(t)-S(t)h(t),
\]
where $S(t)$ is a symmetric, two-by-two matrix.
An eigenvector of $\A_{\gamma^{m},J}$ satisfies a linear, first-order ordinary differential equation
 and therefore never vanishes
since it doesn't vanish identically.
Hence every eigenvector gives a map from $S^{1}\to\R^{2}\setminus\{0\}$
and thus has a well defined winding number.
Since $-i\frac{d}{dt}-S(t)$ is a compact perturbation of $-i\frac{d}{dt}$, it can be shown
that the winding is monotonic in the eigenvalue,
and that to any $k\in\Z$ there is a two dimensional space of eigenvectors with winding $k$.
These results are proved in Section 3 of \cite{hwz:prop2}, and we restate them here as a lemma.

\begin{lemma}\label{l:operator-spectrum}
Let $\gamma\in\Pe^{0}\MH$ be a simple periodic orbit of $\X$, let $\A_{\gamma^{m}, J}$
denote the asymptotic operator of $\gamma^{m}$ for $m\in\Z\setminus\br{0}$, and let
$\mathfrak{T}(\gamma^{*}\xh)$ denote
the set of homotopy classes of unitary trivializations of $(\gamma^{*}\xh, \omega, J)$.  There exists a map
$w:\sigma(\A_{\gamma^{m}, J})\times\mathfrak{T}(\gamma^{*}\xh)\to\Z$ which satisfies
\begin{enumerate}
\item If $e:S^{1}\to (\gamma^{m})^{*}\xh$ is an
eigenvector of $\A_{\gamma^{m}, J}$ with eigenvalue $\lambda$, then $w(\lambda, [\Phi])=\wind(\Phi^{-1}e)$, 
that is, $w(\lambda, [\Phi])$ measures the winding with respect to $\Phi$ of any eigenvector of
$\A_{\gamma^{m}, J}$ with eigenvalue $\lambda$.

\item For any fixed $[\Phi]\in\mathfrak{T}(\gamma^{*}\xh)$ we have that
\[
w(\lambda, [\Phi])<w(\mu, [\Phi])\Rightarrow \lambda<\mu,
\]
that is, the winding of eigenvectors of $\A_{\gamma^{m}, J}$ is (not strictly) monotonic in the eigenvalue.

\item If $m(\lambda)=\dim\ker(\A_{\gamma^{m}, J}-\lambda)$ denotes the multiplicity of $\lambda$ as an eigenvalue we have for every
$k\in\Z$ and $[\Phi]\in\mathfrak{T}(\gamma^{*}\xh)$ that
\[
\sum_{\{\lambda \,|\, w(\lambda, [\Phi])=k\}}m(\lambda)=2,
\]
that is, the space of eigenvectors of $\A_{\gamma^{m}, J}$ with any given winding has dimension $2$.
\end{enumerate}
\end{lemma}

In this paper the relationship between the spectrum of the asymptotic operator associated to a simply
covered orbit $\gamma$ and that of the asymptotic operator associated to its $k$-fold iterate
$\gamma^{k}$ will be important.
First note that there is a $\Z_{k}$ action on sections of $(\gamma^{k})^{*}\xh$ defined by
$[j]_{\Z/k\Z}*_{k}f(t)=f(t+\tfrac{j}{k})$.  It is easily verified that $\A_{\gamma^{k}, J}$ is equivariant with respect
to this action, and therefore that the eigenspaces of $\A_{\gamma^{k}, J}$ are fixed by this action.
We say an eigenvector $e$ of $\A_{\gamma^{k}, J}$ is multiply covered if the isotropy group
\[
G(e)=\br{j\in\Z_{k}\,|\,j*_{k}e=e}
\]
is nontrivial, and we define the covering number of $e$
\[
\cov(e)=|G(e)|
\]
of $e$ to be the order of this group.  Put more simply, $\cov(e)$ is the largest positive integer $m$
dividing $k$ for which $(k/m)*_{k}\cdot$ fixes $e$.
It follows from the definition of $\A_{\gamma,J}$ that if $e$ is an eigenvector of $\A_{\gamma,J}$
with eigenvalue $\lambda$, then the section $e_{k}\in C^{\infty}((\gamma^{k})^{*}\xh)$
defined by $e_{k}(t)=e(kt)$ is an eigenvector of $\A_{\gamma^{k}, J}$ with eigenvalue $k\lambda$.
Furthermore, it is straightforward to see that any eigenvector $e$ of
$\A_{\gamma^{k},J}$ with covering number $m=\cov(e)$ is of the form $e(t)=f(mt)$ for some eigenvector
$f$ of $\A_{\gamma^{k/m}, J}$ with $\cov(f)=1$.

The following lemma is an easy consequence of the preceding discussion and \fullref{l:operator-spectrum}.
\begin{lemma}\label{l:covgcd}
If $e$ is an eigenvector of $\A_{\gamma^{k}, J}$ and $\Phi$ is a trivialization of $\gamma^{*}\xh$, then
\[
\cov(e)=\gcd(\wind(\Phi^{-1}e), k).
\]
\end{lemma}

Following \cite{hwz:prop2} we define the Conley--Zehnder index of a periodic orbit 
$\gamma^{k}$ in terms of the spectrum of the asymptotic operator $\A_{\gamma^{k}, J}$.  
We define a number $\lne(\gamma^{k})\in\sigma(\A_{\gamma^{k}, J})$ by
\begin{equation}\label{e:lne}
\lne(\gamma^{k})=\max(\sigma(\A_{\gamma^{k}, J})\cap\R^{-})
\end{equation}
i.e. so that $\lne(\gamma^{k})$ is the largest negative eigenvalue of $\A_{\gamma^{k}, J}$.
Given a trivialization $\Phi$ of $\gamma^{*}\xh$, we define
\begin{equation}\label{e:alpha}
\alpha^{\Phi}(\gamma^{k})=w(\lne(\gamma^{k}); [\Phi])
\end{equation}
so that $\alpha^{\Phi}(\gamma^{k})$ is the winding relative to $\Phi$ of any eigenvector of
$\A_{\gamma^{k}, J}$
having the largest possible negative eigenvalue.  We define the parity of $p(\gamma^{k})$ of
$\gamma^{k}$ by
\begin{equation}\label{e:parity}
p(\gamma^{k})=
\begin{cases}
0  & \text{ if $\exists\mu\in\sigma(\A_{\gamma^{k}, J})\cap\R^{+}$ with $w(\mu,[\Phi])=\alpha^{\Phi}(\gamma^{k})$ }\\
1  & \text{ otherwise }
\end{cases}
\end{equation}
and we note that this definition does not depend on the choice of $\Phi$.
We recall from \cite{hwz:prop2} that the parity defined here agrees with the designation 
of $\gamma^{k}$ as even or odd as defined in \fullref{ss:ham}.
The Conley--Zehnder index $\mu^{\Phi}(\gamma^{k})$ of the orbit $\gamma^{k}$ relative to the trivialization $\Phi$
is then defined by
\begin{equation}\label{e:cz}
\mu^{\Phi}(\gamma^{k})=2\alpha^{\Phi}(\gamma^{k})+p(\gamma^{k}).
\end{equation}
We will at times suppress the choice of trivialization in our notation for $\alpha$ or $\mu$, but it should always
be understood that a choice of trivialization is necessary to define these quantities.
We also observe that even though the asymptotic operator $\A_{\gamma^{k}, J}$ depends on a
choice of $J\in\J\MH$, the Conley--Zehnder index of an orbit is independent of this choice, as can
be seen by the alternate definitions given in \cite{hwz:prop2}.

Understanding how the Conley--Zehnder index of $\gamma^{k}$
behaves as a function of $k$ for fixed $\gamma$ will be important.
This is well understood and the important facts will be listed in the following lemma, the proof of which follows
from basic facts about the symplectic group $Sp(1)$ which can be found e.g.\ in Appendix 8.1 in \cite{hwz:fol} or
the first chapter of \cite{abb01}.  For our purposes, it will be more convenient to state an iteration formula
for $\alpha$ rather than $\mu$, but the relationship between these two is clear.

\begin{lemma}\label{l:alpha-iterate}
Let $\gamma$ be a (not necessarily simple)
periodic orbit, and assume that $\gamma^{m}$ is nondegenerate for all $m\in\Z\setminus\br{0}$.
Let $k\in\Z\setminus\br{0}$ be a
nonzero integer.
\begin{itemize}
\item
If $\gamma$ is an even orbit, then
\begin{equation}\label{e:alpha-even}
\alpha(\gamma^{k})=k\alpha(\gamma).
\end{equation}
\item
If $\gamma$ is an odd hyperbolic orbit, then
\begin{equation}\label{e:alpha-oddhyperbolic}
\alpha(\gamma^{k})=\fl{k(\alpha(\gamma)+\tfrac{1}{2})}=k\alpha(\gamma)+\tfrac{k-p(k)}{2}
\end{equation}
where $\fl{\cdot}$ is the least integer function, and $p(k)$ denotes the parity of the integer $k$.
\item
If $\gamma$ is an elliptic orbit, then there exists an irrational $\theta\in\R$ so that
\begin{equation}\label{e:alpha-elliptic}
\alpha(\gamma^{k})=\fl{k\theta}.
\end{equation}
\end{itemize}
\end{lemma}

It will be convenient for later to record how $\alpha$ changes with changes of trivialization.
The only subtlety here is that since we always deal with trivializations of the underlying simply
covered orbits, we get a factor which accounts for the covering number of the orbit.
We note that if $E\to S^{1}$ is a Hermitian line bundle with unitary trivializations $\Phi$ and $\Psi$, the map
$\Phi^{-1}\circ\Psi$ determines a map $S^{1}\to U(1)\approx S^{1}$.  We denote the degree of this map
by $\deg(\Phi^{-1}\circ\Psi)$. 
We now state the formula as a lemma and omit the easy proof.

\begin{lemma}\label{l:alpha-triv}
Let $\Phi$ and $\Psi$ be trivializations of $\gamma^{*}\xh$.  Then
\begin{equation}\label{e:alpha-triv}
\alpha^{\Phi}(\gamma^{k})=\alpha^{\Psi}(\gamma^{k})+k\deg(\Phi^{-1}\circ\Psi)
\end{equation}
for any nonzero $k\in\Z$.
\end{lemma}

We close this section by introducing a quantity that will be useful later.  Given any eigenvalue
$\mu\in\sigma(\A_{\gamma^{k}, J})$ it is clear from
\fullref{l:operator-spectrum} and \fullref{l:covgcd} that all eigenvectors
with eigenvalue $\mu$ have the same covering number.  We will therefore write
\[
\cov(\mu)
\]
to denote the covering number of any eigenvector with eigenvalue $\mu$.  Given any
simple periodic orbit $\gamma\in\Pe^{0}\MH$,
we will then define the \emph{spectral covering number} $\cover(\gamma^{k})$
of $\gamma^{k}$ to be
\begin{equation}\label{e:clne}
\cover(\gamma^{k}):=\cov(\lne(\gamma^{k}))
\end{equation}
and note that by \fullref{l:covgcd} that
\[
\cover(\gamma^{k})=\gcd(k, \alpha^{\Phi}(\gamma^{k}))
\]
where $\alpha^{\Phi}(\gamma^{k})$ is computed relative to any trivialization of $\gamma^{*}\xh$.
Note that even though $\alpha^{\Phi}(\gamma^{k})$ depends on a choice of trivialization,
\fullref{l:alpha-triv} above shows that changing the trivialization changes
$\alpha$ by an integer multiple of $k$, so the $\gcd$ on the right hand side of this equation is unchanged.

We note that as a consequence of \fullref{l:alpha-iterate}, the spectral covering number of a hyperbolic
orbit can be computed knowing just the covering number of the orbit, and whether the underlying simple
orbit is even or odd.  We state this result as a lemma and omit the straightforward proof.

\begin{lemma}
Let $\gamma$ be simple periodic orbit, and let $k\in\Z\setminus\br{0}$.
\begin{itemize}
\item If $\gamma$ is an even orbit, then
\[
\cover(\gamma^{k})=|k|.
\]

\item If $\gamma$ is an odd, hyperbolic orbit, then
\[
\cover(\gamma^{k})=
\begin{cases}
|k|/2 & \text{ if $k$ is even} \\
1 & \text{ if $k$ is odd}.
\end{cases}
\]
\end{itemize}
\end{lemma}

\subsubsection{Asymptotically cylindrical ends}\label{ss:ends}
Here we will introduce a notion that will be useful for framing the results of the next several sections.
In this section we continue to consider a $3$-manifold $M$ equipped with a stable Hamiltonian structure
$\Ha=(\lambda, \omega)$ and a compatible $J\in\J\MH$.

Consider quadruples of the form
$(\Sigma, j, z_*, u)$ where $(\Sigma, j)$ is a (not necessarily closed) Riemann surface without boundary,
$z_*\in\Sigma$ is a point, and $u:\Sigma\setminus\br{z_*}\to \R\times M$ is a smooth map.
We will call such a quadruple an \emph{asymptotically cylindrical end model} if there is a
periodic orbit $\gamma^{m}$ so that $u$ is asymptotically cylindrical over $\gamma^{m}$.

We define an equivalence relation on end models in the following way.
We say that asymptotically cylindrical end models
$(\Sigma, j, z_*, u)$ and $(\Sigma', j', w_*, v)$ are equivalent if
there exists an open neighborhood $U\subset\Sigma$ of $z_*$, and a holomorphic embedding
$\psi:U\to \Sigma'$ with $\psi(z_*)=w_*$ so that
\[
u=v\circ\psi.
\]
on $U\setminus\br{z_*}$.  An equivalence class $[\Sigma, j, z_*, u]$ of asymptotically cylindrical 
end models will
be referred to as a \emph{asymptotically cylindrical end}.

An asymptotically cylindrical end is said to be \emph{embedded} if it has a representative model
$(\Sigma, j, z_*, u)$
for which $u$ is an embedding.
A pair of asymptotically cylindrical ends are said to be
\emph{nonintersecting}, if they can be represented by models
 $(\Sigma, j, z_*, u)$ and $(\Sigma', j', w_*, v)$ satisfying
 \[
 u(\Sigma\setminus\br{z_*})\cap v(\Sigma'\setminus\br{w_*})=\emptyset.
 \]
 An asymptotically cylindrical end is said to be \emph{pseudoholomorphic} if it can be represented
 by a model  $(\Sigma, j, z_*, u)$ with $u:(\Sigma\setminus\br{z_{*}}, j)\to(\R\times M, \tildeJ)$
 a pseudoholomorphic map.

By choosing holomorphic coordinates near a given point, it is clear that any asymptotically cylindrical end
can be represented by a model of the form $(\DD, i, 0, u)$, where $\DD$ is the unit disk in $\C$ centered at
$0\in\C$.  We will call such a model, a unit disk model.  Given an asymptotically cylindrical end
$[\DD, i, 0, u]$
and a positive integer $m$, we can define
\emph{the $m$-multiple cover $m\cdot[\DD, i, 0, u]$ of $[\DD, i, 0, u]$} by
\[
m\cdot[\DD, i, 0, u]:=[\DD, i, 0, u\circ\phi_m]
\]
where $\phi_m:\DD\to\DD$ is the map $\phi(z)=z^m$.  It is a straightforward exercise in complex analysis to 
verify that the equivalence class of the end-model produced in this way does not depend on the
choice of unit disk model used in the definition.

It will be convenient for some of the following results to consider asymptotically cylindrical ends with
one additional piece of data.
We will call a quintuple $(\Sigma, j, z_*, \bar v, u)$ a \emph{decorated asymptotically cylindrical end model}
if the quadruple $(\Sigma, j, z_*, u)$ is an asymptotically cylindrical end model, and
\[
\bar v\in\Sigma^0_{z_*}:=(T_{z_*}\Sigma\setminus\br{0})/\R^+
\]
is a ray in the tangent space of $\Sigma$ at $z_*$.
The ray in the data for a decorated end model will be referred to as an
\emph{asymptotic marker} or \emph{decoration} of the end.
Two decorated end models
$(\Sigma, j, z_*, [X]_{\Sigma^0_{z_*}} , u)$ and $(\Sigma', j', w_*, [Y]_{{\Sigma'}^0_{w_*}}, v)$
are said to be equivalent when there is an open neighborhood $U\subset\Sigma$ of
$z_*$ and a holomorphic embedding $\psi:U\to\Sigma'$
with $\psi(z_*)=w_*$, $[d\psi(z_*)X]_{{\Sigma'}^0_{w_*}}=[Y]_{{\Sigma'}^0_{w_*}}$
and
\[
u=v\circ\psi
\]
on $U\setminus\br{z_*}$.  An equivalence class of decorated end models is called
a \emph{decorated asymptotically cylindrical end}.

It is clear that any decorated asymptotically cylindrical end can be represented by a model of the
form $(\DD, i, 0, [1]_{\C^*/\R^+}, u)$.
As with undecorated ends, we can define the $m$-multiple cover by
\[
m\cdot [\DD, i, 0, [1]_{\C^*/\R^+}, u]=[\DD, i, 0, [1]_{\C^*/\R^+}, u\circ\phi_m],
\]
and it is easily verified that
the equivalence class of end models produced in this way does not
depend on the unit disk model chosen to represent the end.

Finally, given a decorated end $[\Sigma, j, z_*, [X]_{\Sigma^0_{z_*}}, u]$, and a
$a\in\C^*$ we define the change of decoration map
taking $[\Sigma, j, z_*, [X]_{\Sigma^0_{z_*}}, u]$ to $a*[\Sigma, j, z_*, [X]_{\Sigma^0_{z_*}}, u]$
by defining
\[
a*[\Sigma, j, z_*, [X]_{\Sigma^0_{z_*}}, u]=[\Sigma, j, z_*, [aX]_{\Sigma^0_{z_*}}, u],
\]
and again, it is easily verified that this is a well-defined operation on ends
independent of the choice of representing model and $X\in [X]_{\Sigma^0_{z_*}}$.

\subsubsection{Asymptotic normal forms and some consequences}\label{ss:asymp}
Underlying all of the results of this paper is a precise description of the asymptotic behavior of
pseudoholomorphic half-cylinders, and of the relative asymptotic behavior of a pair of
half-cylinders limiting to the same periodic orbit.
The key results are proved in \cite{sief08} which builds on previous
results from \cite{hwz:prop1,kr,hwz:prop4,mora}.
In this section, we will review the relevant facts and some consequences.
Here we continue to assume
$(M, \Ha, J)$ to be a $3$-manifold equipped with a stable Hamiltonian structure
and compatible complex multiplication, but we now
assume $\Ha$ to be nondegenerate.

Let $[\Sigma, j, z_*, [X], \tl u=(a,u)]$ be a decorated
asymptotically cylindrical end in $(\R\times M, \tildeJ)$, and assume that
$u$ is asymptotically cylindrical over $\gamma^k$ for some $\tau$-periodic
$\gamma\in\PMH$ and $k\in\Z\setminus\br{0}$.
Moreover, we assume that if $\sigma:[0,\ep)\to\Sigma$ is a smooth curve satisfying
$\sigma(0)=z_*$ and $\sigma'(0)\in[X]$, then $\lim_{t\to 0}u(\sigma(t))=\gamma(0)$,
i.e.\ if we extend the projected map $u$ to the
 $S^{1}$ compactification of $\Sigma\setminus\br{z_{*}}$, the asymptotic marker
$[X]$ hits the point $\gamma(0)$ on the orbit.

The assumption that $\tl u$ is asymptotically cylindrical allows us
to find for some $R\in\R$ an embedding\footnote{
	We remark that the embedding $\psi$ above is not in general holomorphic with respect to the standard
	almost complex structure $j_{0}$ on $\c{R}\subset\R\times S^{1}=\C/i\Z$
	even if the end is pseudoholomorphic.
	However it is easily seen from the results in \cite{sief08}
	and the fact that nondegeneracy implies exponential convergence
	that if the end is pseudoholomorphic, then there exist positive constants
	$d$ and $M$ so that
	\[
	|\dbar_{j,j_{0}} \psi(s,t)|=|\tfrac{1}{2}(d\psi+j\circ d\psi\circ j_{0})(s,t)|
	\le M e^{-(d+2\pi)s}
	\]
	for all $(s,t)\in\c{R}$, where the norm $|\cdot|$ on $\Hom^{0,1}(\R\times S^{1}, T\Sigma)$
	is defined with respect to  the euclidean metric on $\R\times\R/\Z$ and any metric on the
	\emph{unpunctured} surface $\Sigma$.
	}
$\psi:\c{R}\to\Sigma\setminus\Gamma$,
with
\begin{gather}
\lim_{s\to\infty}\psi(s,t)=z_{*} \notag \\
\lim_{s\to\infty}-\tfrac{\partial_s\psi(s,0)}{|\partial_s\psi(s,0)|}\in[X]\in (T_{z_*}\Sigma\setminus\br{0})/\R^+
\label{e:decoration}
\end{gather}
and a map $U:\c{R}\to (\gamma^{k})^{*}\xh$ with
$U(s,t)\in\xh_{\gamma^{k}(t)}$ for all $(s,t)\in\c{R}$ so that
\begin{equation}\label{e:asymp-representative}
\tl u(\psi(s,t))=(k\tau s,\exp_{\gamma^{k}(t)}U(s,t)),
\end{equation}
where $\exp$ denotes the exponential map of the metric
\[
g_{\Ha, J}:=\lambda\otimes\lambda+\omega(\cdot, J\cdot)
\]
on $M$.
We will call a pair $(U, \psi)$ satisfying \eqref{e:decoration}-\eqref{e:asymp-representative} an
\emph{asymptotic representative of $[\Sigma, j, z_*, [X], u]$}.
It is clear from the $C^1$ convergence of $\tl u$ to $\R\times\gamma^k$ that
\eqref{e:decoration} and
\eqref{e:asymp-representative} uniquely determine $(U, \psi)$ up to restriction of the domain.

Now let
$[\Sigma', j', w_*, [Y], \tl v]$ be a second decorated asymptotically cylindrical end,
which is also asymptotically cylindrical over $\gamma^k$.
Then we can find an asymptotic representative $(V, \phi)$ of $\tl v$ near $w_*$ so that
we can write
\[
\tl v(\phi(s,t))=(k\tau s,\exp_{\gamma^{k}(t)}V(s,t))
\]
If the ends being considered are pseudoholomorphic, 
the asymptotic behavior of the difference of the maps $U$ and $V$ is given by the following
theorem.  The proof can be found in \cite{sief08}.

\begin{theorem}\label{t:asymp-diff}
Assume that $[\Sigma, j, z_*, [X], \tl u]$ and $[\Sigma', j', w_*, [Y], \tl v]$
are distinct, decorated pseudoholomorphic ends asymptotically cylindrical over $\gamma^{k}$,
with asymptotic representatives $(U, \psi)$ and $(V, \phi)$ respectively.
Then the difference $U-V$, where defined, can be written
\[
U(s,t)-V(s,t)=e^{\lambda s}[e(t)+r(s,t)]
\]
where $e$ is an eigenvector of the asymptotic operator $\A_{\gamma^{k}, J}$
with eigenvalue $\lambda<0$ and where $r$ satisfies
\[
|\nabla_{s}^{i}\nabla_{t}^{j}r(s,t)|\le M_{ij}e^{-ds}
\]
for all $(s,t)\in\c{R}$, $(i,j)\in\N^{2}$, and some appropriate constants $d>0$, $M_{ij}>0$.
\end{theorem}

Some remarks about this theorem are in order.
We first remark that if $(U, \phi)$ is an asymptotic representative of
$[\Sigma, j, z_{*}, [X], \tl u]$, then $(U(\cdot, \cdot+\tfrac{1}{k}), \phi(\cdot, \cdot+\tfrac{1}{k}))$
is an asymptotic representative of the end
\[
e^{i2\pi/k}*[\Sigma, j, z_{*}, [X], \tl u]=[\Sigma, j, z_{*}, [e^{i2\pi/k}X], \tl u].
\]
Thus if $e^{i2\pi l/k}*[\Sigma, j, z_{*}, [X], \tl u]\neq e^{i2\pi j /k}*[\Sigma', j', w_{*}, [Y], \tl v]$,
the previous theorem then lets us write
\[
U(s, t+\tfrac{l}{k})-V(s, t+\tfrac{j}{k})
=e^{\lambda_{lj}}\bbr{e_{lj}(t)+r_{lj}(s,t)}
\]
with each $\lambda_{lj}$, $e_{lj}$ and $r_{lj}$ satisfying that same properties
as the $\lambda$, $e$, and $r$ appearing in the theorem.  When $|k|>1$, we can apply this observation
to the case where $\tl v=\tl u$ so the theorem gives an asymptotic description of each of the maps
$U(s, t+\tfrac{i}{k})-U(s, t+\tfrac{j}{k})$.
Additionally, we remark that in the case that $\tl v$ is a cylinder over $\gamma^{k}$, we get that
$V\equiv 0$, so this theorem reduces to give an asymptotic description of a single half-cylinder
as in \cite{hwz:prop1, mora}.
Combining these observations, we get the following description of the 
asymptotic behavior of $U$ which is a refinement of the results of \cite{hwz:prop1, mora} in the
event that $|k|>1$.  For proof, see \cite{sief08}.

\begin{theorem}\label{t:asymp}
With $(U, \phi)$ as defined above,
there exists an $s_{0}\in\R$ so that for $(s, t)\in [s_{0}, \infty)\times S^{1}$, either
$U$ vanishes identically or
\[
U(s,t)=\sum_{i=1}^{N}e^{\lambda_{i}s}(e_{i}(t)+r_{i}(s,t))
\]
where
\begin{itemize}
\item The $\lambda_{i}$ are a sequence of negative eigenvalues of $\A_{\gamma^{k},J}$ which is strictly
   decreasing in $i$ (i.e.\ $\lambda_{j}<\lambda_{i}$ for $j>i$).
\item Each $e_{i}(\ne 0)$ is an eigenvector of $\A_{\gamma^{k},J}$ with eigenvalue $\lambda_{i}$.
\item The sequence of positive integers defined by setting $k_{1}=\cov(e_{1})$ and
   $k_{i}=\gcd(k_{i-1}, \cov(e_{i}))$, is strictly decreasing in $i$.
\item The $r_{i}$ satisfy $r_{i}(s, t)=r_{i}(s, t+\tfrac{1}{k_{i}})$.
   Moreover, each $r_{i}$ satisfies exponential decay estimates of the form
   \[
   \left|\nabla_{s}^{l}\nabla_{t}^{m}r_{i}(s,t)\right|<M_{lm}e^{-ds}.
   \]
\end{itemize}
\end{theorem}

The following two corollaries will be important for our later results.  For proofs see \cite{sief08}.

\begin{corollary}\label{c:no-isect-puncture}
Let $[\Sigma, j, z_*, \tl u]$ and $[\Sigma', j', w_*, \tl v]$ be pseudoholomorphic ends
in $\R\times M$.  Then precisely one of the following is true:
\begin{itemize}
\item The ends are equal, i.e.\ $[\Sigma, j, z_*, \tl u]=[\Sigma', j', w_*, \tl v]$.

\item One end is a multiple cover of the other, i.e.\  there exists an integer $m\ge 2$ so that
either
\[
m\cdot [\Sigma, j, z_*, \tl u]=[\Sigma', j', w_*, \tl v]
\] 
or
\[
[\Sigma, j, z_*, \tl u]=m\cdot [\Sigma', j', w_*, \tl v].
\]

\item The ends $[\Sigma, j, z_*, \tl u]$ and $[\Sigma', j', w_*, \tl v]$ are nonintersecting.
\end{itemize}
\end{corollary}

\begin{corollary}\label{c:factor-puncture}
Let $[\Sigma, j, z_*, \tl u]$ be a pseudoholomorphic end in $\R\times M$.
Then either $[\Sigma, j, z_*, \tl u]$ is embedded, or there exists an embedded pseudoholomorphic end
$[\Sigma', j', w_*, \tl v]$ and an integer $m\ge 2$ so that
\[
[\Sigma, j, z_*, \tl u]=m\cdot [\Sigma', j', w_*, \tl v].
\]
\end{corollary}

Combining these corollaries with results about the local behavior of pseudoholomorphic curves
due to Micallef and White \cite{mw} or McDuff \cite{mcd94} gives the following generalizations
of results that are well-known for closed curves.  

The first such result is that two connected curves either have the same image or intersect in at most
a finite set.
\begin{corollary}\label{c:finite-isect}
Let $(W, \barJ)$ be an almost complex $4$-manifold with Hamiltonian structured ends,
and let $[\Sigma_{i}, j_{i}, \Gamma_{i}, \tl u_{i}]\in\M(W, \barJ, \Ha^{+}, \Ha^{-})$
be connected pseudoholomorphic curves.
Then the sets $\tl u_{i}^{-1}(\tl u_{j}(\Sigma_{j}\setminus\Gamma_{j}))$
are either finite or equal to $\Sigma_{j}\setminus\Gamma_{j}$.
\end{corollary}
 
We next address the question of double points and singularities for a connected curve.
A pseudoholomorphic curve $[\Sigma, j, \Gamma, \tl u]\in\M(W, \barJ, \Ha^{+}, \Ha^{-})$ is called
\emph{multiply covered} it factors through another curve
$[\Sigma', j', \Gamma', \tl v]$ via a holomorphic map
$\psi:\Sigma\to \Sigma'$ with $\deg\psi\ge 2$.
If a curve is not multiply covered, it is said to be \emph{simple}.
We can now state the following corollary.

\begin{corollary}\label{c:finite-sing}
Let $(W, \barJ)$ be an almost complex $4$-manifold with Hamiltonian structured ends,
and let $[\Sigma, j, \Gamma, \tl u]\in\M(W, \barJ, \Ha^{+}, \Ha^{-})$ a connected,
simple pseudoholomorphic curve.
Then the set of double points of $\tl u$
\[
D(\tl u)=\br{(p,q)\in(\Sigma\setminus\Gamma)^{2}\,|\,\tl u(p)=\tl u(q),\,p\ne q}
\]
and the set of singular points\footnote{
Note that the fact that $\tl u$ is pseudoholomorphic implies that either
$\operatorname{rank} d\tl u(p)=2$ or $d\tl u(p)=0$ for every $p\in\Sigma\setminus\Gamma$.
}
of $\tl u$
\[
S(\tl u)=\br{p\in\Sigma\setminus\Gamma\,|\,d\tl u(p)=0 }
\]
are both finite.
\end{corollary}

A version of the following corollary is proved in \cite{hwz:prop1} in the case that the stable Hamiltonian structure
arises from a contact form.  The proof given there readily generalizes to any stable Hamiltonian structure.

\begin{corollary}\label{c:pi-du-punctures}
Let $[\Sigma, j, \Gamma, \tl u]\in\MW$
be a pseudoholomorphic curve, let $z_{*}\in\Gamma$ be a puncture,
and let
$\pi:T(\R^{\pm}\times M^{\pm})\to\xi^{\Ha^{\pm}}$ be the projection onto $\xi^{\Ha^{\pm}}$
over the cylindrical ends
determined by the splitting $T(\R^{\pm}\times M^{\pm})\approx\R\oplus\R X_{\Ha^{\pm}}\oplus\xi^{\Ha^{\pm}}$.
Choose an open neighborhood $U'\subset\Sigma$ of $z_{*}$ so that $\tl u(U'\setminus\br{z_{*}})$
lies entirely within one of the cylindrical ends so that $\pi\circ d\tl u$ is defined
on $U'\setminus\br{z_{*}}$.  Then there exists an open neighborhood
$U'\subset U$ of $z_{*}$ so that $\pi\circ d\tl u$ either vanishes identically on $U'$
or vanishes nowhere on $U'$.
\end{corollary}

\subsubsection{Asymptotic winding numbers}\label{ss:asymp-winding-nos}
In this section we will define some invariants associated to the ends of an asymptotically cylindrical
map.  Throughout $(W, \barJ)$ will denote a $4$-dimensional almost complex cobordism
with positive/negative ends
$(\R^{\pm}\times M^{\pm}, \tildeJ_{\pm}, \Ha_{\pm}=(\lambda_{\pm}, \omega_{\pm}))$.

Let $[\Sigma, j, z, [X], \tl u]$ be a decorated (not necessarily pseudoholomorphic)
end which is asymptotically cylindrical over
$\gamma^m$ for some  $\gamma\in\PMH$ and some $m\in\Z\setminus\br{0}$.
We will abbreviate the end by $[\tl u; z]$ for simplicity.  Define:
\begin{align}
& \alpha^{\Phi}(\tl u; z):=\alpha^{\Phi}(\gamma^{m})
\text{ with the right hand side
as defined in \eqref{e:alpha}.} \\
&\mu^{\Phi}(\tl u; z):=\mu^{\Phi}(\gamma^{m}) \text{ with the right hand side as defined in
\eqref{e:cz}.} \\
&\label{e:spec-cover-end}
\cover(\tl u; z):=\cover(\gamma^{m}) \text{ with the right hand side as defined in \eqref{e:clne}.}
\intertext{We note that each of these quantities depends only on the asymptotic limit of the end, and not on
the map, or the decoration.  In the case that the end $[\tl u; z]$ is pseudoholomorphic and the latter alternative
of \fullref{t:asymp} holds,
we define
}
\label{e:leading-evect}
& e_{1}(\tl u; z):=
\text{the leading eigenvector in the asymptotic expression of} \\
&\notag
\hskip.75in\text{$\tl u$ at $z$ from
\fullref{t:asymp}}.
\end{align}
If $|m|>1$ then $e_1(\tl u; z)$ depends on the choice of asymptotic marker, but
applying the change of decoration map
$(e^{i2\pi/m})\cdot\in\Hom(T_{z}\Sigma)$ has the effect of replacing the eigenvector
$e_1=e_1(\tl u; z)$ with $e_{1}(\cdot+\frac{1}{m})$.
Thus
$e_{1}(\tl u; z)$ determines a well-defined $m$-fold multisection of $\xh|_{\gamma}$
which doesn't depend on the choice of marker.

We say that the end $[\tl u; z]$ \emph{winds} if it can be represented by a model
$(\Sigma, j, z, [X], \tl u=(a, u))$ in which the map $\tl u$ does not intersect the cylinder over $\gamma^{m}$,
or equivalently, so that the projection $u$ of the map to $M$ does not intersect the asymptotic limit $\gamma$.
In this case, we define the \emph{asymptotic winding} $\winfty^{\Phi}(\tl u; z)$ of $\tl u$ at $z$
relative to the trivialization $\Phi$ of $\gamma^{*}\xh$ by choosing a small clockwise loop around
$z$, and defining $\winfty^{\Phi}(\tl u; z)$ to be the winding of the image of this loop around $\gamma$
computed in the trivialization $\Phi$.  The assumption that the image of $u$ does not intersect $\gamma$
implies that the resulting quantity is independent of the loop chosen.
In the case that $[\tl u; z]$ is a pseudoholomorphic end, it is an easy consequence of \fullref{t:asymp}
that
\begin{equation}\label{e:wind-infinity}
\winfty^{\Phi}(\tl u; z):=\wind(\Phi^{-1}e_{1}(\tl u; z))
\end{equation}
Note that the comments of the previous paragraph show that
$\winfty^{\Phi}(\tl u; z)$ depends only on the undecorated end $[\Sigma, j, z, \tl u]$ and not on the choice of
asymptotic marker.
The following lemma collects some useful facts about $\winfty$.

\begin{lemma}\label{l:wind-infinity}
Let $\Phi$ and $\Psi$ be trivializations of $\gamma^{*}\xh$,
and let $[\tl u; z]$ denote an end which winds and
which is asymptotically cylindrical over $\gamma^{m}$
with $m\in\Z\setminus\br{0}$.
Then
\begin{equation}\label{e:wind-infinity-triv}
\winfty^{\Phi}(\tl u; z)=\winfty^{\Psi}(\tl u; z)+m \deg(\Phi^{-1}\circ\Psi).
\end{equation}
and the quantity
\begin{equation}\label{e:def}
d_{0}(\tl u; z):=\alpha^{\Phi}(\tl u; z)-\winfty^{\Phi}(\tl u; z)
\end{equation}
is independent of the choice of trivialization.
Moreover, if $[\tl u; z]$ is pseudoholomorphic
then $d_{0}(\tl u; z)\ge 0$ and\footnote{
The significance of the quantity $\Done(\tl u; z)$ is, at this point, likely not apparent to the reader,
so we give a brief preview here.  As mentioned in the introduction, we will later introduce a quantity
$\delta_{\infty}(\tl u; z)$ for embedded ends which
is an upper bound on the number of self-intersections that can appear 
at that end if pseudoholomorphic.  The quantity $\Done(\tl u; z)$ will be twice the contribution
to $\delta_{\infty}(\tl u, z)$ (if any) arising from the failure of
the leading eigenvector in the asymptotic formula to have extremal winding,
i.e.\ failure of $d_{0}(\tl u; z)$ to be zero.
}
\begin{equation}\label{e:delta1}
\Done(\tl u; z):=(|m|-1)d_{0}(\tl u; z)-\cover(\tl u; z)+\cov(e_{1}(\tl u; z))
\end{equation}
is even, nonnegative, and
$\Done(\tl u; z)=0$ if and only if at least
one of the following holds:
\begin{itemize}
\item $|m|=1$
\item $d_{0}(\tl u; z)=0$
\item $d_{0}(\tl u; z)=1$ and $\cover(\tl u; z)=|m|$.
\end{itemize}
\end{lemma}

\begin{proof}
The change of trivialization formula follows as in \eqref{e:alpha-triv} and we omit the trivial proof.
The trivialization independence of
$\alpha^{\Phi}(\tl u; z)-\winfty^{\Phi}(\tl u;z)$ is then a consequence of the respective
change-of-trivialization formulas.
In the event that $[\tl u; z]$ is pseudoholomorphic, the nonnegativity of this quantity
has been previously proved in \cite{hwz:prop2}, and
is an easy consequence 
the monotonicity of winding relative to eigenvalue (\fullref{l:operator-spectrum}) and the
fact that $e_{1}(\tl u; z)$ has negative eigenvalue (\fullref{t:asymp}), and therefore must have
have winding no greater than $\alpha^{\Phi}(\tl u; z)$ (by definition of $\alpha$).

Next, still assuming that $[\tl u; z]$ is pseudoholomorphic, we consider the quantity
\[
\Done(\tl u; z)=(|m|-1)[\alpha^{\Phi}(\tl u; z)-\winfty^{\Phi}(\tl u; z)]
-\cover(\tl u; z)+\cov(e_{1}(\tl u; z)),
\]
and observe that if
$\winfty^{\Phi}(\tl u; z)=\alpha^{\Phi}(\tl u; z)$ then
$\cov(e_{1}(\tl u; z))=\cover(\tl u; z)$ as a consequence of \fullref{l:covgcd},
so then $\Done(\tl u; z)=0$.  On the other hand, if $|m|=1$, then we must have
$\cov(e_{1}(\tl u; z))=\cover(\tl u; z)=1$ since both of these quantities must divide $m$.
Again we have $\Done(\tl u; z)=0$, and we can conclude
that $\Done(\tl u; z)=0$ whenever $(|m|-1)d_{0}(\tl u; z)=0$.

Assuming then that $|m|\ge 2$ and $d_{0}(\tl u; z)\ge 1$, we get that
\[
\Done(\tl u; z)\ge |m|-1-\cover(\tl u; z)+\cov(e_{1}(\tl u; z))
\]
which must be nonnegative since $\cover(\tl u; z)\le |m|$
while $\cov(e_{1}(\tl u; z))\ge 1$.
Moreover, the only possibility for $\Done(\tl u; z)=0$ is if
$d_{0}(\tl u; z)=1$, $\cover(\tl u; z)=|m|$
and $\cov(e_{1}(\tl u; z))=1$.
If we assume that $\cover(\tl u; z)=|m|$, we can conclude from
\fullref{l:covgcd} that
\[
\gcd(|m|, \alpha^{\Phi}(\tl u; z))=\cover(\tl u; z)=|m|
\]
so $|m|$ divides $\alpha^{\Phi}(\tl u; z)$.
If we further assume that $d_{0}(\tl u; z)=1$, we can use this with \fullref{l:covgcd}
to argue
\begin{align*}
\cov(e_{1}(\tl u; z))&=\gcd(|m|, \winfty^{\Phi}(\tl u; z)) \\
&=\gcd(|m|, \alpha^{\Phi}(\tl u; z)-d_{0}(\tl u; z)) \\
&=\gcd(|m|, \alpha^{\Phi}(\tl u; z)-1) \\
&=1.
\end{align*}
We can conclude that if $d_{0}(\tl u; z)=1$ and $\cover(\tl u; z)=|m|$, then
$\cov(e_{1}(\tl u; z))=1$, and hence $\Done(\tl u; z)=0$.  We have thus shown
that $\Done(\tl u; z)\ge 0$ and $\Done(\tl u; z)=0$ if and only if
$(|m|-1)d_{0}(\tl u; z)=0$,  or $d_{0}(\tl u; z)=1$ and $\cover(\tl u; z)=|m|$.

Finally, to prove the evenness of $\Done(\tl u; z)$, we observe that, as a result of
\fullref{l:covgcd}, the quantity is of the form
\begin{equation}\label{e:evenness}
(m-1)(a-b)-\gcd(m, a)+\gcd(m, b)
\end{equation}
(where here $a=\alpha^{\Phi}(\tl u; z)$ and $b=\winfty^{\Phi}(\tl u; z)$).
If the term $(m-1)(a-b)$ is odd then it must be the case that $m$ is even and $a$ and $b$ have opposite parity.
In this case, it follows that $\gcd(m, a)$ and $\gcd(m, b)$ must have opposite parity, so the sum of the three terms
is even.  In the case that $(m-1)(a-b)$ is even, we have either that $m$ is odd or that $a$ and $b$ have the same
parity (or both).  In either case, it follows that $\gcd(m, a)$ and $\gcd(m, b)$ have the same parity, and
so the sum of the three terms must be even.
\end{proof}

Next we will introduce a quantity associated to an embedded end
which records some of the finer information about
the asymptotic behavior in the event that the asymptotic limit of the puncture is multiply covered.
Let $[\tl u; z]=[\Sigma, j, z, [X], \tl u]$ still denote a
decorated end in $(\R\times M, \tildeJ)$
which is asymptotically cylindrical over $\gamma^m$, but now assume
in addition that $[\tl u; z]$ is an embedded end and that $|m|\ge 2$.
Let $(U, \psi)$ be an asymptotic representative
of $[\tl u; z]$ as defined in the previous section, i.e.
\[
\tl u(\psi(s,t))=(m\tau s, \exp_{\gamma^{m}(t)}U(s,t)).
\]
The assumption that $[\tl u; z]$ is an embedded end implies that
$U(s,t)-U(s, t+\frac{j}{m})=0$ if and only if $j$ is a multiple of $m$.
Consequently, the winding of the map
\[
t\in S^{1}\mapsto \Phi^{-1}_{\gamma^{m}(t)}\bbr{U(s,t)-U(s, t+\tfrac{j}{m})}\in\R^{2}\setminus\br{0}
\]
is well defined and independent of $s$.
We define the
\emph{secondary winding}
of $\tl u$ at $z$ relative to the trivialization $\Phi$ by
\[
\wtwo^{\Phi}(\tl u; z):=\sum_{j=1}^{|m|-1}\wind \Phi^{-1}\bbr{U(s,\cdot )-U(s, \cdot+\tfrac{j}{m})}.
\]
As with the asymptotic winding, it is straightforward to show that the secondary winding
depends only on the undecorated end $[\Sigma, j, z, \tl u]$, and not on the choice of asymptotic marker.

In the event that the embedded end $[\tl u; z]$ is pseudoholomorphic,
the secondary winding can be computed in terms of the eigenfunctions from
\fullref{t:asymp-diff}.
Letting $(U, \psi)$ still denote an asymptotic representative of $[\tl u; z]$, we have
from \fullref{t:asymp-diff},
that
\[
U(s,t)-U(s,t+\tfrac{j}{m})=e^{\lambda_{j}s}[e_{j}(t)+r_{j}(s,t)]
\]
 with $\lambda_{j}<0$ an eigenvector of $\A=\A_{\gamma^{m}, J}$, $e_{j}\ne 0$ an eigenvector of $\A$
 with eigenvalue
$\lambda_{j}$, and $r_{j}$ converging exponentially to $0$.  We will write
\begin{equation}\label{e:e-delta-j}
e^{\Delta}_{j}(\tl u; z)=e_{j}
\end{equation}
to denote the eigenvector appearing in the formula for $U(s,t)-U(s+\tfrac{j}{m})$.
Since the $r_{j}$ converge to $0$ as $s\to\infty$,
it follows that
\[
\wtwo^{\Phi}(\tl u; z)=\sum_{j=1}^{|m|-1}\wind(\Phi^{-1}e^{\Delta}_{j}(\tl u; z)).
\]

We collect some useful properties of $\wtwo$ in the following lemma.

\begin{lemma}\label{l:wind2}
Let $[\tl u, z]=[\Sigma, j, z, \tl u]$ be an embedded end in
$(\R\times M, \tildeJ)$ and assume that $[\tl u; z]$ is asymptotically cylindrical
over $\gamma^m$.
If $\Phi$ and $\Psi$ are trivializations of $\gamma^{*}\xh$ then
\[
\wtwo^{\Phi}(\tl u; z)=\wtwo^{\Psi}(\tl u; z)+m(|m|-1)\deg(\Phi^{-1}\circ\Psi).
\]
If $[\tl u; z]$ is pseudoholomorphic, then the quantities\footnote{
As with $\Done(\tl u; z)$ before, the significance of the quantities $\Dtwo(\tl u; z)$
and $\Dtot(\tl u; z)$
is likely not apparent to the reader at this point,
so we give a brief preview here.  As mentioned before, we will later introduce a quantity
$\delta_{\infty}(\tl u; z)$ for embedded ends which
is an upper bound on the number of self-intersections that can appear 
at that end if pseudoholomorphic.  The quantity $\Dtwo(\tl u; z)$ will be twice the contribution
to $\delta_{\infty}(\tl u, z)$ (if any) arising from the failure of
the eigenvectors in the asymptotic formula 
from \fullref{t:asymp}
to be simply covered, while $\Dtot(\tl u; z)$ is the sum of $\Done(\tl u; z)$ and $\Dtwo(\tl u; z)$.
}
\begin{gather}
\Dtwo(\tl u; z):=(|m|-1)\winfty^{\Phi}(\tl u; z)-\cov(e_{1}(\tl u; z))+1-\wtwo^{\Phi}(\tl u; z) \label{e:delta2}  \\
\Dtot(\tl u; z)=(|m|-1)\alpha^{\Phi}(\tl u; z)-\cover(\tl u; z)+1-\wtwo^{\Phi}(\tl u; z) \notag
\end{gather}
are even, nonnegative, and independent of the choice of trivialization.
Moreover, we have that $\Dtwo(\tl u; z)=0$ if and only if the asymptotic representation of $\tl u$ near $z$
(from \fullref{t:asymp}) has at most two terms and the winding of the eigenvectors appearing in 
this formula differ precisely by $1$.
In particular
\[
\wtwo^{\Phi}(\tl u; z)=(|m|-1)\winfty^{\Phi}(\tl u; z)
\]
if and only if $\cov(e_{1}(\tl u; z))=1$.
\end{lemma}

\begin{proof}
The change of trivialization formula follows as in \eqref{e:alpha-triv} and we again omit the easy proof.
The trivialization independence of the two quantities
$\Dtwo(\tl u; z)$ and $\Dtot(\tl u; z)$
is then a straightforward consequence of the change of trivialization
formulas.

Assuming now that $[\tl u; z]$ is pseudoholomorphic, the nonnegativity and evenness of $\Dtot(\tl u; z)$
will follow from the nonnegativity and evenness of $\Dtwo(\tl u; z)$
along with \fullref{l:wind-infinity}, since
$\Dtot(\tl u; z)-\Dtwo(\tl u; z)=\Done(\tl u; z)$,
which was already shown to be nonnegative and even.

To prove that $\Dtwo(\tl u; z)$ is even and nonnegative,
we first write down an alternate formula
for $\wind_{2}$ in terms of the eigenvectors appearing in \fullref{t:asymp}.  According to the theorem,
we can write
\begin{equation}\label{e:asymp-wind2}
\tl u(\psi(s,t))=\left(m\tau s,\exp_{\gamma^m(t)}\sum_{i=1}^{N}e^{\lambda_{i}s}[e_{i}(t)+r_{i}(s,t)]\right)
\end{equation}
with $\lambda_{i}<\lambda_{i-1}$, the sequence
\[
k_{1}=\cov(e_{1}) \hskip.25in k_{i}=\gcd(k_{i-1}, \cov(e_{i}))
\]
strictly decreasing in $i$, and the $r_{i}$ satisfying $r_{i}(s,t+\tfrac{1}{k_{i}})=r_{i}(s,t)$
and converging exponentially to zero.
We observe that \fullref{l:operator-spectrum} implies that
$\wind \Phi^{-1}e_{i}\le \wind\Phi^{-1}e_{i-1}$ since $\lambda_{i}<\lambda_{i-1}$.
Since $k_{i-1}$ divides $\cov(e_{i-1})$ and
$\gcd(k_{i-1}, \cov(e_{i}))=k_{i}<k_{i-1}$ it follows that
that $\cov(e_{i})\ne\cov(e_{i-1})$, so \fullref{l:covgcd}
lets us conclude that we have the strict inequality
$\wind \Phi^{-1}e_{i}< \wind\Phi^{-1}e_{i-1}$.

Abbreviating $k_{0}=|m|$ and $w_{i}=\wind\Phi^{-1}e_{i}$,
we claim that we can conclude from \eqref{e:asymp-wind2} that
\begin{equation}\label{e:wind2-alt-form}
\wtwo^{\Phi}(\tl u; z)=\sum_{i=1}^{N}(k_{i-1}-k_{i})w_{i}.
\end{equation}
To see this, observe that if $e_{i}(\cdot)-e_{i}(\cdot+\tfrac{j}{k_{0}})$ is nonzero then
\fullref{l:operator-spectrum} implies that
\[
\wind\Phi^{-1}e_{i}=\wind\Phi^{-1}[e_{i}(\cdot)-e_{i}(\cdot+\tfrac{j}{k_{0}})]
\]
since $e_{i}$ and $e_{i}(\cdot)-e_{i}(\cdot+\tfrac{j}{k_{0}})$ are eigenvectors of $\A$ with the same
eigenvalue.
Then, we can use that $r_{i}(s,t+\tfrac{1}{k_{i}})=r_{i}(s,t)$
and that
\[
e_{i}(t+\tfrac{1}{k_{i}})=e_{i}(t+\tfrac{\cov(e_{i})/k_{i}}{\cov(e_{i})})=e_{i}(t)
\]
to conclude that if $(U, \psi)$ is an asymptotic representative of
$\tl u$, then as $j$ varies over $\br{1, \dots, k_{0}-1}$,
precisely $k_{i-1}-k_{i}$ of the terms $\Phi^{-1}[U(s, \cdot) -U(s, \cdot+\tfrac{j}{k_{0}})]$
have winding $w_{i}$ for $s$ large.
The formula \eqref{e:wind2-alt-form} follows immediately from this observation
and the definition of $\wtwo$.

Using the formula \eqref{e:wind2-alt-form}, the quantity $\Dtwo(\tl u; z)$
of interest can be written
\[
\Dtwo(\tl u; z)=(k_{0}-1)w_{1}-k_{1}+1-\sum_{i=1}^{N}(k_{i-1}-k_{i})w_{i}.
\]
Using that $k_{N}=1$ since $[u; z]$ is an embedded end, we can rewrite this as
\begin{align}
\Dtwo(\tl u; z)&=\sum_{i=1}^{N}(k_{i-1}-k_{i})w_{1}-k_{1}+1-\sum_{i=1}^{N}(k_{i-1}-k_{i})w_{i} \notag \\
&=1-k_{1}+\sum_{i=1}^{N}(k_{i-1}-k_{i})(w_{1}-w_{i}) \notag \\
\intertext{and since the $i=1$ term of the sum vanishes, we continue}
&=1-k_{1}+\sum_{i=2}^{N}(k_{i-1}-k_{i})(w_{1}-w_{i}) \notag \\
&=1-k_{1}+\sum_{i=2}^{N}\sum_{j=2}^{i}(k_{i-1}-k_{i})(w_{j-1}-w_{j}) \notag \\
&=1-k_{1}+\sum_{j=2}^{N}\sum_{i=j}^{N}(k_{i-1}-k_{i})(w_{j-1}-w_{j}) \notag \\
&=1-k_{1}+\sum_{j=2}^{N}(k_{j-1}-1)(w_{j-1}-w_{j}) \notag \\
&=\sum_{j=2}^{N}(k_{j-1}-1)(w_{j-1}-w_{j})+k_{j}-k_{j-1} \label{e:Dtwo-sum}.
\end{align}
Using that $w_{j-1}-w_{j}\ge 1$ for each $j$ and that the
$k_{j}$ are a decreasing sequence of positive integers,
we find that each term in this sum in nonnegative (in fact
positive except possibly when $j=N$).
Using the definition of $k_{j}$ with \fullref{l:covgcd}, we have that
\begin{equation}\label{e:kj}
k_{j}=\gcd(k_{j-1}, \cov(e_{j}(\tl u; z)))=\gcd(k_{j-1}, \gcd(k_{0},w_{j}))=\gcd(k_{j-1}, w_{j})
\end{equation}
since $k_{j-1}$ divides $k_{0}$.
Using further that $k_{j}$ divides
$w_{j}$, we can also write $k_{j}=\gcd(k_{j}, w_{j})$.  These observations allow us to rewrite this sum
as
\begin{equation}
\sum_{j=2}^{N}(k_{j-1}-1)(w_{j-1}-w_{j})+\gcd(k_{j-1}, w_{j})-\gcd(k_{j-1}, w_{j-1}).
\end{equation}
Comparing this with \eqref{e:evenness}, it easily follows that each term in the sum is even, and consequently
that the sum is even.

Finally, to see the last claims are true, we observe that it follows from
\eqref{e:Dtwo-sum}
that $\Dtwo(\tl u; z)=0$ if and only if $w_{j-1}-w_{j}=1$ and $k_{j}=1$ for all values of $j\in\Z\cap[2, N]$
for which they are defined.
Since the $k_{j}$ form a decreasing sequence of positive
integers, 
we conclude that $\Dtwo(\tl u; z)=0$ precisely when either
$k_{2}=1$ and $w_{1}-w_{2}=1$
or when  $k_{1}=1$, i.e.\ precisely when the
asymptotic representation \eqref{e:asymp-wind2}
has at most two terms with the winding of the two eigenvectors (in the case of two terms)
appearing in that formula differing by $1$.
Moreover, there is one term in this representation if and only if $k_{1}=1$, which happens precisely
when the leading eigenvector $e_{1}(\tl u; z)$ is simply covered, and in this case
$\Dtwo(\tl u; z)=0$ is equivalent to
\[
\wtwo^{\Phi}(\tl u; z)=(|m|-1)\winfty^{\Phi}(\tl u; z).
\proved\]
\end{proof}

We close this section by defining a quantity that records relative asymptotic information about
two nonintersecting ends limiting to a cover of the same orbit with the same covering number.
We again let
$[\tl u; z]=[\Sigma, j, z, [X], \tl u]$
be a decorated end which is asymptotically cylindrical
over $\gamma^m$, we let $[\tl v; w]=[\Sigma', j', w, [X], \tl v]$ be a second decorated
end which is also asymptotically cylindrical over $\gamma^m$,
and assume that the end $[\tl u; z]$ and $[\tl v; w]$ are nonintersecting.
If we let $(U, \psi)$ and $(V, \phi)$ be asymptotic representatives of $[\tl u;z]$ and
$[\tl v, w]$ respectively,
then the assumption that $[\tl u; z]$ and $[\tl v; w]$ are nonintersecting
implies that $U(s, t)-V(s, t+\frac{j}{m})$ is never zero.
Consequently, the maps
\[
t\in S^{1}\mapsto \Phi^{-1}_{\gamma^{m}(t)}\bbr{U(s, t)-V(s, t+\tfrac{j}{m})}\in\R^{2}\setminus\br{0}
\]
have well-defined winding which is independent of $s$.
We then define the 
\emph{total relative winding}
of $[\tl u; z]$ and $[\tl v; w]$ 
relative to a trivialization $\Phi$ of $\gamma^*\xh$ by
\[
\wrel^{\Phi}([\tl u; z], [\tl v; w]):=\sum_{j=0}^{|m|-1}\wind\Phi^{-1}\bbr{U(s, \cdot)-V(s, \cdot+\tfrac{j}{m})}.
\]
As with the secondary winding, the total relative winding of two nonintersecting ends is
easily verified to be independent of choice of asymptotic marker, and therefore gives a
well-defined invariant of the pair of undecorated ends
$[\Sigma, j, z, \tl u]$ and $[\Sigma', j', w, \tl v]$.

In the event that the ends $[\tl u; z]$ and $[\tl v; w]$ are pseudoholomorphic,
the total relative winding can be computed in terms of the eigenfunctions appearing in
\fullref{t:asymp-diff}.
Indeed the assumption that $[\tl u; z]$ and $[\tl v; w]$ are nonintersecting
tells us that for each $j\in \br{0, \dots, |m|-1}$ we can write
\[
U(s,t)-V(s,t+\tfrac{j}{k})=e^{\lambda_{j}s}[e_{j}(t)+r_{j}(s,t)]
\]
where $e_{j}$ is an eigenvector of the asymptotic operator with eigenvalue
$\lambda_{j}<0$ and where $r_{j}$ converges exponentially to zero.
The fact that the $r_{j}$ converge to zero as $s\to\infty$ lets us conclude that
\[
\wrel^{\Phi}([\tl u; z], [\tl v; w])=\sum_{j=0}^{|m|-1}\wind(\Phi^{-1}e_{j}).
\]

We collect some useful properties of $\wrel$ in the following lemma.

\begin{lemma}\label{l:wind-rel}
Let $[\tl u; z]=[\Sigma, j, z, \tl u]$ and $[\tl v; w]=[\Sigma', j', w, \tl v]$ be 
nonintersecting ends in $(\R\times M, \tildeJ)$ which are asymptotically
cylindrical over $\gamma^m$.
If $\Phi$ and $\Psi$ are trivializations of $\gamma^{*}\xh$, then
\[
\wrel^{\Phi}([\tl u; z], [\tl v; w])=\wrel^{\Psi}([\tl u; z], [\tl v; w])+m|m|\deg(\Phi^{-1}\circ\Psi)
\]
and
\[
\wrel^{\Phi}([\tl u; z], [\tl v; w])=\wrel^{\Phi}([\tl v; w], [\tl u; z]).
\]
If $[\tl u; z]$ and $[\tl v; w]$ are pseudoholomorphic ends,
the quantity
\[
|m|\max\br{\winfty^{\Phi}(\tl u; z), \winfty^{\Phi}(\tl v; w)}-\wrel^{\Phi}([\tl u; z], [\tl v; w])
\]
is nonnegative and independent of choice of trivialization, and we have the strict inequality
\[
|m|\max\br{\winfty^{\Phi}(\tl u; z), \winfty^{\Phi}(\tl v; w)}-\wrel^{\Phi}([\tl u; z], [\tl v; w])>0
\]
only if $e_{1}(\tl u; z)=j*_{m}e_{1}(\tl v; w)$ for some $j\in\Z_{|m|}$, where $*_{m}$ denotes the $\Z_{|m|}$-action
on $(\gamma^{m})^{*}\xh$
(or equivalently if $e_{1}(\tl u; z)$ and $e_{1}(\tl v; w)$ determine the same multisection
of $\gamma^{*}\xh$).
\end{lemma}

\begin{proof}
As with \fullref{l:wind-infinity} and \fullref{l:wind2}, the change of trivialization formula follows just
as in \eqref{e:alpha-triv} and we omit the trivial proof.
The symmetry of the total relative winding is an easy consequence of the fact that
winding is fixed by negation and the $\Z_{m}$ action,
while the fact that the quantity
\[
|m|\max\br{\winfty^{\Phi}(\tl u; z), \winfty^{\Phi}(\tl v; w)}-\wrel^{\Phi}([\tl u; z], [\tl v; w])
\]
is independent of choice of trivialization follows from the change-of-trivialization formulas
for $\wrel$ and $\winfty$.

Now assume that $[\tl u; z]$ and $[\tl v; w]$ are pseudoholomorphic.
To see that
\[
|m|\max\br{\winfty^{\Phi}(\tl u; z), \winfty^{\Phi}(\tl v; w)}-\wrel^{\Phi}([\tl u; z], [\tl v; w])
\ge 0
\]
we observe that the differences
\[
U(s,t)-V(s,t+\tfrac{j}{k})
\]
cannot decay slower than both of $U$ and $V$.  This observation with the definition of $\wrel$ and the
monotonicity of winding with respect to eigenvalue imply that
\[
\wrel^{\Phi}([\tl u; z], [\tl v; w])\le |m|\max\br{\winfty^{\Phi}(\tl u; z), \winfty^{\Phi}(\tl v; w)}
\]
as claimed.
To prove the final claim, we first assume that $e_{1}(\tl u; z)$ and $e_{1}(\tl v; w)$ 
have different eigenvalues.  Using the symmetry of the relative winding,
we can assume without loss of generality
that the eigenvalue of $e_{1}(\tl u; z)$ is strictly larger than that of $e_{1}(\tl v; w)$.
In this case, it is immediate that $e_{1}(\tl u; z)$ must be the eigenvector appearing in the asymptotic
formula for $U(s, t)-V(s,t+\tfrac{j}{m})$ for all $j\in\Z_{m}$.  It follows from
the definition of the total relative winding, and the monotonicity of winding in the eigenvalue
(\fullref{l:operator-spectrum})
that
\[
\wrel^{\Phi}([\tl u; z], [\tl v; w])=|m|\winfty^{\Phi}(\tl u; z)
=|m|\max\br{\winfty^{\Phi}(\tl u; z), \winfty^{\Phi}(\tl v; w)}
\]
as claimed.  Next assuming that $e_{1}(\tl u; z)$ and $e_{1}(\tl v; w)$ have the same eigenvalue
(and hence the same winding), but
belong to different $\Z_{m}$-orbits.  Then it is again straightforward to argue that
the eigenvector appearing in the asymptotic formula of $U(s, t)-V(s,t+\tfrac{j}{m})$
is
\[
e_{1}(\tl u; z)-j*_{m}e_{1}(\tl v; w),
\]
and hence
\begin{align*}
\wrel^{\Phi}([\tl u; z], [\tl v; w])
&=\sum_{j=0}^{|m|-1}\wind\Phi^{-1}\bbr{e_{1}(\tl u; z)-j*_{m}e_{1}(\tl v; w)} \\
&=|m|\winfty^{\Phi}(\tl u; z)\bp{=|m|\winfty^{\Phi}(\tl v; w)} \\
&=|m|\max\br{\winfty^{\Phi}(\tl u; z), \winfty^{\Phi}(\tl v; w)}.
\end{align*}
The only remaining possibility is that $e_{1}(\tl u; z)$ and $e_{1}(\tl v; w)$
have the same eigenvalue
and that $e_{1}(\tl u; z)$ is in the $\Z_{m}$-orbit of $e_{1}(\tl v; w)$.
Hence, this must be true if
the strict inequality is to occur.
\end{proof}

\subsection{Local intersection theory at $\infty$}\label{ss:asympt-int}
We now study local intersection properties of asymptotically cylindrical maps near their punctures.
The key idea, considered in a special case by Kriener in \cite{kr} and more generally 
by Hutchings
in \cite{hut02},
is to perturb the maps near the ends
and compute an intersection number which depends on the direction of the perturbation.
We will see that for pseudoholomorphic curves,
this intersection number is always bounded from below by a number which also depends on the direction
of perturbation,
but that the difference between the intersection number and its lower bound is independent of
this choice.  We can thus associate to each pseudoholomorphic end or pair of ends, an
``asymptotic self-intersection/intersection number"
which, roughly speaking, measures tangencies (or self-tangencies)
at infinity.
Many of our results in this section are analogous to those in Section 6 of \cite{hut02}.
We present full proofs here because our point of view and
conventions are somewhat different than those used
in \cite{hut02}, and because in one case we achieve a stronger bound than what is given in
\cite{hut02} which is important for our applications.

We start by considering an embedded (not necessarily pseudoholomorphic)
end model $(\Sigma, j, z, \tl u)$ in $\R\times M$, which we assume is
asymptotically cylindrical over $\gamma^m$.
Let $(U, \psi)$ denote an asymptotic representative for $\tl u$, so that
\[
\tl u(\psi(s,t))=(m\tau s, \exp_{\gamma^{m}(t)}U(s,t))
\]
with $\psi:\c{R}\to\Sigma\setminus\br{z}$ an embedding converging to
$z$ as $s\in [R, \infty)$ approaches infinity.

Letting $\Phi:S^{1}\times\C\to\gamma^{*}\xh$ be a unitary trivialization
of $\gamma^{*}\xh$ and $\ep\in\R^{+}\subset\C$ be a positive real number
(thought of as lying in $\C$),
we define a map
$\tl u_{z,\ep, \Phi}:\c{R}\to\R\times M$ by
\[
\tl u_{z,\ep, \Phi}(s,t)=(m\tau s, \exp_{\gamma^{m}(t)}[U(s,t)+\Phi(mt)\ep])
\]
so that $\tl u_{z, \ep, \Phi}$ is a perturbation of $\tl u$ near $z$ in a direction determined by $\Phi$.

\begin{lemma}\label{l:self-int-end}
There exists an $\ep_{0}>0$ so that the intersection number
\[
\inum(\tl u|_{\psi(\c{R})}, \tl u_{z, \ep, \Phi} )
\]
is well-defined, independent of $\ep\in (0, \ep_{0})$, and invariant under homotopies
$\Phi_{\tau}:S^{1}\times\C\to\gamma^{*}\xh$ of unitary trivializations of $\gamma^{*}\xh$.
\end{lemma}

\begin{proof}
Abbreviating $C_{R}=\c{R}$ and $I_{R}=[R, R+1]\times S^{1}$, it follows from the asymptotic
behavior of $\tl u$ and the assumption that $\tl u|_{\psi(\c{R})}$ is an embedding that
we can find a number $\ep_{0}>0$ so that for any $\ep\in(0, \ep_{0})$ 
there exists an $R(\ep)$ so that
\[
\tl u(\psi(I_{R}\cup C_{R(\ep)}))\cap\tl u_{z, \ep, \Phi}(C_{R})=\emptyset
\quad\text{ and }\quad
\tl u(\psi(C_{R}))\cap\tl u_{z, \ep, \Phi}(I_{R}\cup C_{R(\ep)})=\emptyset
\]
for any unitary trivialization $\Phi$; that is, the preimages of all intersections of $\tl u$ and $\tl u_{z, \ep, \Phi}$
are contained in the cylinder $[R+1, R(\ep)]\times S^{1}$.

With this condition satisfied, we can apply standard transversality arguments to
find smooth homotopies
$v_{\cdot}, w_{\cdot}:[0,1]\times\c{R}\to\R\times M$ so that
$v_{0}=\tl u\circ\psi$ and $w_{0}=\tl u_{z, \ep, \Phi}$,
\begin{gather*}
v_{\tau}|_{I_{R}\cup C_{R(\ep)}}=\tl u\circ\psi|_{I_{R}\cup C_{R(\ep)}}
\quad\text{ and }\quad
w_{\tau}|_{I_{R}\cup C_{R(\ep)}}=\tl u_{z,\ep,\Phi}|_{I_{R}\cup C_{R(\ep)}}, \\
v_{\tau}(I_{R}\cup C_{R(\ep)})\cap w_{\tau}(C_{R})=\emptyset
\quad\text{ and }\quad
v_{\tau}(C_{R})\cap w_{\tau}(I_{R}\cup C_{R(\ep)})=\emptyset
\end{gather*}
for all $\tau\in[0, 1]$,
and so that $v_{1}$ and $w_{1}$ have only transverse intersections.
These conditions imply that the set
$S=\br{(p,q)\in C_{R}\times C_{R}\,|\, v_{1}(p)=w_{1}(q)}$ is contained in
$([R+1, R(\ep)]\times S^{1})^{2}$, and the transversality assumption implies that the set $S$ is finite.
We can therefore define the intersection number $\inum(v_{1}, w_{1})$ by counting with sign
the intersections of $v_{1}$ and $w_{1}$, and we define
\[
\inum(\tl u|_{\psi(\c{R})}, \tl u_{z, \ep, \Phi} )=\inum(v_{1}, w_{1}).
\]
Considering a generic path of homotopies shows that this number does not depend on the choice
of homotopy.  Moreover, if $\tl u^{\tau}_{z}$ and $\tl u^{\tau}_{z, \ep, \Phi}$
are any homotopies starting at $\tl u|_{\psi(\c{R})}$ and $\tl u_{z, \ep, \Phi}$ respectively, we will
have that
\[
\inum(\tl u|_{\psi(\c{R})}, \tl u_{z, \ep, \Phi} )=
\inum(\tl u^{\tau}_{z}, \tl u^{\tau}_{z, \ep, \Phi} )
\]
provided that there exists an $R_{1}$ so that
\[
\tl u^{\tau}_{z}(I_{R}\cup C_{R_{1}})\cap\tl u^{\tau}_{z, \ep, \Phi}(C_{R})=\emptyset
\quad\text{ and }\quad
\tl u^{\tau}_{z}(C_{R})\cap\tl u^{\tau}_{z, \ep, \Phi}(I_{R}\cup C_{R_{1}})=\emptyset
\]
for all $\tau\in[0,1]$.  In particular, this homotopy invariance implies that
the intersection number
$\inum(\tl u|_{\psi(\c{R})}, \tl u_{z, \ep, \Phi} )$ is independent of $\ep\in(0, \ep_{0})$ and
invariant under homotopies of unitary trivializations $\Phi_{\tau}:S^{1}\times\C\to\gamma^{*}\xh$.
\end{proof}

We use the above lemma to define an invariant of an embedded, asymptotically cylindrical
end model.
Let $(\Sigma, j, z, \tl u)$ still denote an embedded end model, let
$O$ be some open neighborhood of $z$, and let
$\tl u_{z, \ep, \Phi}$ denote a perturbation of $\tl u|_{O}$ near $z$
defined as above.
The preceding lemma shows that the intersection number
\[
\inum(\tl u|_{O}, \tl u_{z, \ep, \Phi})
\]
is independent of sufficiently small $\ep>0$ and depends only on the
map $\tl u$, and the homotopy class of
unitary trivialization $\Phi$.  Using this observation, we define
the \emph{relative asymptotic self-intersection number} $i_{\infty}^{\Phi}(\tl u; z)$ of
the end model $(\Sigma, j, z, \tl u)$
relative to $\Phi$ by
\[
\ain^{\Phi}(\tl u; z):=\inum(\tl u|_{O}, \tl u_{z, \ep, \Phi})
\]
for any sufficiently small $\ep>0$ so that the conclusions of the preceding lemma hold.

Using \fullref{t:asymp}, we can compute $\ain^{\Phi}(\tl u; z)$ in terms of the
winding invariants defined in the previous section.  We state this as a lemma.

\begin{lemma}\label{l:asint-compute}
Let $(\Sigma, j, z, \tl u)$ be an embedded, asymptotically cylindrical end model,
and let $\Phi$ be a unitary trivialization of the plane-field
$\xh$ along the asymptotic limit of $\tl u$ at $z$.  
Then the asymptotic self-intersection number
$\ain^{\Phi}(\tl u; z)$ of $\tl u$ at $z$ relative to $\Phi$ is given by
\[
\ain^{\Phi}(\tl u; z)=-\wtwo^{\Phi}(\tl u; z).
\]
\end{lemma}

We remark that this Lemma demonstrates that the asymptotic intersection number is an invariant
of the embedded pseudoholomorphic end $[\Sigma, j, z, \tl u]$ and does not depend on the
choice of embedded model $(\Sigma, j, z, \tl u)$ representing $[\Sigma, j, z, \tl u]$.

\begin{proof}
Assume that $\gamma^m$ is the asymptotic limit of the end, and let $(U, \psi)$ be an asymptotic representative
of the end, so we can write
\[
\tl u(\psi(s,t))=(m\tau s, \exp_{\gamma^{m}(t)}U(s,t))
\]
Note that the assumption that $\tl u$ is an embedding implies
$U(s,t)-U(s, t+\tfrac{j}{m})$ has no zeroes for $j\ne 0\mod m$.
By definition, the asymptotic self-intersection number $\ain^{\Phi}(\tl u; z)$ is computed by
counting signed intersections of $\tl u$ and the map
\[
\tl u_{z, \ep, \Phi}(s,t)=(m\tau s, \exp_{\gamma^{m}(t)}[U(s,t)+\Phi(mt)\ep]).
\]
Considering the representations given for these maps, we see that
\[
\tl u\circ\psi(s,t)=\tl u_{z. \ep, \Phi}(s', t')
\]
if and only if
\begin{gather*}
s=s', \\
t=t'+\tfrac{j}{m} \text{ for some $j\in\Z_{m}$, and} \\
U(s,t)=U(s,t+\tfrac{j}{m})+\Phi(mt)\ep,
\end{gather*}
where in the last condition we have applied the previous two and have used the $1$-periodicity of $\Phi$
and the $\tfrac{1}{m}$-periodicity of $\gamma^{m}$.
It clear then that the algebraic count of intersections of $\tl u$ with $\tl u_{z, \ep, \Phi}$ is equal
to the sum of 
the algebraic count of zeroes of
\[
U(s,t)-U(s,t+\tfrac{j}{m})-\Phi(mt)\ep
\]
as $j$ varies over $\Z_{m}$,
which in turn is equal to the algebraic count of zeroes of the function
$F_{j}:\c{R}\to\C$ defined by
\[
F_{j}(s,t)=\Phi(mt)^{-1}[U(s,t)-U(s,t+\tfrac{j}{m})]-\ep\in\C
\]
as $j$ varies over $\Z_{m}\setminus\br{0}$.

Since all zeroes of the $F_{j}$ lie in the interior of a compact cylinder of the form
$[R_{1}, R_{2}]\times S^{1}$, the algebraic count of zeroes of $F_{j}$ will be equal to
the difference in winding of the loop $t\mapsto F_{j}(s, t)$ for $s\ge R_{2}$ and $s\le R_{1}$.
For large $s'\ge R_{2}$, the assumed decay of $U-j*U$ implies that
\[
\wind(F_{j}(s',\cdot))
=\wind(\Phi(m\,\cdot)^{-1}[U(s,\cdot)-U(s,\cdot+\tfrac{j}{m})]-\ep)
=\wind(-\ep)=0.
\]
Next, assuming we've chosen $\ep<\inf_{t\in S^{1}}|U(R_{1},t)-U(R_{1}, t+\tfrac{j}{m})|$,
we have that
\begin{align*}
\wind(F_{j}(R_{1},\cdot))
&=\wind(\Phi(m\,\cdot)^{-1}[U(R_{1},\cdot)-U(R_{1},\cdot+\tfrac{j}{m})]-\ep) \\
&=\wind(\Phi(m\,\cdot)^{-1}[U(R_{1},\cdot)-U(R_{1},\cdot+\tfrac{j}{m})]).
\end{align*}
Since $\Phi(mt)^{-1}[U(s,t)-U(s,t+\tfrac{j}{m})]$
has no zeroes, we know that the winding is independent of $s$.
We these observations with the definition of $\wind_{2}$ from the previous section to conclude that
\begin{align*}
\ain^{\Phi}(\tl u; z)
&=\sum_{j=1}^{|m|-1} \wind(F_{j}(R_{2}, \cdot))- \wind(F_{j}(R_{1}, \cdot)) \\
&=-\sum_{j=1}^{|m|-1} \wind(\Phi(m\,\cdot)^{-1}[U(R_{1},\cdot)-U(R_{1},\cdot+\tfrac{j}{m})])\\
&=-\wtwo^{\Phi}(\tl u; z),
\end{align*}
as claimed.
\end{proof}

As an immediate corollary of this computation and \fullref{l:wind2} we get that the asymptotic
self-intersection number satisfies a number of properties.
Note that \eqref{e:ain-wind-bound} below is an improvement on
Lemma 6.7 in \cite{hut02}, and reduces to the result given there when
$\cov(e_{1}(\tl u; z))=1$.

\begin{corollary}
Let $[\tl u; z]=[\Sigma, j, z, \tl u]$ be an embedded end,
with asymptotic limit $\gamma^m$,
and let $\Phi$ and $\Psi$ be unitary trivializations
of $\gamma^*\xh$.
Then the relative asymptotic self-intersection number
of the end $[\tl u; z]$ satisfies
\[
\ain^{\Phi}(\tl u; z)
=\ain^{\Psi}(\tl u; z)-m(|m|-1)\deg(\Phi^{-1}\circ\Psi)
\]
If $[\tl u; z]$ is a pseudoholomorphic end, then:
\begin{itemize}
\item We have that
\begin{equation}\label{e:ain-wind-bound}
\ain^{\Phi}(\tl u; z)\ge -(|m|-1)\winfty^{\Phi}(\tl u; z)+\cov(e_{1}(\tl u; z))-1
\end{equation}
and that
\[
\ain^{\Phi}(\tl u; z)= -(|m|-1)\winfty^{\Phi}(\tl u; z)
\]
if and only if $\cov(e_{1}(\tl u; z))=1$.
\item In particular, we have that
\[
\ain^{\Phi}(\tl u; z)\ge -(|m|-1)\alpha^{\Phi}(\tl u; z)+\cover(\tl u; z)-1
\]
and that
\[
\ain^{\Phi}(\tl u; z)= -(|m|-1)\alpha^{\Phi}(\tl u; z)
\]
if and only if $d_{0}(\tl u; z)=0$ and $\cover(\tl u; z)=1$.
\end{itemize}
\end{corollary}

The lower bounds on $\ain^{\Phi}(\tl u; z)$ should be thought of as a sort of
``positivity of intersections at infinity'' for embedded pseudoholomorphic ends.
Indeed, we have shown that the quantity
\[
\ain^{\Phi}(\tl u; z)+(|m|-1)\alpha^{\Phi}(\tl u; z)
\]
is independent of choice of trivialization, and that
\begin{align*}
\ain^{\Phi}(\tl u; z)+(|m|-1)\alpha^{\Phi}(\tl u; z)
&=\wtwo^\Phi(\tl u; z)+(|m|-1)\alpha^{\Phi}(\tl u; z) \\
&=\Done(\tl u; z)+\Dtwo(\tl u; z)+[\cover(\tl u; z)-1]
\end{align*}
where each of the three terms on the right hand side of this equation are nonnegative quantities
when $[\tl u; z]$ is pseudoholomorphic.
We will define the \emph{asymptotic self-intersection index} $\delta_{\infty}(\tl u; z)$ of the end
$[\tl u; z]$ to be
\begin{align}\label{e:sint-index}
\delta_{\infty}(\tl u; z)&:=\tfrac{1}{2}\bbr{\ain^{\Phi}(\tl u; z)+(|m|-1)\alpha^{\Phi}(\tl u; z)-\cover(\tl u; z)+1} \\
&=\tfrac{1}{2}\Dtot(\tl u; z). \notag
\end{align}
We note that for pseudoholomorphic ends
the asymptotic self-intersection index takes values in the nonnegative integers
as a result of \fullref{l:wind2} and \fullref{l:asint-compute}.
The asymptotic self-intersection index can be thought of as a measure of 
the self-tangency at infinity of an embedded pseudoholomorphic end.  In any given trivialization, it counts 
one half the difference
between the number of intersections the end has with a perturbed copy of itself and the minimum
number that must occur for any embedded pseudoholomorphic end with the same asymptotic data.

We now move on to considering the asymptotic intersection properties of a pair of
ends of curves which are asymptotic to coverings of the same simple periodic orbit with the same sign.
Our assumptions will be that $(\Sigma, j, z, \tl u)$ and $(\Sigma', j', w, \tl v)$ are nonintersecting
(not necessarily pseudoholomorphic) end models in $\R\times M$.
We further assume that $\tl u$ is asymptotically
cylindrical near $z$ over $\gamma^{m_z}$, and $\tl v$ is asymptotically cylindrical near
$w$ to $\gamma^{m_w}$, and that $m_z$ and $m_w$ have the same sign.

We assume that $(U, \psi)$ and $(V, \phi)$ are asymptotic representatives of $\tl u$ and $\tl v$ respectively
so that we can write
\begin{align*}
\tl u\circ\psi(s,t)&=(m_z\tau s, \exp_{\gamma^{m_z}(t)}U(s,t)) \\
\tl v\circ\phi(s,t)&=(m_w\tau s, \exp_{\gamma^{m_w}(t)}V(s,t))
\end{align*}
for $(s, t)\in\c{R}$ for some $R\in\R$.
If $\Phi$ is a trivialization of $\gamma^{*}\xh$, we define a map
$\tl v_{w,\ep, \Phi}:\c{R}\to W$ as above by perturbing the map $\tl v$ on a neighborhood of $w$
in a direction determined by $\Phi$.  As above we have the following lemma.

\begin{lemma}\label{l:ain-twoends}
There exists an $\ep_{0}>0$ so that the intersection number
\[
\inum(\tl u, \tl v_{w, \ep, \Phi})
\]
is well-defined and independent of $\ep\in(0, \ep_{0})$ and invariant under homotopies
of unitary trivializations $\Phi$.
\end{lemma}

\begin{proof}
Again denoting $C_{R}=\c{R}$ and $I_{R}=[R, R+1]\times S^{1}$, it follows from the asymptotic
behavior of $(\tl u; z)$ and $(\tl v; w)$ that there exists an $\ep_{1}>0$ so that
\[
\tl u\circ\phi(I_{R})\cap \tl v_{w, \ep, \Phi}(C_{R})=\emptyset
\hskip10pt\text{ and }\hskip10pt
\tl u\circ\phi(C_{R})\cap \tl v_{w, \ep, \Phi}(I_{R})=\emptyset
\]
for and $\ep\in(0, \ep_{1})$ and any unitary trivialization $\Phi$ of $\gamma^{*}\xh$.
Moreover, the asymptotic behavior of $\tl u$ and $\tl v$ implies that there exists an $\ep_{2}$ so that
for any $\ep\in(0, \ep_{2})$ we can find an $R(\ep)$ so that
$\tl u\circ\phi(C_{R(\ep)})\cap\tl v_{w, \ep, \Phi}(C_{R})=\emptyset$ and
$\tl u\circ\phi(C_{R})\cap\tl v_{w, \ep, \Phi}(C_{R(\ep)})=\emptyset$
for any unitary $\Phi$.
Therefore choosing $\ep_{0}>0$ less than $\min\br{\ep_{1}, \ep_{2}}$, it follows that
for any $\ep\in (0, \ep_{0})$ we have
\[
\tl u\circ\phi(I_{R}\cup C_{R(\ep)})\cap \tl v_{w, \ep, \Phi}(C_{R})=\emptyset
\]
and
\[
\tl u\circ\phi(C_{R})\cap \tl v_{w, \ep, \Phi}(I_{R}\cup C_{R(\ep)})=\emptyset
\]
for any unitary $\Phi$.
Given this, the claims of the lemma follow from standard transversality and homotopy invariance
arguments as in \fullref{l:self-int-end}.
\end{proof}

We use this lemma to define a local invariant of a pair of nonintersecting asymptotically cylindrical ends.
Let $(\Sigma, j, z, \tl u)$ and $(\Sigma', j', w, \tl v)$ still
denote nonintersecting end models asymptotically cylindrical over
$\gamma^{m_z}$ and $\gamma^{m_w}$ respectively, with $m_z$ and $m_w$ having the same sign.
Let $O$ be some open neighborhood of $w$, and let $v_{\ep, \Phi}$ still denote a perturbation
of $\tl v|_O$ defined as above.
The preceding lemma shows that the intersection number
\[
\inum(\tl u, \tl v_{\ep, \Phi})
\]
is independent of sufficiently small $\ep>0$ and depends only on the maps $\tl u$ and $\tl v$ and the homotopy
class of the unitary trivialization $\Phi$.
Using this observation, we define the
\emph{relative asymptotic intersection number} of the end models
$(\Sigma, j, z, \tl u)$ and $(\Sigma', j', w, \tl v)$ relative to $\Phi$ by
\[
\ain^{\Phi}([\tl u; z], [\tl v; w]):=\inum(\tl u|_{O}, \tl v_{\ep, \Phi}).
\]
where $\ep>0$ is chosen small enough so that the conclusions of the previous lemma apply.

As with the relative asymptotic self-intersection number, 
we can use \fullref{t:asymp-diff} to compute the relative asymptotic intersection number 
in terms of the local winding-related invariants studied in the previous section.

\begin{lemma}\label{l:rain-compute}
Let $(\Sigma, j, z, \tl u)$ and $(\Sigma', j', w, \tl v)$ be nonintersecting asymptotically cylindrical end models,
and assume that $(\Sigma, j, z, \tl u)$ and $(\Sigma', j', w, \tl v)$ are asymptotically cylindrical over
$\gamma^{m_z}$ and $\gamma^{m_w}$ respectively with $m_z m_w>0$.
Then the
relative asymptotic intersection number of $(\tl u; z)$ and $(\tl v; w)$ is given by
\[
\ain^{\Phi}([\tl u; z], [\tl v; w])
=-\tfrac{1}{m_z m_w}\wrel^{\Phi}(|m_w|\cdot[\tl u; z], |m_z|\cdot[\tl v; w])
\]
where $\Phi$ is a unitary trivialization of $\gamma^{*}\xh$, and $m\cdot$ is the operation of multiply covering 
an asymptotically cylindrical end.
\end{lemma}

We remark that, as with the relative asymptotic self-intersection number, this lemma demonstrates that
the relative asymptotic intersection number is an invariant of the ends
$[\Sigma, j, z, \tl u]$ and $[\Sigma', j', w, \tl v]$, and does not depend on the choice of
models representing these ends.

\begin{proof}
We initially assume that the ends in question cover $\gamma$ with the same covering number,
i.e.\ that $m_z=m_w$.
In this case, an argument analogous to that in \fullref{l:asint-compute} shows that
\[
\ain^{\Phi}([\tl u; z], [\tl v; w])=-\wrel^{\Phi}([\tl u; z], [\tl v; w]).
\]
Moreover, it is an easy consequence of the definition of $\wrel$ that
\[
\wrel^{\Phi}(m\cdot[\tl u; z], m\cdot[\tl v; w])=m^{2}\wrel^{\Phi}([\tl u; z], [\tl v; w])
\]
for any positive integer $m$ so this is equivalent to the formula given in the statement of the lemma.

In the case that $[\tl u; z]$ and $[\tl v; w]$ cover $\gamma$ with different covering numbers
(i.e.\ $m_z\ne m_w$) then we can replace $[\tl u; z]$ and $[\tl v; w]$ with the
multiply covered ends $|m_w|\cdot [\tl u; z]$ and $|m_z|\cdot[\tl v; w]$ which both
have covering number $m_z|m_w|=|m_z|m_w$.  Then the reasoning of the
previous paragraph gives
\[
\ain^{\Phi}(|m_{2}|\cdot[\tl u; z], |m_{z}|\cdot[\tl v; w])
=-\wrel^{\Phi}(|m_{2}|\cdot[\tl u; z], |m_{w}|\cdot[\tl v; w]).
\]
Moreover, it is an easy consequence of the definition of
the relative asymptotic intersection number that
\[
\ain^{\Phi}(k_{1}\cdot[\tl u; z], k_{2}\cdot[\tl v; w])
=k_{1}k_{2}\,\ain^{\Phi}([\tl u; z], [\tl v; w])
\]
for any positive integers $k_{1}$ and $k_{2}$.
The claim of the lemma follows immediately from these two equations.
\end{proof}

This lemma used with the basic results about the total relative winding yields the following corollary.
Note that \eqref{e:ain-bound-wind-infty} below is analogous to Lemma 6.9 in \cite{hut02}.

\begin{corollary}\label{c:ain}
Let $[\tl u; z]$ and $[\tl v; w]$ be nonintersecting ends
asymptotically cylindrical over $\gamma^{m_z}$ and $\gamma^{m_w}$ respectively,
and let $\Phi$ and $\Psi$ be unitary trivializations of
$\gamma^{*}\xh$.  Then the relative asymptotic intersection
number satisfies
\[
\ain^{\Phi}([\tl u; z], [\tl v; w])=\ain^{\Phi}([\tl v; w], [\tl u; z])
\]
and
\[
\ain^{\Phi}([\tl u; z], [\tl v; w])=\ain^{\Psi}([\tl u; z], [\tl v; w]) -m_z |m_w|\deg(\Phi\circ\Psi).
\]
If $[\tl u; z]$ and $[\tl v; w]$ are pseudoholomorphic, then
\begin{gather}
\ain^{\Phi}([\tl u; z], [\tl v; w])\ge -m_z m_w
\max\br{\tfrac{\winfty^{\Phi}(\tl u; z)}{|m_z|}, \tfrac{\winfty^{\Phi}(\tl v; w)}{|m_w|}} \label{e:ain-bound-wind-infty} \\
\ain^{\Phi}([\tl u; z], [\tl v; w])\ge -m_z m_w 
\max\br{\tfrac{\alpha^{\Phi}(\gamma^{m_z})}{|m_z|}, \tfrac{\alpha^{\Phi}(\gamma^{m_w})}{|m_w|}} \label{e:ain-bound-alpha}.
\end{gather}
\end{corollary}

\begin{proof}
The first two statements follow directly from
\fullref{l:rain-compute} and \fullref{l:wind-rel}.  The inequality \eqref{e:ain-bound-wind-infty}
follows from these same Lemmas and the additional elementary fact that
\[
\winfty^{\Phi}(m\cdot[\tl u; z])=m\winfty^{\Phi}(\tl u; z)
\]
for any pseudoholomorphic end $[\tl u; z]$ and any positive integer $m$.
Finally, the inequality \eqref{e:ain-bound-alpha} follows from \eqref{e:ain-bound-wind-infty}
and that
\[
\winfty^{\Phi}(\tl u; z)\le \alpha^{\Phi}(\tl u; z)
\]
for pseudoholomorphic ends, as observed in \fullref{l:wind-infinity}.
\end{proof}

As with the case of a single end, we define the
\emph{asymptotic intersection index}, $\delta_{\infty}([\tl u; z], [\tl v; w])$,
of the pair of ends $[\tl u; z]$ and $[\tl v; w]$ satisfying the assumptions of the previous lemmas
by
\begin{equation}\label{e:def-asymp-int-index}
\delta_{\infty}([\tl u; z], [\tl v; w]):=
\ain^{\Phi}([\tl u ; z], [\tl v; w])+
m_z m_w
\max\br{\tfrac{\alpha^{\Phi}(\gamma^{m_z})}{|m_z|}, \tfrac{\alpha^{\Phi}(\gamma^{m_w})}{|m_w|}}
\end{equation}
and we note that
when $[\tl u; z]$ and $[\tl v; w]$ are pseudoholomorphic,
this quantity is nonnegative as a result of the preceding corollary.
The asymptotic intersection index can be thought of as a measure of the tangency at infinity of the two
pseudoholomorphic ends.
It counts the difference between the number of intersections that appear when one
end is perturbed and the minimum number that must appear for any two nonintersecting
pseudoholomorphic ends with the
given asymptotic data.

\section{Global intersection theory}\label{s:global}
\subsection{The generalized intersection number and positivity of intersections}\label{ss:gen-int}
\subsubsection{Relative intersection number}

In this section we introduce the relative intersection number, and establish its basic properties.
We consider two smooth asymptotically cylindrical maps
$(\Sigma, j, \Gamma, \tl u)$, $(\Sigma', j', \Gamma', \tl v)\in C^{\infty}(W, \Ha^{+}, \Ha^{-})$
where $W$ is a $4$-manifold with Hamiltonian structured cylindrical ends.
Given a trivialization $\Phi$ of the stable Hamiltonian hyperplane field 
along the periodic orbits, we define the relative intersection number
$\rin^{\Phi}(\tl u, \tl v)$ of $\tl u$ and $\tl v$ relative to $\Phi$ by
perturbing $\tl v$ near its punctures in a direction determined by
$\Phi$ and computing the algebraic intersection number
of $\tl u$ with the perturbation of $\tl v$.
More precisely, given a puncture $w\in\Gamma'$ of $\tl v$, we can find a neighborhood
$U$ of $w$ which gets mapped entirely within one of the cylindrical ends, and
so that
\[
\tl v(\phi(s,t))=(m s, \exp_{\gamma^{m}(t)}h(s,t))\in\R^{\pm}\times M^{\pm}
\]
for some simple periodic orbit $\gamma$ and $m\in\Z\setminus\br{0}$,
where $\phi:[R, \infty)\times S^{1}\to U\setminus\br{w}$ is a diffeomorphism.
Choosing some smooth cutoff function $\beta$ supported in $U$ and equal to $1$ in
a neighborhood of $w$, we define the perturbed end
$\tl v_{\Phi, \ep}^{U}:U\setminus\br{w}\to\R^{\pm}\times M^{\pm}$
by
\begin{equation}\label{e:local-perturbation}
\tl v_{\Phi, \ep}^{U}=(m s, \exp_{\gamma^{m}(t)}\bbr{h(s,t)
+\beta(\phi(s,t))\Phi_{\gamma}(mt)\ep})\in\R^{\pm}\times M^{\pm}
\end{equation}
where $\Phi_{\gamma}:S^{1}\times\C\to\gamma^{*}\xi^{\Ha^\pm}$ 
is a unitary trivialization of $\gamma^{*}\xi^{\Ha^\pm}$ and $\ep\in\R^{+}\in\C$
is a small, positive real number.
We define the map $\tl v_{\Phi, \ep}:\Sigma'\setminus\Gamma'\to W$
by replacing $\tl v$ in a neighborhood of each puncture with a perturbation
as in \eqref{e:local-perturbation}.
It then follows from standard arguments that there exists an $\ep_{0}>0$
so that for $\ep\in(0, \ep_{0})$ the algebraic intersection number
\[
\inum(\tl u, \tl v_{\Phi, \ep})
\]
is independent of all choices made except for that of the homotopy class of trivialization
$\Phi$.
We thus define the \emph{relative intersection number} $\rin^{\Phi}(\tl u, \tl v)$
of $\tl u$ and $\tl v$ relative to $\Phi$ by
\[
\rin^{\Phi}(\tl u, \tl v):=\inum(\tl u, \tl v_{\Phi, \ep})
\]
for some sufficiently small $\ep>0$.

We summarize important properties of the relative intersection number in the following proposition.

\begin{proposition}\label{p:rel-inum}
Let $\tl u$ and $\tl v$ satisfy the assumptions of the previous paragraph, and let
$\Phi$ be a trivialization of $\xi^{\Ha^\pm}$ along the asymptotic limits of $\tl u$ and $\tl v$.
Moreover, assume that at $z\in\Gamma$, $\tl u$ is asymptotic to
$\gamma_{z}^{m_{z}}$ and at $w\in\Gamma'$, $\tl v$ is asymptotic to $\gamma_{w}^{m_{w}}$.
Then the
relative intersection number $i^{\Phi}(\tl u, \tl v)$ satisfies:
\begin{enumerate}
\item $\rin^{\Phi}(\tl u, \tl v)$ depends only on the homotopy classes of $\Phi$, $\tl u$, and $\tl v$
\item $\rin^{\Phi}(\tl u, \tl v)=\rin^{\Phi}(\tl v, \tl u)$
\item If $\Psi$ is another trivialization of $\xi^{\Ha^{\pm}}$ along the
periodic orbit set of the ends, then
\begin{align*}
\rin^{\Phi}(\tl u, \tl v)&=\rin^{\Psi}(\tl u, \tl v)
+\sum_{\substack{(z, w)\in\Gamma_{+}\times \Gamma'_{+} \\ \gamma_{z}=\gamma_{w}}}
-m_{z}m_{w}\deg(\Phi_{\gamma_{z}}^{-1}\circ\Psi_{\gamma_{w}})\\
&\hskip.25in
+\sum_{\substack{ (z, w)\in\Gamma_{-}\times \Gamma'_{-} \\ \gamma_{z}=\gamma_{w} }}
m_{z}m_{w}\deg(\Phi_{\gamma_{z}}^{-1}\circ\Psi_{\gamma_{w}})
\end{align*}
\end{enumerate}
\end{proposition}

\begin{proof}
Part $(1)$ follows from standard transversality and homotopy-invariance
arguments, and we omit the straightforward details.

To see that part $(2)$ is true, it suffices construct homotopies
$\tl u_{\tau}$ and $\tl v_{\tau}$, so that
\begin{align*}
\tl u_{0}&=\tl u \\
\tl u_{1}&=\tl u_{\Phi, \ep} \\
\tl v_{0}&=\tl v_{\Phi, \ep} \\
\tl v_{1}&=\tl v
\end{align*}
and so that there exist neighborhoods $U$ of $\Gamma$ and $V$ of $\Gamma'$ so that
\[
\tl u_{\tau}(U\setminus\Gamma)\cap \tl v_{\tau}(\Sigma'\setminus\Gamma')=\emptyset
\]
and
\[
\tl u_{\tau}(\Sigma\setminus\Gamma)\cap \tl v_{\tau}(V\setminus\Gamma')=\emptyset.
\]
for all $\tau\in[0,1]$.  This again is straightforward and we omit the details.

Finally $(3)$, follows from an argument analogous to that in
Lemma 2.5 (b) in \cite{hut02}, and we refer the reader there for more detail.
\end{proof}

For pseudoholomorphic maps, it follows from \fullref{c:finite-isect} that the relative intersection number
can be computed by summing local intersection numbers, and asymptotic intersection numbers as defined in the
preceding section.
Before stating the result, we establish some notation and terminology.
Consider pseudoholomorphic curves $[\Sigma, j, \Gamma, \tl u]$, $[\Sigma', j', \Gamma', \tl v]\in\MW$.
A small neighborhood of a puncture $z\in\Gamma$ determines a pseudoholomorphic end
(as defined in \fullref{ss:ends})
which we will abbreviate $[\tl u; z]$.
Thus if $[\tl u; z]$ and $[\tl v; w]$ 
are nonintersecting and asymptotic to a covering of the same orbit with the same sign,
we can define the
relative asymptotic intersection number
\[
\ain^{\Phi}([\tl u; z], [\tl v; w])
\]
as in the previous section.  Otherwise, we define
\[
\ain^{\Phi}([\tl u; z], [\tl v; w])=0
\]
and we define the \emph{total relative asymptotic intersection number} $\ain^{\Phi}(\tl u, \tl v)$
of $[\Sigma, j, \Gamma, \tl u]$ and $[\Sigma', j', \Gamma', \tl v]$ by
\[
\ain^{\Phi}(\tl u, \tl v):=\sum_{(z, w)\in\Gamma\times \Gamma'}\ain^{\Phi}([\tl u; z], [\tl v; w]).
\]

Next, we say that 
$[\Sigma, j, \Gamma, \tl u]$ and $[\Sigma', j', \Gamma', \tl v]\in\MW$
have no common components if $\tl u^{-1}(\tl v(\Sigma'\setminus\Gamma'))$ does
not contain an open set.  In this case, it follows from
\fullref{c:finite-isect} that $\tl u$ and $\tl v$ intersect in a finite number of points.
We can thus define the algebraic intersection number of $\tl u$ and $\tl v$ by summing local intersection numbers.
We now state the result giving the relative intersection number of two
pseudoholomorphic curves.  This result is analogous to
Lemma 8.5 in \cite{hut02} in the case that the $S_{1}$ and $S_{2}$ considered there are pseudoholomorphic.
The result follows immediately from local intersection properties of pseudoholomorphic curves,
\fullref{c:finite-isect}, and the definitions of the terms involved,
and we omit the straightforward proof.

\begin{theorem}\label{t:rin-intasymp}
Let $[\Sigma, j, \Gamma, \tl u]$, $[\Sigma', j', \Gamma', \tl v]\in\MW$
be a pair of pseudoholomorphic curves having no common components.
Then
\begin{equation}\label{e:rel-int-compute}
\rin^{\Phi}(\tl u, \tl v)
     =\inum(\tl u, \tl v)
     +\ain^{\Phi}(\tl u, \tl v).
\end{equation}

\end{theorem}

\subsubsection{The generalized intersection number}
We now define the generalized intersection number of two asymptotically cylindrical maps.
This quantity will be an integer-valued symmetric product on the space of
(homotopy classes of)
smooth asymptotically cylindrical maps,
which for pseudoholomorphic curves will bound the algebraic intersection number from
above.

Motivated by the bounds on relative asymptotic intersection numbers from
\fullref{c:ain} above, we define a homotopy-invariant product
on asymptotically cylindrical maps by adding a term to the relative intersection pairing that
balances the trivialization dependence.
Consider two smooth asymptotically cylindrical maps
$(\Sigma, j, \Gamma, \tl u)$ and $(\Sigma', j', \Gamma', \tl v)\in\CW$
with $W$ as usual denoting a $4$-dimensional manifold with Hamiltonian structured ends.
Assume that at $z\in \Gamma_\pm$ that $\tl u$ is asymptotic to
$\gamma_z^{m_z}$, and at $w\in\Gamma'_\pm$ that
$\tl v$ is asymptotic to over $\gamma_w^{m_w}$.
Given a trivialization of $\xi^{\Ha^\pm}$ along the asymptotic limits of $\tl u$ and $\tl v$,  define a quantity
$\Omega^\Phi(\tl u, \tl v)$ by\footnote{
Note that this is the same quantity appearing in the sum on the right hand side of
formula \eqref{e:gin-def-intro} from the introduction, except here we
use that $\alpha^{\Phi}(\gamma^{k})=\fl{\mu^{\Phi}(\gamma^{k})/2}$ (see \eqref{e:cz})
and we explicitly separate the parts of the sum coming respectively from
positive punctures and negative punctures.
}
\begin{align*}
\Omega^\Phi(\tl u, \tl v)
&:=\sum_{\substack{(z,w)\in\Gamma_+\times\Gamma'_+ \\ \gamma_z=\gamma_w}}
m_z m_w \max\br{\tfrac{\alpha^\Phi(\gamma_z^{m_z})}{m_z}, \tfrac{\alpha^\Phi(\gamma_w^{m_w})}{m_w}} \\
&\hskip.25in+\sum_{\substack{(z,w)\in\Gamma_-\times\Gamma'_- \\ \gamma_z=\gamma_w}}
m_z m_w \max\br{\tfrac{\alpha^\Phi(\gamma_z^{m_z})}{-m_z}, \tfrac{\alpha^\Phi(\gamma_w^{m_w})}{-m_w}}.
\end{align*}
Note that this quantity depends only on the asymptotic data associated to the two maps, and
not on the maps themselves.
We then define the \emph{generalized intersection number} of two asymptotically cylindrical maps
by
\begin{equation}\label{e:gin-def}
\gin{\tl u}{\tl v}=\rin^{\Phi}(\tl u, \tl v)+\Omega^\Phi(\tl u, \tl v).
\end{equation}
We note that as an immediate consequence of the change of trivialization formulas for
$\rin^\Phi$ and $\alpha^\Phi$
from \fullref{p:rel-inum} and \fullref{l:alpha-triv}
the generalized intersection number does not depend on the choice of
trivialization.

We now prove \fullref{t:gin-prop} which gives the basic properties of the generalized
intersection number.  We restate the result here for the convenience of the reader.

\begin{proposition}[Properties of the generalized intersection number]
Let $W$, $W_{1}$ and $W_{2}$ be $4$-manifolds with Hamiltonian structured
cylindrical ends, and assume we can form the concatenation $W_{1}\odot W_{2}$.  Then:
\begin{enumerate}
\item  If $(\Sigma, j, \Gamma, \tl u)$ and $(\Sigma', j', \Gamma', \tl v)\in\CW$
are asymptotically cylindrical maps
then the generalized intersection
number $\gin{\tl u}{\tl v}$ depends only on the homotopy classes of $\tl u$ and $\tl v$.

\item For any $(\Sigma, j, \Gamma, \tl u)$ and $(\Sigma', j', \Gamma', \tl v)\in\CW$
\[
\gin{\tl u}{\tl v}=\gin{\tl v}{\tl u}.
\]

\item If $(\Sigma, j, \Gamma, \tl u)$, $(\Sigma', j', \Gamma', \tl v)$,
$(\Sigma'', j'', \Gamma'', \tl w)\in\CW$
then
\[
\gin{\tl u+\tl v}{\tl w}=\gin{\tl u}{\tl w}+\gin{\tl v}{\tl w}
\]
where ``$+$'' on the left hand side denotes the disjoint union of the maps $\tl u$ and $\tl v$.

\item
If $u_{1}\odot u_{2}$ and $v_{1}\odot v_{2}$ are asymptotically cylindrical buildings in
$W_{1}\odot W_{2}$ then
\[
\gin{u_{1}\odot u_{2}}{v_{1}\odot v_{2}}\ge \gin{u_{1}}{v_{1}}+\gin{u_{2}}{v_{2}}.
\]
Moreover, strict inequality occurs if and only if there is a periodic orbit $\gamma$ so
that $\tl u_{1}$ has a negative puncture asymptotic to $\gamma^m$, $\tl v_1$ has a negative puncture asymptotic
to $\gamma^n$, and both $\gamma^m$ and $\gamma^n$ are odd orbits.

\end{enumerate}
\end{proposition}

\begin{proof}
The first claim follows immediately from the fact that
$\ain^{\Phi}(\tl u, \tl v)$ and $\Omega^{\Phi}(\tl u, \tl v)$
are determined entirely by the homotopy classes of $\tl u$, $\tl v$, and $\Phi$.
Since their sum does not depend on a choice of trivialization, the generalized intersection number
only depends on the homotopy class of $\tl u$ and $\tl v$ in $\CW$.

The second claim follows immediately from the symmetry of the relative intersection number and
the quantity $\Omega^\Phi(\tl u, \tl v)$, and the third claim follows immediately from the definition
and basic properties of intersection numbers.

Before proving the final claim, we first establish some notation.
We first observe that, since we assume that the concatenation
$u_1\odot u_2$ is defined, the asymptotic data at the negative punctures of $u_1$ must
correspond with and match the asymptotic data at the positive punctures of $u_2$.
We will let $\Gamma_{1,-}=\br{z^-_{i}}$ denote the set of negative punctures of $u_1$
and $\Gamma_{2, +}=\br{z^+_{i}}$ denote the set of positive punctures of $u_2$.
Moreover we assume that at $z^+_{i}$, $u_2$ is asymptotic to the periodic orbit
$\gamma_{i}^{m_i}$, while at $z^-_{i}$, $u_1$ is asymptotic to the periodic orbit
$\gamma_i^{-m_i}$.  Similarly, let
$\Gamma'_{1,-}=\br{w^-_j}$ ($\Gamma'_{2,+}=\br{w^+_j}$) be the set
of negative (positive) punctures of $v_1$ ($v_2$), and assume that at
$w^+_j$ ($w^-_j$) that $v_2$ ($v_1$) is asymptotic to the periodic orbit
$\hat \gamma_j^{n_j}$ ($\hat \gamma_j^{-n_j}$).

With this notation set, it follows from the definition of the generalized intersection
number that
\begin{equation}\label{e:gen-int-building}
\begin{aligned}
&\gin{u_{1}\odot u_{2}}{v_{1}\odot v_{2}}-\gin{u_{1}}{v_{1}}-\gin{u_{2}}{v_{2}} \\
&=\sum_{\substack{i, j \\ \gamma_i=\hat\gamma_j }}-m_{i}n_{j}
\bp{
\max\br{\tfrac{\alpha^\Phi(\gamma_i^{m_i})}{m_i},\tfrac{\alpha^\Phi(\hat\gamma_j^{n_j})}{n_j}}
+\max\br{\tfrac{\alpha^\Phi(\gamma_i^{-m_i})}{m_i},\tfrac{\alpha^\Phi(\hat\gamma_j^{-n_j})}{n_j}}
}
\end{aligned}
\end{equation}
Consider then the expression
\[
\Delta(\gamma, m, n):=-
mn
\bp{
\max\br{\tfrac{\alpha^\Phi(\gamma^{m})}{m},\tfrac{\alpha^\Phi(\gamma^{n})}{n}}
+\max\br{\tfrac{\alpha^\Phi(\gamma^{-m})}{m},\tfrac{\alpha^\Phi(\gamma^{-n})}{n}}
}
\]
with $m$ and $n$ positive integers.  If $\gamma$ is an even orbit, then \fullref{l:alpha-iterate} implies
that
\begin{align*}
\Delta(\gamma, m, n)&=-mn\bp{
\max\br{\tfrac{m\alpha^\Phi(\gamma)}{m},\tfrac{n\alpha^\Phi(\gamma)}{n}}
+\max\br{\tfrac{-m\alpha^\Phi(\gamma)}{m},\tfrac{-n\alpha^\Phi(\gamma)}{n}}
}\\
&=0.
\end{align*}
If $\gamma$ is odd and hyperbolic, we get
\begin{align*}
\Delta(\gamma, m, n)&=
-mn\Bigg(
\max\br{\tfrac{m\alpha^\Phi(\gamma)+\tfrac{1}{2}(m-p(m))}{m},\tfrac{n\alpha^\Phi(\gamma)+\tfrac{1}{2}(n-p(n))}{n}} \\
&\hskip.25in
+\max\br{\tfrac{-m\alpha^\Phi(\gamma)+\tfrac{1}{2}(-m-p(-m))}{m},\tfrac{-n\alpha^\Phi(\gamma)+\tfrac{1}{2}(-n-p(-n))}{n}}
\Bigg) \\
&=
-mn\bp{
\max\br{-\tfrac{p(m)}{m},-\tfrac{p(n)}{n}}
} \\
&=\min\br{n p(m), mp(n)}
\ge 0,
\end{align*}
and we get a strict inequality here if and only if $m$ and $n$ are both odd.
Finally, if $\gamma$ is elliptic, we get that there is an irrational $\theta$ so that
\begin{align*}
\Delta(\gamma, m, n)&=-mn\bp{
\max\br{\tfrac{\fl{m\theta}}{m}, \tfrac{\fl{n\theta}}{n}}
+\max\br{\tfrac{\fl{-m\theta}}{m}, \tfrac{\fl{-n\theta}}{n}}
} \\
&=\min\br{-n\fl{m\theta}, -m\fl{n\theta}}+\min\br{-n\fl{-m\theta}, -m\fl{-n\theta}} \\
&\ge-\fl{nm\theta}-\fl{-nm\theta} \\
&=\lceil nm\theta \rceil-\fl{nm\theta}=1.
\end{align*}
Thus for each pair $(i,j)$ with $\gamma_i=\hat\gamma_j$ and
$\gamma_i^{m_i}$ and $\hat\gamma_j^{n_j}$ both having odd Conley--Zehnder index,
we get a positive term in the sum \eqref{e:gen-int-building}, which proves our claim.
\end{proof}

We next prove \fullref{t:gin-positivity-intro}
which is one of the main motivations
for defining the generalized intersection number as we have.  It says that
for a pair of asymptotically cylindrical pseudoholomorphic curves with no common components,
the generalized intersection number is equal to the algebraic intersection number plus the total measure of
tangency at infinity between the two curves.  Since the generalized intersection number is a homotopy
invariant quantity, while the algebraic intersection number in general is not, this result demonstrates the utility of
this quantity in situations where one is hoping to 
obtain topological control over the total count of intersections between two curves.

Before restating and proving the theorem we establish some notation.
As previously noted, each puncture $z\in\Gamma$
of an asymptotically cylindrical pseudoholomorphic curve
$[\Sigma, j, \Gamma, u]\in\MW$ determines a pseudoholomorphic end $[u; z]$
(as defined in \fullref{ss:ends}) in one of the ends of $W$.  Given a second curve
$[\Sigma', j', \Gamma', v]$ having no common components with $\tl u$ and
so that the ends $[u; z]$ and $[v; w]$ are asymptotic to a covering of the same orbit with the same
sign, we can define the asymptotic intersection index
\[
\delta_{\infty}([u; z], [v; w])
\]
as before in \eqref{e:def-asymp-int-index}.
Otherwise, we define
\[
\delta_{\infty}([u; z], [v; w])=0,
\]
and we then define the \emph{total asymptotic intersection index} $\delta_{\infty}(u, v)$
of $[\Sigma, j, \Gamma, u]$ and $[\Sigma', j', \Gamma', v]$
by
\[
\delta_{\infty}(u, v):=\sum_{(z, w)\in\Gamma\times \Gamma'}\delta_{\infty}([u; z], [v; w]).
\]

We now restate and prove \fullref{t:gin-positivity-intro}.

\begin{theorem}[Positivity of the generalized intersection number]\label{t:gin-positivity}
Let $(W^{4}, \barJ)$ be an almost complex cobordism with cylindrical ends
$(\R^\pm\times M^\pm, J^\pm, \Ha^\pm)$, and let
$[\Sigma, j, \Gamma, u]$, $[\Sigma', j', \Gamma', v]\in\M(W, \barJ, \Ha^+, \Ha^-)$ be pseudoholomorphic curves in $W$ with no common components.
Then
\begin{equation}\label{e:gin-compute}
 \gin{u}{v}=\inum(u, v)+\delta_\infty(u, v).
\end{equation}
In particular
\[
\gin{u}{v}\ge\inum(u, v)\ge 0,
\]
and
\[
\gin{u}{v}=0
\]
if and only if $u$ and $v$ don't intersect, and the total asymptotic intersection number vanishes,
i.e.\
$\delta_\infty(u, v)=0$.
\end{theorem}

\begin{proof}
The equation \eqref{e:gin-compute} follows immediately from adding $\Omega^\Phi(u, v)$ to both sides
of \eqref{e:rel-int-compute} and then applying the definition \eqref{e:def-asymp-int-index} of the asymptotic intersection index.

The final two claims are direct consequences of \eqref{e:gin-compute},
local positivity of
intersections for pseudoholomorphic curves, and the nonnegativity
of the asymptotic intersection index from \eqref{e:ain-bound-alpha} and \eqref{e:def-asymp-int-index}.
\end{proof}

\subsection{The adjunction formula}\label{ss:adjunction}
\subsubsection{Relative Chern numbers and the total Conley--Zehnder index}\label{ss:chern}
In order to state the appropriate generalization of \eqref{e:adj-closed}, we will need to introduce
a relative version of the first Chern number of a complex vector bundle over a punctured Riemann surface.
This will then be used to define the total Conley--Zehnder index of an asymptotically cylindrical map.

Let $(\Sigma,j)$ be a closed Riemann surface, and let $\Gamma\subset\Sigma$ be a finite set.
We will call an open neighborhood $U$ of
$\Gamma$ disk-like if the closure $\bar U$ of $U$ is diffeomorphic to a disjoint union of closed disks,
and if each component of $U$ contains precisely one element of $\Gamma$.
Let $E\to\Sigma\setminus\Gamma$ be a complex vector bundle.
Since for any disk-like neighborhood $U$ of $\Gamma$ the set
$U\setminus\Gamma$ has the homotopy type of a disjoint union of
circles, the restriction $E|_{U\setminus\Gamma}$ is trivial.
Let $\Phi:U\setminus\Gamma\times\C^{n}\to E$ be a trivialization of $E|_{U\setminus\Gamma}$.
If $E$ is a line bundle, then we will define the first Chern number
 of $E$ relative to $\Phi$, denoted $c_{1}(E; \Phi)$ or $c_{1}^{\Phi}(E)$, to be equal to the algebraic count
of zeroes of a generic extension of the section
$\Phi(z)1$
of $E|_{U\setminus\Gamma}$ to a section of $E$.
If $\dim_{\C}E=n>1$, we define the first Chern number of $E$ relative to $\Phi$ 
by
\[
c_{1}(E; \Phi)=c_{1}(\Lambda_{\C}^{n}E; \Lambda_{\C}^{n}\Phi)
\]
where $\Lambda_{\C}^{n}E$ is the determinant line bundle of $E$
and $\Lambda_{\C}^{n}\Phi$ is the trivialization of
$\Lambda_{\C}^{n}E|_{U\setminus\Gamma}$ induced by $\Phi$.
We can define an equivalence relation on the set of pairs $(U,\Phi)$
of disk-like neighborhoods of $\Gamma$ and trivializations
$\Phi$ of $E|_{U\setminus\Gamma}$ be defining $(U,\Phi)\sim(V, \Psi)$
provided there exists a disk-like neighborhood
$U'\subset U\cap V$ of $\Gamma$ so that $\Phi|_{U'\setminus\Gamma}$
is homotopic to $\Psi|_{U'\setminus\Gamma}$.  The
relative first Chern number $c_{1}(E; \Phi)$ clearly only depends
on the equivalence class of the pair $(U,\Phi)$.

We collect some useful properties of the relative first Chern number in the lemma below,
but before stating the lemma we will establish some
notation and terminology.
Given a punctured Riemann surface $\Sigma\setminus\Gamma$,
we will call a pair $(U,\psi)$ a holomorphic cylindrical coordinate system 
around $\Gamma$ if $U$ is a disk-like neighborhood of $\Gamma$ and $\psi$ is a biholomorphic map
$\psi:\amalg_{i=1}^{\#\Gamma}(\c{0})_{i}\to\bar U\setminus\Gamma$.  
Using the coordinate fields $\ds$ and $\dt$ to identify $T_{(s,t)}(\c{0})$ with $\C$,
we note that the derivative of $\psi$ is a map
\[
d\psi:T(\amalg_{i=1}^{\#\Gamma}(\c{0})_{i})\approx (\amalg_{i=1}^{\#\Gamma}(\c{0})_{i})\times \C\to T\Sigma|_{U\setminus\Gamma}
\]
and therefore determines a trivialization of $(T\Sigma,j)$ over $U\setminus\Gamma$.
If $E\to\Sigma\setminus\Gamma$ is a complex bundle, and $\Phi$ and $\Psi$ are trivializations of
$E|_{U\setminus\Gamma}$, then
$\Phi^{-1}\circ\Psi$ determines a map $U\setminus\Gamma\to GL_{n}(\C)$.
If we denote $\psi_{i}=\psi|_{(\c{0})_{i}}$, then $\det((\Phi^{-1}\circ\Psi)(\psi_{i}(s',\cdot)))$ is an oriented
loop in $\C\setminus\{0\}$ for any fixed $s'\in\R^{+}$. We define $\deg(\Phi^{-1}\circ\Psi)$ by
\[
\deg(\Phi^{-1}\circ\Psi)=\sum_{i=1}^{\#\Gamma}\wind\left(\det\left((\Phi^{-1}\circ\Psi)(\psi_{i}(s',\cdot))\right)\right)
\]
and note that this definition does not depend on the choice of holomorphic cylindrical coordinates, or on $s'$.

\begin{lemma}\label{l:chern}
Let $(\Sigma,j)$ be a closed Riemann surface and let $\Gamma\subset\Sigma$ be a finite set.
Let $E$, $E_{1}$, $E_{2}\to\Sigma\setminus\Gamma$ be complex vector bundles over $\Sigma\setminus\Gamma$,
and let $U$ be a disk-like neighborhood of $\Gamma$.  Then
\begin{itemize}
\item If $\Phi$ and $\Psi$ are each trivializations of $E|_{U\setminus\Gamma}$ then
\begin{equation}\label{chern-triv}
c_{1}(E;\Phi)=c_{1}(E;\Psi)-\deg(\Phi^{-1}\circ\Psi).
\end{equation}
\item If $\Phi_{i}$ is a trivialization of $E_{i}|_{U\setminus\Gamma}$ for $i\in\{1,2\}$, then
\begin{equation}\label{chern-sum}
c_{1}(E_{1}\oplus E_{2}; \Phi_{1}\oplus\Phi_{2})=c_{1}(E_{1};\Phi_{1})+c_{1}(E_{2};\Phi_{2})
\end{equation}
and
\begin{equation}\label{chern-tensor}
c_{1}(E_{1}\otimes E_{2}; \Phi_{1}\otimes\Phi_{2})=c_{1}(E_{1};\Phi_{1})+c_{1}(E_{2};\Phi_{2}).
\end{equation}
\item If $(U,\psi)$ is a holomorphic cylindrical coordinate system, then
\begin{equation}\label{chern-euler}
c_{1}(T(\Sigma\setminus\Gamma); d\psi)=\chi(\Sigma\setminus\Gamma)=\chi(\Sigma)-\#\Gamma
\end{equation}
where $\chi(S)$ denotes the Euler characteristic of $S$.
\end{itemize}
\end{lemma}

\begin{proof}
To see that \eqref{chern-triv} is true, we pick a holomorphic cylindrical coordinate system $(U,\psi)$, and choose a generic section $h$ of
$E|_{U\setminus\Gamma}$ which is equal to $\Lambda_{\C}^{n}\Phi(z)1$ for $z\in\cup_{i}\psi_{i}([0,1]\times S^{1})$ and equal to
$\Lambda_{\C}^{n}\Psi(z)1$
for $z\in\cup_{i}\psi_{i}(\c{2})$.  Extending $h$ to a generic section of $\Lambda_{\C}^{n}E$, the algebraic count of zeroes of $h$ in 
$\Sigma\setminus\Gamma$ is equal to $c_{1}(E;\Psi)$ while the algebraic count of zeroes in
$\Sigma\setminus U$ is equal
to $c_{1}(E;\Phi)$.  Therefore, letting $m$ represent the algebraic count of zeroes of $h$ in the necks $\cup_{i}\psi_{i}([1,2]\times S^{1})$,
we must have
\[
c_{1}(E;\Phi)=c_{1}(E;\Psi)+m
\]
so it suffices to show $m=-\deg(\Phi^{-1}\circ\Psi)$.
This indeed follows from 
\begin{align*}
m&=\sum_{i=1}^{\#\Gamma}\wind((\Lambda_{\C}^{n}\Phi(\psi_{i}(2,\cdot)))^{-1}h(\psi_{i}(2,\cdot)))-\wind((\Lambda_{\C}^{n}\Phi(\psi_{i}(1,\cdot)))^{-1}h(\psi_{i}(1,\cdot))) \\
&=\sum_{i=1}^{\#\Gamma}\wind(1)-\wind((\Lambda_{\C}^{n}\Phi(\psi_{i}(1,\cdot)))^{-1}(\Lambda_{\C}^{n}\Psi(\psi_{i}(1,\cdot)))1) \\
&=\sum_{i=1}^{\#\Gamma}-\wind(\Lambda_{\C}^{n}(\Phi^{-1}\circ\Psi)(\psi_{i}(1,\cdot))1)=\sum_{i=1}^{\#\Gamma}-\wind(\det((\Phi^{-1}\circ\Psi)(\psi_{i}(1,\cdot)))) \\
&=-\deg(\Phi^{-1}\circ\Psi).
\end{align*}

To prove \eqref{chern-sum} and \eqref{chern-tensor}, we first observe that there
are natural isomorphisms
\[
\Lambda_{\C}^{n_{1}+n_{2}}(E_{1}\oplus E_{2})\to (\Lambda_{\C}^{n_{1}}E_{1})\otimes(\Lambda_{\C}^{n_{2}}E_{2})
\]
and
\[
\Lambda_{\C}^{n_{1}n_{2}}(E_{1}\otimes E_{2})\to (\Lambda_{\C}^{n_{1}}E_{1})\otimes(\Lambda_{\C}^{n_{2}}E_{2})
\]
where $n_{i}=\dim_{\C}E_{i}$.  Furthermore, the trivializations
$\Lambda_{\C}^{n_{1}+n_{2}}(\Phi_{1}\oplus\Phi_{2})$ of $\Lambda_{\C}^{n_{1}+n_{2}}(E_{1}\oplus E_{2})|_{U\setminus\Gamma}$
and
$\Lambda_{\C}^{n_{1}n_{2}}(\Phi_{1}\otimes\Phi_{2})$ of $\Lambda_{\C}^{n_{1}n_{2}}(E_{1}\otimes E_{2})|_{U\setminus\Gamma}$
both induce the trivialization
$(\Lambda_{\C}^{n_{1}}\Phi_{1})\otimes(\Lambda_{\C}^{n_{2}}\Phi_{2})$ on $(\Lambda_{\C}^{n_{1}}E_{1})\otimes(\Lambda_{\C}^{n_{2}}E_{2})|_{U\setminus\Gamma}$
under the given isomorphisms.  Therefore both \eqref{chern-sum} and \eqref{chern-tensor} are equivalent to the special case
of \eqref{chern-tensor} when $E_{1}$ and $E_{2}$ are line bundles.  In this case, we let $h_{i}$ be a generic section of $E_{i}$ agreeing with
$\Phi_{i}1$ over $U\setminus\Gamma$, we assume that the zero loci of $h_{1}$ and $h_{2}$ are disjoint.  Then $h_{1}\otimes h_{2}$ is 
a section of $E_{1}\otimes E_{2}$ agreeing with $(\Phi_{1}\otimes\Phi_{2})1$ over $U\setminus\Gamma$, and the algebraic count of zeroes
of $h_{1}\otimes h_{2}$ is given by the sum of the algebraic counts of zeroes of $h_{1}$ and $h_{2}$.

Finally, to see that \eqref{chern-euler} holds, we note that if $(U,\psi)$ is holomorphic cylindrical coordinate system then the section
$d\psi(z)1$ of $U\setminus\Gamma$ extends to a smooth section of $U$ with a simple zero of positive index at each point of $\Gamma$.  Any extension
of this section to a generic section of $T\Sigma$ will have $\chi(\Sigma)$ zeroes counted with sign,
of which $\chi(\Sigma)-\#\Gamma$ lie within $\Sigma\setminus\Gamma$.
\end{proof}

Having defined the relative first Chern number,
we now define the total Conley--Zehnder index.
Let $(\Sigma, j, \Gamma, u)\in\CW$ be an asymptotically cylindrical map
in a $4$-manifold $W$ with Hamiltonian structured ends, and assume that
$W$ is equipped with a compatible almost complex structure $\barJ$.
A choice of unitary trivialization of
$(\xi^{\Ha^\pm}, J^\pm)$ along each asymptotic limit of $u$
induces a complex trivialization of the pull-back bundle $(u^*TW, \barJ)$
near the punctures of $u$ since on the cylindrical ends we have
\[
(TW, \barJ)
\approx (\R\oplus\R X_{\Ha^{\pm}}\oplus\xi^{\Ha^\pm}, \tildeJ^{\pm})
\approx\C X_{\Ha^{\pm}}\oplus(\xi^{\Ha^\pm}, J^{\pm}).
\]
We will use $\Phi$ to denote both the chosen trivialization of $(\xi^{\Ha^\pm}, J^\pm)$ along the
asymptotic limits, and the induced trivialization of $(u^*TW, \barJ)$ near the punctures.
We then define the \emph{total Conley--Zehnder index} of $\mu(u)$ of
$u$ by
\[
\mu(u)=2 c_{1}^{\Phi}(u^{*}TW, u^{*}{\barJ})+\sum_{z\in\Gamma}\mu^{\Phi}(\gamma_{z}^{m_{z}}).
\]
where $u$ is assumed to be asymptotic to $\gamma_{z}^{m_{z}}$ at $z\in\Gamma$.
We note that this quantity is independent of choice of trivialization as a result of
\fullref{l:alpha-triv} and \fullref{l:chern}, and depends only on the homotopy classes
of $u$ and compatible $\barJ$.  We further note that in applications, the $4$-manifold $W$ is
typically equipped with a symplectic form with which $\barJ$ must be compatible
on all of $W$ (as opposed to being just compatible on the ends as we assume here).
In this case, it is well known
(see e.g.\ Proposition 5 and discussion following in Section 1.3 of \cite{hz})
that the space of such $\barJ$ is contractible, so
the total Conley--Zehnder index would depend only on the homotopy class of the map $u$.

\subsubsection{The adjunction formula}\label{sss:adjunction}
In this section we prove 
\fullref{t:adj-main-intro} which is
a generalization of \eqref{e:adj-closed} for punctured curves
which is stated in terms of the generalized intersection number.

Before proving the result, we review some basic
facts about the local intersection properties of pseudoholomorphic curves which we will need.
Proofs of these facts can be found in McDuff \cite{mcd94} or Micallef and White \cite{mw}.
Let $u:(\Sigma, j)\to (W^4, J)$ be a pseudoholomorphic map in an almost
complex $4$-manifold.  Since $du:T_z\Sigma\to T_{u(z)}W$
is complex linear, the derivative of $u$ is always either rank $2$ or rank $0$.
Define the singular set $S(u)$ of $u$ to be the set of points where the
derivative of $u$ is zero, i.e.\
\[
S(u):=\br{z\in\Sigma\,|\, du(z)=0}.
\]
Each point $z\in S(u)$ has an open neighborhood $U$ so that either
$u|_{U\setminus\br{z}}$ is an embedding, or $u|_{U\setminus\br{z}}$
factors through an embedding via a holomorphic covering.
Consequently, the singular set of a pseudoholomorphic map is an isolated subset of the domain.
For an asymptotically cylindrical pseudoholomorphic map,
this fact combined with the assumed convergence to an orbit cylinder at the punctures implies that the
singular set is finite.

Given a point $z\in S(u)$, assume that $u|_{U\setminus\br{z}}$ is an embedding for some
open neighborhood $U$ of $z$.
We can define a quantity $\delta(u; z)$ in the following way: there exists a
$C^1$-small perturbation $J'$ of $J$ supported near $u(z)$ and $C^1$-small perturbation $u'$ of
$u$ supported in a neighborhood $U$ of $z$ so that $u'|_{U}$ is an immersed $J'$-holomorphic map with
precisely $\delta(u; z)$ transverse double points.  This quantity is independent of the choice of perturbation.

Now, let $v:(\Sigma', j')\to (W^4, J)$ be another pseudoholomorphic map, and assume
that  $u(p)=v(q)$ for some $(p,q)\in\Sigma\times\Sigma'$.
Then there exist open neighborhoods $U$ of $p$ and $V$ of $q$ so that either
\[
V\subset v^{-1}(u(U))
\]
or
\[
u(U\setminus\br{p})\cap v(V\setminus\br{q})=\emptyset.
\]
In the latter case, when $u$ and $v$ have an isolated intersection at $u(p)=v(q)$, we will denote the local intersection
number of $u$ at $p$ with $v$ at $q$ by
\[
\delta((u; p), (v; q)).
\]
The local intersection number $\delta((u; p), (v; q))$ is always greater than or equal to $1$, and is
equal to $1$ if and only if $u$ is immersed at $p$, $v$ is immersed at $q$, and $u$ and $v$
intersect transversely.

Assume next that $u:(\Sigma, j)\to(W^4, J)$ is a simple map, that is that $u$ does not factor through
a branched cover of Riemann surfaces.  Define the double point set
$D_u\subset \Sigma\times\Sigma\setminus\Delta\Sigma$ of $u$ by
\[
D(u)=\br{(p,q)\,|\,u(p)=u(q), p\ne q}.
\]
As observed in \fullref{c:finite-sing}, the
preceding discussion along with \fullref{c:no-isect-puncture} and \fullref{c:factor-puncture} imply that an 
asymptotically cylindrical curve has at most a finite number of double points.
For such a map, we define the \emph{self-intersection index} $\delta(u)$ of $u$ by
\begin{equation}\label{e:self-int-index}
\delta(u)
=\sum_{z\in S(u)}\delta(u;z)
+\frac{1}{2}\sum_{(z,w)\in D(u)}\delta((u; z), (u; w)).
\end{equation}
We note that this quantity is an integer since the second sum counts
$(z, w)$ and $(w, z)$ separately.
A pseudoholomorphic map $u$ is an embedding if and only if $\delta(u)=0$.
Furthermore,  $\delta(u)$ has the following significance: if $J'$ is a compactly supported,
perturbation of $J$, and
$u'$ is an immersed, $J'$-holomorphic, compactly supported perturbation of $u$ having only 
transverse double points, then $u'$ has precisely $\delta(u)$ double points, that is
$\delta(u)=\frac{1}{2}\# D(u)$.

Now consider an almost complex cobordism of Hamiltonian structures
$(W, \barJ, \Ha^{+}, \Ha^{-})$
and let $[\Sigma, j, \Gamma, u]\in\MW$ be a
simple pseudoholomorphic curve.
We define the \emph{total asymptotic self-intersection index} of $u$ by
\[
\delta_{\infty}(u)
=\sum_{z\in\Gamma}\delta_{\infty}(u;z)
+\frac{1}{2}
\sum_{\substack{(z,w)\in\Gamma\times\Gamma \\ z\ne w}}\delta_{\infty}([u; z], [u; w]),
\]
where $\delta_{\infty}(u; z)$ is as defined in \eqref{e:sint-index} and
$\delta_{\infty}([u; z], [u; w])$ is as defined in \eqref{e:def-asymp-int-index}.
This quantity will be nonnegative integer-valued since the asymptotic intersection index
$\delta_{\infty}([u; z], [u; w])$ is symmetric, and each pair of distinct punctures is counted twice
in the second sum.

Letting $[\Sigma, j, \Gamma, u]\in\MW$ still denote a simple curve,
we now define the \emph{singularity index}\footnote{
	N.B.\ This definition of $\sing(u)$ is different from the one used in \cite{sief05}
	}
$\sing(u)$ of $u$ by
\[
\sing(u):=\delta(u)+\delta_{\infty}(u).
\]
The singularity index is nonnegative-integer-valued, and equals zero for a given curve if and only if
that curve is embedded and has total asymptotic self-intersection index equal to zero.
Finally,
recalling the definition \eqref{e:spec-cover-end} of the spectral covering number
$\cover(u; z)$ of a pseudoholomorphic end $[u; z]$,
we define the \emph{total spectral covering number} of $u$ by
\[
\cover(u)=\sum_{z\in\Gamma}\cover(u; z).
\]
The spectral covering number depends only on the asymptotic data of the map, and not on the map itself.

We now prove \fullref{t:adj-main-intro} which we restate here for the convenience
of the reader (cf.\ Remark 3.2 in \cite{hut02} and Theorem 3.13 in \cite{sief05}).

\begin{theorem}[\fullref{t:adj-main-intro}]\label{t:adj-main}
Let $(W, \barJ)$ be an almost complex $4$-manifold with cylindrical ends
$(\R^{\pm}\times M^{\pm}, \Ha^{\pm}, J^{\pm})$
and let
$[\Sigma, j, \Gamma, u]\in\MW$ be a connected, simple pseudoholomorphic curve.
Then
\begin{equation}\label{e:adj-main}
\gin{u}{u}-\tfrac{1}{2}\mu(u)+\tfrac{1}{2}\#\Gamma_{\text{odd}}
+\chi(\Sigma)-\cover(u)=2\sing(u)
\end{equation}
where $\#\Gamma_{\text{odd}}$ denotes the number of punctures of $u$
limiting to orbits with odd Conley--Zehnder indices.
In particular, 
\[
\gin{u}{u}-\tfrac{1}{2}\mu(u)+\tfrac{1}{2}\#\Gamma_{\text{odd}}
+\chi(\Sigma)-\cover(u)\ge 0
\]
and equality occurs if and only if $u$ is an embedding and $u$ has
total asymptotic self-intersection index equal to zero.
\end{theorem}

\begin{proof}
We first observe that it suffices to assume that $u$ is an immersion having only transverse double points.
If not we can apply results of Micallef and White \cite{mw} to find an immersion
$u_{1}:\Sigma\setminus\Gamma\to W$
which is homotopic to $u$ via a homotopy that is constant outside a compact subset of
$\Sigma\setminus\Gamma$,
which has precisely
$\delta(u)$ transverse double points and no other singularities, and
which is $\tildeJ_{1}$-holomorphic for some $\tildeJ_{1}$ which is homotopic to $\tildeJ$ via a homotopy
that is constant outside a compact set.  Such a homotopy leaves all the terms
in \eqref{e:adj-main} constant, and it thus suffices to prove the result for such a $u_{1}$.

Proceeding with these assumptions, we compute the relative intersection number of $u$ with itself for a given
trivialization $\Phi$ by pushing $u$ off of itself with a section of the normal bundle to $u$
which is asymptotic at each puncture to $\Phi(\ep)$.
Using standard homotopy invariance and transversality arguments,
it can be shown that the zeroes of the section contribute $c_1^\Phi(N_{u})$ to $\rin^\Phi(u, u)$,
and the double points of $u$ contribute $2\delta(u)$ to $\rin^{\Phi}(u, u)$.
Moreover, the local asymptotic analysis in \fullref{ss:asympt-int} shows that multiply covered ends,
and distinct ends approaching coverings of the same orbits contributes a total of
\[
\sum_{z\in\Gamma}\ain^{\Phi}(u; z)
+\sum_{\substack{(z,w)\in\Gamma\times\Gamma \\ z\ne w}}\ain^{\Phi}([u; z], [u; w])
\]
so consequently we have that
\begin{equation}\label{e:rel-adj}
\rin^{\Phi}(u, u)=c_{1}^{\Phi}(N_{u})+2\delta(u)
+\sum_{z\in\Gamma}\ain^{\Phi}(u; z)
+\sum_{\substack{(z,w)\in\Gamma\times\Gamma \\ z\ne w}}\ain^{\Phi}([u; z], [u; w]).
\end{equation}
(This is equivalent to formula (4.1) in \cite{hut09} and to the first formula on page 110 of \cite{sief05}).

Assuming that at a puncture $z\in\Gamma^\pm$, $u$ is asymptotic to $\gamma_z^{m_z}$, we add
\begin{align*}
&\sum_{\substack{(z,w)\in\Gamma\times\Gamma \\ \gamma_z=\gamma_w \\ m_z m_w>0 }}
m_{z}m_{w}\max\br{\tfrac{\alpha^\Phi(\gamma_z^{m_z})}{|m_z|},
\tfrac{\alpha^\Phi(\gamma_w^{m_w})}{|m_w|}} \\
&\hskip1in=\sum_{z\in\Gamma}|m_z|\alpha^\Phi(\gamma_z^{m_z})
+\sum_{\substack{(z,w)\in\Gamma\times\Gamma \\ z\ne w \\ \gamma_z=\gamma_w \\ m_z m_w>0 }}
m_{z}m_{w}\max\br{\tfrac{\alpha^\Phi(\gamma_z^{m_z})}{|m_z|},
\tfrac{\alpha^\Phi(\gamma_w^{m_w})}{|m_w|}}
\end{align*}
to \eqref{e:rel-adj}
and use the definitions of $\gin{\cdot}{\cdot}$,  $\delta_{\infty}(\cdot)$, and
$\delta_{\infty}(\cdot, \cdot)$
to get
\begin{align*}
\gin{u}{u}&=c_{1}^{\Phi}(N_{u})+2\delta(u)
+\sum_{z\in\Gamma}[2\delta_{\infty}(u; z)+\alpha(u; z)+\cover(u; z)-1] \\
&\hskip.5in
+\sum_{\substack{(z,w)\in\Gamma\times\Gamma \\ z\ne w}}\delta_{\infty}([u; z], [u; w]) \\
&=c_{1}^{\Phi}(N_{u})+2\sing(u)+\cover(u)-\#\Gamma+
\sum_{z\in\Gamma}\alpha(u; z)\\
&=c_{1}^{\Phi}(N_{u})+2\sing(u)+\cover(u)-\#\Gamma
-\tfrac{1}{2}\#\Gamma_{\text{odd}}+\sum_{z\in\Gamma}\tfrac{1}{2}\mu^\Phi(u; z).
\end{align*}

Finally, using the properties of the relative Chern number from \fullref{l:chern}, we have that
\begin{align*}
c_{1}^{\Phi}(u^{*}TW)&=\chi(\Sigma\setminus\Gamma)+c_{1}^{\Phi}(N_{u}) \\
&=\chi(\Sigma)-\#\Gamma+c_{1}^{\Phi}(N_{u})
\end{align*}
which, combined with the above formula for $\gin{u}{u}$, gives that
\[
\gin{u}{u}-\tfrac{1}{2}\mu(u)+\chi(\Sigma)
+\tfrac{1}{2}\#\Gamma_{\text{odd}}-\cover(u)=2\sing(u).
\]
as claimed.
\end{proof}

It is useful in certain contexts to rewrite \eqref{e:adj-main} in terms of the virtual dimension
(or Fredholm index) of the moduli space of the curve being considered.
We define the index of a pseudoholomorphic curve
$[\Sigma, j, \Gamma, u]\in\MW$ by the formula
\begin{equation}
\ind(u)=\mu(u)-\chi(\Sigma)+\#\Gamma.
\end{equation}
If $u$ is simple, the index is equal to the dimension of the moduli space of unparametrized
curves near $u$ for a generic choice of admissible almost complex structures on $W$
(see \cite{hwz:prop3, dr}).
Using this formula we can restate the adjunction formula as follows.

\begin{corollary}
Let $(W, \barJ)$ be an almost complex $4$-manifold with cylindrical ends
$(\R^{\pm}\times M^{\pm}, \Ha^{\pm}, J^{\pm})$
and let
$[\Sigma, j, \Gamma, u]\in\MW$ be a connected, simple pseudoholomorphic curve.
Then
\begin{equation}\label{e:adj-index}
\gin{u}{u}-\tfrac{1}{2}\bbr{\ind(u)-\chi(\Sigma)+\#\Gamma_{\text{even}}}
+\#\Gamma-\cover(u)=2\sing(u).
\end{equation}
\end{corollary}

Using the homotopy invariance of the generalized intersection number, we have the following corollary
which allows us to bound the number of intersections between two pseudoholomorphic
curves which are homotopic.  This result is useful in the study of finite energy foliations
(see e.g. \cite{hwz:fol}).

\begin{corollary}
Let
$[\Sigma, j, \Gamma, u]\in\MW$ be a connected, simple pseudoholomorphic curve,
and let
$[\Sigma', j', \Gamma', v]\in\MW$ be
a second curve, distinct from $u$, which is homotopic
to $u$ in $C^{\infty}(W, \Ha^{+}, \Ha^{-})$.
Then
\[
\inum(u, v)\le \tfrac{1}{2}\bbr{\ind(u)-\chi(\Sigma)+\#\Gamma_{\text{even}}}
+\cover(u)-\#\Gamma+2\sing(u).
\]
\end{corollary}

\section{Cylindrical cobordisms}\label{s:cylindrical}
In this section we will consider some applications of the intersection theory we've developed
in previous sections to curves in a cylindrical cobordism.
Throughout the remainder of this section $M$ will denote a $3$-manifold equipped with the nondegenerate
stable Hamiltonian structure $\Ha=(\lambda, \omega)$, and a compatible almost complex structure
$J\in\J\MH$.
As in previous sections, we will denote the associated $\R$-invariant almost complex structure on
$\R\times M$ by $\tildeJ$.

It will be convenient for many of the results in this section to consider a special class of smooth
asymptotically cylindrical maps.
Given an asymptotically cylindrical map $(\Sigma, j, \Gamma, (a,u))\in\CM$, we say that
$(\Sigma, j, \Gamma, (a,u))$ has \emph{ends which wind} if there is an open neighborhood
$U$ of $\Gamma$ so that $u|_{U\setminus\Gamma}$ does not intersect any of the asymptotic limits of $u$,
i.e.\ in the language of \fullref{ss:asymp-winding-nos} a neighborhood of each puncture
determines an asymptotically cylindrical end which winds.
If $(\Sigma, j, \Gamma, (a,u))$ is an asymptotically cylindrical map with ends that wind,
we can define for each puncture $z\in\Gamma$ the asymptotic winding $\winfty^{\Phi}(u; z)$ as we do
in \fullref{ss:asymp-winding-nos}
by choosing a small clockwise loop $\rho$ around $z$ and
computing the winding of $u\circ\rho$ relative to a trivialization $\Phi$ of
$\xh$ along the asymptotic limit of $u$ at $z$.

\subsection{Tangency to the flow and $\wind_\pi$}
Let $\tl u=(a, u):\Sigma\setminus\Gamma\to\R\times M$ be a pseudoholomorphic curve, and
let $\pi_{\xh}:TM=\R\X\oplus\xh\to\xh$ be the projection onto $\xh$ determined by the splitting \eqref{e:splitting}.
Then the projection $\pi_{\xh}\circ du$ of the derivative of $u$ is a complex linear map
\[
\pi_{\xh}\circ du(z):(T_z(\Sigma\setminus\Gamma), j)\to(\xh_{u(z)}, J)
\]
according to the second equation in \eqref{e:j-hol-cylindrical}.
Thus for any $z\in\Sigma\setminus\Gamma$,  $\pi_{\xh}\circ du(z)$ has rank $0$ or $2$.
In the case that the Hamiltonian structure on $M$ arises from a contact form, it is shown in
\cite{hwz:prop2} that 
$\pi_{\xh}\circ du$ either vanishes identically, or has isolated zeroes
of finite positive order.
The proof given there readily generalizes to the case of an arbitrary stable Hamiltonian structure
(see \fullref{a:windpi}).  Moreover, assuming that $\pi_{\xh}\circ du$ doesn't vanish identically,
it follows from \fullref{c:pi-du-punctures} that $\pi_{\xh}\circ du$ has
at most a finite number of zeroes.  Following \cite{hwz:prop2}, we will
denote the algebraic count of the zeroes of $\pi_{\xh}\circ du$ by $\wind_\pi(u)$.

Since $\pi_{\xh}\circ du$ is a section of the complex line bundle,
$\Hom_\C(T_z(\Sigma\setminus\Gamma), u^*\xh)$,
the algebraic count of zeroes of $\pi_{\xh}\circ du$ is determined by the behavior of $u$ near the punctures.
Indeed in \cite{hwz:prop2} it is shown that
\begin{equation}\label{e:wind-pi-ind}
\wind_{\pi}(u)+d_{0}(\tl u)=\tfrac{1}{2}(\ind(\tl u)-\chi(\Sigma)+\#\Gamma_{\text{even}})
\end{equation}
where
\[
d_{0}(\tl u)=\sum_{z\in\Gamma}d_{0}(\tl u; z), 
\]
and $d_0(\tl u; z)$ is as defined in \eqref{e:def}.
Substituting into \eqref{e:adj-index} gives the following corollary to \fullref{t:adj-main} which will be useful
for some results later in this section.

\begin{corollary}\label{c:adj-wind-pi}
Let $[\Sigma, j, \Gamma, \tl u=(a, u)]\in\MM$ be a connected, simple
pseudoholomorphic curve.
If $\pi_{\xh}\circ du$ doesn't vanish identically, then 
\begin{equation}\label{e:adj-wind-pi}
\gin{\tl u}{\tl u}-\bbr{\wind_{\pi}(\tl u)+d_{0}(\tl u)}
+\#\Gamma-\cover(\tl u)=2\sing(\tl u).
\end{equation}
\end{corollary}

The results in this section will mostly be concerned with pseudoholomorphic curves that have no component
with image contained in an orbit cylinder.  We collect in the following lemma some conditions that
are equivalent to this.  This fact is well-known and can be deduced using
\fullref{c:pi-du-punctures} and \fullref{c:pi-du}.
For further discussion and proof see \cite{hwz:prop2}.

\begin{lemma}
Consider a pseudoholomorphic curve $[\Sigma, j, \Gamma, \tl u=(a, u)]\in\MM$.
The following are equivalent.
\begin{enumerate}
\item $\pi_{\xh}\circ du$ does not vanish identically on any component of $\Sigma\setminus\Gamma$.
\item No component of the map has image contained in an orbit cylinder.
\item No component of the map is fixed by the $\R$-action.
\item The map $[\Sigma, j, \Gamma, \tl u=(a, u)]$ has ends which wind, i.e.\
for each $z\in\Gamma$ and trivialization $\Phi$ of $\xh$ along the asymptotic limit of $\tl u$ at $z$,
the asymptotic winding $\winfty^{\Phi}(\tl u; z)$ is defined.
\end{enumerate}
\end{lemma}

\subsection{Intersections with orbit cylinders}
In this section we compute intersection numbers
of an orbit cylinder with an asymptotically cylindrical map with ends that wind.
Throughout this section, we will be dealing with an asymptotically cylindrical map
$(\Sigma, j, \Gamma, \tl u)$ with ends that wind, and we will use the notation
$\gamma_{z}^{m_{z}}$, with $m_{z}$ a nonzero integer and $\gamma_{z}\in\Pe^{0}\MH$
a simple periodic orbit,
to indicate that at the puncture $z\in\Gamma$, $\tl u$ is asymptotic to $\tl \gamma_{z}^{m_{z}}$.

We first compute the relative intersection number.

\begin{lemma}\label{l:rin-cylinder}
Let $(\Sigma, j, \Gamma, \tl u=(a,u))\in\CM$
be an asymptotically cylindrical map with ends that wind, and let
$\tl\gamma^k$
be the orbit cylinder over $\gamma^{k}$,
with $\gamma$ a simple periodic orbit and $k$ a positive integer.
Then the algebraic intersection number $\inum(\gamma, u)$ is well-defined and
\begin{equation}\label{e:rin-cylinder}
\rin^\Phi(\tl u, \tl\gamma^k)=
k\bp{
\inum(\gamma, u)-\sum_{\substack{z\in\Gamma \\ \gamma_z=\gamma}}
\winfty^\Phi(u; z)
}.
\end{equation}
\end{lemma}

\begin{proof}
Since we assume that the map $\tl u$ has ends that wind,
there is a neighborhood of each puncture on which $\tl u$ doesn't intersect
$\tl \gamma$ and hence on which the projected map $u$ doesn't intersect $\gamma$.
Applying standard transversality and homotopy-invariance arguments,
$u$ can therefore be perturbed on a compact set to a map with a finite number of isolated,
transverse intersections with $\gamma$, and the total algebraic count of these intersections
is independent of the perturbation.  The intersection number $\inum(\gamma, u)$ can then be defined
to be the algebraic count of these intersections.

Using the definition of the relative intersection number, and using the same reasoning that leads
to \fullref{t:rin-intasymp}
we have
\[
\rin^\Phi(\tl u, \tl\gamma^k)=\inum(\tl \gamma^k, \tl u)
+\sum_{\substack{z\in\Gamma^+ \\ \gamma_z=\gamma}}\ain^\Phi([\tl\gamma^k; \infty], [\tl u; z])
+\sum_{\substack{z\in\Gamma^- \\ \gamma_z=\gamma}}\ain^\Phi([\tl\gamma^k; 0], [\tl u; z]).
\]
where we consider $\tl\gamma^k$ as a map from $S^2\setminus\br{0, \infty}$ with $\infty$
the positive puncture and $0$ the negative puncture.  
Since $\tl\gamma$ is fixed by the $\R$-action
we have that\footnote{
Note that the orientation of $M$ matters in the right-most expression since $M$ is $3$-dimensional.
We always orient $M$ so that $\lambda\wedge\omega$ is a positive volume form, and
$\R\times M$ is oriented so that $dx\wedge\lambda\wedge\omega>0$ where
$x$ is the coordinate along $\R$.
}
\begin{equation}\label{e:int-cylinder}
\inum(\tl\gamma^k, \tl u)=\inum(\gamma^k, u)=k\inum(\gamma, u).
\end{equation}
Moreover, it follows from \fullref{l:rain-compute} and the definitions of
$\wrel$ and $\winfty$ that for a puncture $z\in\Gamma^+$ with $\gamma_z=\gamma$ we have
\begin{align*}
\ain^\Phi([\tl\gamma^k; \infty], [\tl u; z])
&=-\tfrac{1}{k m_z}\wrel^\Phi(m_z\cdot[\tl\gamma^k; \infty], k\cdot[\tl u; z]) \\
&=-\tfrac{1}{k m_z} k m_z \winfty^\Phi(k\cdot [\tl u; z]) \\
&=-k\winfty^\Phi(\tl u; z).
\end{align*}
Similarly, for $z\in\Gamma^-$ with $\gamma_z=\gamma$ we have
\[
\ain^\Phi([\tl\gamma^k; 0], [\tl u; z])
=-k\winfty^\Phi(\tl u; z)
\]
which completes the proof.
\end{proof}

Next we compute the generalized intersection number of an orbit cylinder
with an asymptotically cylindrical map with ends that wind.
A variation on this is used in \cite{mom10-2} 
to study a variation on contact homology 
in which the curves in the differential are required to not intersect
a prescribed collection of elliptic orbits.

\begin{corollary}\label{c:int-cylinder}
Let $(\Sigma, j, \Gamma, \tl u=(a,u))\in\CM$
be an asymptotically cylindrical map with ends that wind
and let $\tl\gamma$ be a cylinder over a simply covered periodic orbit $\gamma$.
Then for $k>0$ we have
\begin{align*}
\gin{\tl u}{\tl \gamma^{k}}
&=
k\Big(\inum\bp{\gamma, u} \\
     &\hskip.25in
     +\sum_{\substack{z\in\Gamma_{+} \\ \gamma_{z}=\gamma}}
     m_{z}\left[\max\br{\tfrac{\alpha(\gamma^{m_z})}{m_z}, \tfrac{\alpha(\gamma^{k})}{k}}
     -\tfrac{\alpha(\gamma^{m_z})}{m_z}\right]
     +d_{0}(\tl u; z) \\
     &\hskip.25in
     +\sum_{\substack{z\in\Gamma_{-} \\ \gamma_z=\gamma }}
     |m_z|\left[\max\br{\tfrac{\alpha(\gamma^{m_z})}{|m_z|}, \tfrac{\alpha(\gamma^{-k})}{k}}
     -\tfrac{\alpha(\gamma^{m_z})}{|m_z|}\right]
     +d_{0}(\tl u; z)\Big).
\end{align*}
\end{corollary}

\begin{proof}
The result follows from adding
\[
\sum_{\substack{z\in\Gamma_{+} \\ \gamma_{z}=\gamma}}
     km_{z}\max\br{\tfrac{\alpha(\gamma^{m_z})}{m_z}, \tfrac{\alpha(\gamma^{k})}{k}}
     +\sum_{\substack{z\in\Gamma_{-} \\ \gamma_z=\gamma }}
     k|m_z| \max\br{\tfrac{\alpha(\gamma^{m_z})}{|m_z|}, \tfrac{\alpha(\gamma^{-k})}{k}}
\]
to \eqref{e:rin-cylinder} and using the definition of $d_0(\tl u; z)$ and $\gin{\tl u}{\tl \gamma^k}$.
\end{proof}

We observe that in the event that the map $(\Sigma, j, \Gamma, \tl u)$ in this Corollary is pseudoholomorphic,
then the formula given above expresses the generalized intersection number
$\gin{\tl u}{\tl\gamma^{k}}$ as a sum of nonnegative terms.  Indeed, the terms of the form
$\max\br{\frac{\alpha(\gamma^{m})}{|m|}, \frac{\alpha(\gamma^{k})}{|k|}}-\frac{\alpha(\gamma^{m})}{|m|}$ are
nonnegative for any map.  In the case that $\tl u$ is pseudoholomorphic, the
terms $d_{0}(\tl u; z)$ are nonnegative by \fullref{l:wind-infinity}.
Moreover, when $\tl u$ is pseudoholomorphic the fact that $\inum(\gamma, u)$ is nonnegative
follows from the fact from \eqref{e:int-cylinder} that $\inum(\gamma, u)=\inum(\tl \gamma, \tl u)$ 
and positivity of intersections for pseudoholomorphic curves.

In the case that $\gamma$ is a hyperbolic orbit, the iteration formulas from
\fullref{l:alpha-iterate} simplify the formula 
from \fullref{c:int-cylinder} somewhat.  We state this result here and omit the straightforward proof.

\begin{corollary}
Let $\tl u=(a,u)$ and $\gamma$ satisfy the assumptions of \fullref{c:int-cylinder}, and let $k>0$.
If $\gamma$ is an even orbit, then
\[
\gin{\tl u}{\tl \gamma^{k}}=
k\bp{\inum\bp{\gamma, u}
     +\sum_{\substack{z\in\Gamma \\ \gamma_{z}=\gamma}}d_{0}(\tl u; z)}.
\]
If $\gamma$ is an odd, hyperbolic orbit, and $k$ is even, then
\begin{align*}
\gin{\tl u}{\tl \gamma^{k}}=
k\bp{\inum\bp{\gamma, u} 
     +\sum_{\substack{z\in\Gamma \\ \gamma_{z}=\gamma}}
    \bbr{\tfrac{p(m_z)}{2}
     +d_{0}(\tl u; z)}}
\end{align*}
where $p(m)$ denotes the parity of $m$.
If $\gamma$ is an odd, hyperbolic orbit, and $k$ is odd, then
\begin{align*}
\gin{\tl u}{\tl \gamma^{k}}
=
k\inum\bp{\gamma, u}
     +\sum_{\substack{z\in\Gamma \\ \gamma_{z}=\gamma}}\bbr{
     \tfrac{p(m_z)}{2}\bp{k-\min\br{ k, |m_z|}
     }
     +kd_{0}(\tl u; z)}.
\end{align*}
\end{corollary}

Finally, we compute the generalized intersection number of an orbit cylinder with itself.

\begin{proposition}
Let $\gamma$ be a simple periodic orbit, and let $k>0$.  Then
\[
\gin{\tl\gamma^k}{\tl\gamma^k}=-kp(\gamma^{k})
\]
where $p(\gamma^k)$ is the parity of periodic orbit $\gamma^k$.
\end{proposition}

\begin{proof}
It follows immediately from the definition of the relative intersection number that
$\rin^\Phi(\tl\gamma^k, \tl\gamma^k)$ vanishes for any trivialization $\Phi$ of $\gamma^*\xh$.
Consequently, the definition \eqref{e:gin-def} of the generalized intersection number gives us
that
\begin{align*}
\gin{\tl\gamma^k}{\tl\gamma^k}&=\Omega^\Phi(\tl\gamma^k, \tl\gamma^k) \\
&=k\bp{\alpha^\Phi(\gamma^k)+\alpha^\Phi(\gamma^{-k})}.
\end{align*}
If $\gamma$ is an even orbit, we have from \eqref{e:alpha-even} that
\[
\alpha^\Phi(\gamma^k)+\alpha^\Phi(\gamma^{-k})=
k\alpha^\Phi(\gamma)-k\alpha^\Phi(\gamma)=0
\]
so $\gin{\tl\gamma^k}{\tl\gamma^k}=0=-kp(\gamma^k)$ in this case.
If $\gamma$ is an odd hyperbolic orbit, we have from \eqref{e:alpha-oddhyperbolic} that
\begin{align*}
\alpha^\Phi(\gamma^k)+\alpha^\Phi(\gamma^{-k})
&=k\alpha^\Phi(\gamma)+\tfrac{k-p(k)}{2}-k\alpha^\Phi(\gamma)+\tfrac{-k-p(-k)}{2} \\
&=-p(k)=-p(\gamma^k),
\end{align*}
so $\gin{\tl\gamma^k}{\tl\gamma^k}=-kp(\gamma^k)$ in this case as well.
Finally, if $\gamma$ is elliptic, it follows from \eqref{e:alpha-elliptic}
that there exists an irrational number $\theta$ so that
\begin{align*}
\alpha^\Phi(\gamma^k)+\alpha^\Phi(\gamma^{-k})
&=\fl{k\theta}+\fl{-k\theta} \\
&=\fl{k\theta}-\ceil{k\theta}=-1 \\
\end{align*}
so we have $\gin{\tl\gamma^k}{\tl\gamma^k}=-k=-kp(\gamma^k)$ for elliptic orbits as well.
\end{proof}

This computation has an application to curves in $\R\times M$ equipped
with a non-$\R$-invariant almost complex structure.  We consider a stable Hamiltonian structure
$\Ha=(\lambda, \omega)$ on $M$, and let
$c\cdot\Ha$ denote the stable Hamiltonian structure $c\cdot\Ha=(c\cdot\lambda, c\cdot\omega)$.
Consider $\R\times M$ equipped with an almost complex structure
$\barJ$ which agrees on $[1, +\infty)\times M$
with a standard cylindrical $\tildeJ_1$ for $c_1\cdot\Ha$
and agrees on $(-\infty, -1]\times M$ with a standard cylindrical $\tildeJ_2$ for $c_2\cdot\Ha$.
If $\gamma^k$ is a periodic orbit for the structure $\Ha$, it is also a periodic orbit for the structure
$c_i\cdot\Ha$, and the Conley--Zehnder index of $\gamma^k$ is the same relative to all three structures.
We can consider the orbit cylinder $\tl\gamma^k$, which in general will not be $\barJ$-holomorphic,
but is still a smooth asymptotically cylindrical map.  The preceding proposition still applies and implies that
$\tl\gamma^k$ has a negative (generalized) self-intersection number if $\gamma^k$ is an odd orbit.
By the positivity of intersections for the generalized intersection number (\fullref{t:gin-positivity}), and
the homotopy invariance of the generalized intersection number, we can conclude that
there is at most one $\barJ$-holomorphic curve in the homotopy class of the orbit cylinder $\tl\gamma^k$.

\subsection{Counting intersections}
In this section we will compute generalized intersection number of two maps in terms of
data that is entirely determined by the $M$-components of the maps:
namely the intersection numbers of each map with the asymptotic limits of the other,
the asymptotic winding invariants, and the asymptotic data.
We then identify necessary and sufficient conditions for the generalized intersection number of
two pseudoholomorphic curves to be zero.
Finally, we examine a set of necessary and sufficient conditions for the projections of two
pseudoholomorphic curves into the three manifold to not intersect.

We recall that there is an $\R$-action on maps in cylindrical cobordisms, defined by shifting
in the $\R$-direction of $\R\times M$.
We will denote this action by $c\cdot\tl u=(a+c, u)$
or by $\tl u_{c}=(a+c, u)$ where
$\tl u=(a,u):\Sigma\setminus\Gamma\to\R\times M$ is an asymptotically cylindrical map.

We also recall that \cite{behwz}
defines what it means for a pseudoholomorphic curve to converge to a pseudoholomorphic
building (see CHC1-CHC2 in Section 7.3 of \cite{behwz}).
The definition of convergence of pseudoholomorphic curves given there can be generalized in
straightforward way to asymptotically cylindrical maps converging to
asymptotically cylindrical buildings.

\begin{lemma}\label{l:rin-sft-compactness}
Let $(\Sigma_{k}, j_{k}, \Gamma_{k}, \tl u_{k})$, $(\Sigma'_{k}, j'_{k}, \Gamma'_{k}, \tl v_{k})\in \CM$
be sequences of asymptotically cylindrical maps
so that the disjoint union $\tl u_{k}+\tl v_{k}$ converges in the sense of \cite{behwz} to
a $2$-level building $(\tl u_{\infty, 1}\odot\tl u_{\infty, 2})+(\tl v_{\infty, 1}\odot\tl v_{\infty, 2})$.
Then
\[
\lim_{k\to\infty}\rin^\Phi(\tl u_{k},\tl v_{k})
=\rin^\Phi(\tl u_{\infty, 1},\tl v_{\infty, 1})+\rin^\Phi(\tl u_{\infty, 2},\tl v_{\infty, 2}).
\]
\end{lemma}

We remark that the result and proof here easily adapt to the case where the limit building has more than $2$
levels and when the limiting curves are mapped into a splitting symplectic manifold.

\begin{proof}
According to definition of convergence in \cite{behwz}, our assumptions imply that there exist
constants $c_k$ and $d_k$ 
and holomorphic reparametrizations
$\phi_{k, i}:\Sigma_{k}\amalg \Sigma'_{k}\to \Sigma_{k}\amalg \Sigma'_{k}$
so that
$c_k\cdot(\tl u_k +\tl v_k)\circ \phi_{k, 1}$ converges in $C^\infty_{loc}$ to
$\tl u_{\infty, 1} + \tl v_{\infty, 1}$, and so that
$d_k\cdot(\tl u_k +\tl v_k)\circ \phi_{k, 2}$ converges in $C^\infty_{loc}$ to
$\tl u_{\infty, 2}+ \tl v_{\infty, 2}$.  Moreover,
we can compactify $(\R\times M)\odot(\R\times M)$ and identify it with
$[0, 1]\times M$ (equipped with a $C^{0}$ structure)
and can choose a sequence of identifications of the 
compactification of $\R\times M$ with $[0, 1]\times M$
and diffeomorphisms $\psi_{k}:\Sigma_{k}\amalg \Sigma'_{k}\to \Sigma_{k}\amalg \Sigma'_{k}$
in such a way that
$(\tl u_k +\tl v_k)\circ\psi_{k}$ converges uniformly to
$(\tl u_{\infty, 1}\odot\tl u_{\infty, 2})+(\tl v_{\infty, 1}\odot\tl v_{\infty, 2})$,
all viewed as maps into $[0, 1]\times M$.
It then follows from the properties of intersection numbers, and the definition
of the relative intersection number that for large $k$ we have
$\rin^\Phi(\tl u_{k},\tl v_{k})=\rin^\Phi(\tl u_{\infty, 1},\tl v_{\infty, 1})+\rin^\Phi(\tl u_{\infty, 2},\tl v_{\infty, 2})$
as claimed.
\end{proof}

As an immediate application of this, we have the following computation of
the generalized intersection number of two asymptotically cylindrical maps with ends that wind.
\begin{theorem}\label{t:gin-cylindrical}
Let 
$(\Sigma, j, \Gamma, \tl u=(a, u))$, $(\Sigma', j', \Gamma', \tl v=(b, v))\in\CM$
be asymptotically cylindrical maps with ends that wind.
Assume that at $z\in\Gamma$ that $\tl u$ is asymptotic to $\gamma_z^{m_z}$ and
that at $w\in\Gamma'$ $\tl v$ is asymptotic to $\gamma_w^{m_w}$
with $m_{z}$, $m_{w}$ nonzero integers and $\gamma_{z}$, $\gamma_{w}\in\Pe^{0}\MH$,
simple periodic orbits.
Then
\begin{equation}\label{e:gin-cylindrical}
\begin{aligned}
\gin{\tl u}{\tl v}
&
=\sum_{w\in\Gamma'_+}
|m_w|\bp{
\inum(\gamma_w, u)+\sum_{\substack{z\in\Gamma_+ \\ \gamma_z=\gamma_w}}
|m_z|\;\Mxx{ \tfrac{\alpha(\gamma_w^{m_w})}{|m_w|}, \tfrac{\alpha(\gamma_z^{m_z})}{|m_z|}}
+d_0(\tl u; z)
} \\
&
\hskip.125in+\sum_{z\in\Gamma_-}
|m_z|\bp{
\inum(\gamma_z, v)+\sum_{\substack{w\in\Gamma'_- \\ \gamma_w=\gamma_z}}
|m_w|\; \Mxx{\tfrac{\alpha(\gamma_z^{m_z})}{|m_z|}, \tfrac{\alpha(\gamma_w^{m_w})}{|m_w|}}
+d_0(\tl v; w)
} \\
&
\hskip.125in+\sum_{\substack{(z, w)\in\Gamma_{-}\times\Gamma'_{+} \\ \gamma_z=\gamma_w}}
|m_w|d_0(\tl u; z)+|m_z|d_0(\tl v; w)-|m_{w}m_{z}|\bbr{\tfrac{\alpha(\gamma_z^{m_z})}{|m_z|}+\tfrac{\alpha(\gamma_w^{m_w})}{|m_w|}}
\end{aligned}
\end{equation}
where
\begin{equation}\label{e:mxx-define}
\Mxx{x, y}:=\max\br{x, y}-y.
\end{equation}
\end{theorem}

\begin{proof}
We apply the $\R$-action to $\tl u$, and note that
\[
\lim_{c\to\infty}(c\cdot\tl u)+\tl v
=(\tl u\odot[+_{z\in\Gamma_-}\tl\gamma_z^{-m_z}])
+
([+_{w\in\Gamma'_{+}}\tl\gamma_w^{m_w}]\odot\tl v)
\]
where the limit is to be understood in the sense of \cite{behwz},
and ``$+$'' denotes the disjoint union of the maps.
Consequently, the homotopy invariance of the relative intersection number
along with \fullref{l:rin-sft-compactness} imply that
\begin{equation}\label{e:gin-cylindrical-pf}
\rin^\Phi(\tl u, \tl v)=\lim_{c\to\infty}\rin^\Phi(c\cdot\tl u, \tl v)
=\sum_{w\in\Gamma'_{+}}\rin^\Phi(\tl u, \tl\gamma_{w}^{m_{w}})
+\sum_{z\in\Gamma_{-}}\rin^\Phi(\tl v, \tl\gamma_{z}^{-m_{z}}).
\end{equation}
Applying \fullref{l:rin-cylinder} and the definition of $d_0(\tl u; z)$, we have that
\begin{align*}
i^\Phi(\tl u, \tl\gamma_w^{m_w})&=m_{w}
\bp{
\inum(\gamma_w, u)
-\sum_{\substack{z\in\Gamma \\ \gamma_z=\gamma_w}}\winfty^\Phi(u; z)
} \\
&=m_{w}
\bp{
\inum(\gamma_w, u)
+\sum_{\substack{z\in\Gamma \\ \gamma_z=\gamma_w}}[d_0(\tl u; z)-\alpha^\Phi(\gamma_z^{m_z})]
} \\
\end{align*}
and similarly
\[
i^\Phi(\tl v, \tl\gamma_z^{-m_z})
=|m_{z}|
\bp{
\inum(\gamma_z, v)
+\sum_{\substack{w\in\Gamma' \\ \gamma_w=\gamma_z}}[d_0(\tl v; w)-\alpha^\Phi(\gamma_w^{m_w})]
}
\]
Substituting these formulas in \eqref{e:gin-cylindrical-pf} and adding
\[
\Omega^\Phi(\tl u, \tl v)=\sum_{\substack{(z, w)\in\Gamma\times\Gamma' \\ \gamma_z=\gamma_w \\ m_z m_w>0}}
m_z m_w\max\br{\tfrac{\alpha(\gamma_z^{m_z})}{|m_z|}, \tfrac{\alpha(\gamma_w^{m_w})}{|m_w|}}
\]
to both sides gives
\begin{align*}
\gin{\tl u}{\tl v}&=\rin^\Phi(\tl u, \tl v)+\Omega^\Phi(\tl u, \tl v) \\
&=\sum_{w\in\Gamma'_+}
|m_w|\bp{
\inum(\gamma_w, u)+\sum_{\substack{z\in\Gamma_+ \\ \gamma_z=\gamma_w}}
|m_z|\;\Mxx{ \tfrac{\alpha(\gamma_w^{m_w})}{|m_w|}, \tfrac{\alpha(\gamma_z^{m_z})}{|m_z|}}
+d_0(\tl u; z)
} \\
&
\hskip.125in+\sum_{z\in\Gamma_-}
|m_z|\bp{
\inum(\gamma_z, v)+\sum_{\substack{w\in\Gamma'_- \\ \gamma_w=\gamma_z}}
|m_w|\; \Mxx{\tfrac{\alpha(\gamma_z^{m_z})}{|m_z|}, \tfrac{\alpha(\gamma_w^{m_w})}{|m_w|}}
+d_0(\tl v; w)
} \\
&
\hskip.125in+\sum_{\substack{(z, w)\in\Gamma_{-}\times\Gamma'_{+} \\ \gamma_z=\gamma_w}}
|m_w|d_0(\tl u; z)+|m_z|d_0(\tl v; w)-|m_{w}m_{z}|\bbr{\tfrac{\alpha(\gamma_z^{m_z})}{|m_z|}+\tfrac{\alpha(\gamma_w^{m_w})}{|m_w|}}
\end{align*}
as claimed,
where we've abbreviated $\Mxx{x,y}=\max\br{x, y}-y$.
\end{proof}

We note that in the case that the maps $\tl u$ and $\tl v$ from this theorem are pseudoholomorphic
that \eqref{e:gin-cylindrical} expresses the generalized intersection number $\gin{\tl u}{\tl v}$
as a sum of nonnegative terms.  Indeed the nonnegativity of each of the terms is discussed
following \fullref{c:int-cylinder}, with the exception of the terms of the form
\[
-|m_{w}m_{z}|\bbr{\tfrac{\alpha(\gamma_z^{m_z})}{|m_z|}+\tfrac{\alpha(\gamma_w^{m_w})}{|m_w|}}.
\]
The fact that these terms contribute nonnegatively to $\gin{\tl u}{\tl v}$
is addressed in the proof of \fullref{c:gin-zero} below.

The fact that \eqref{e:gin-cylindrical} expresses
$\gin{\tl u}{\tl v}$  for pseudoholomorphic curves as a sum of nonnegative terms allows us to prove 
the following result
which gives necessary and sufficient conditions for
two curves to have generalized intersection number equal to zero.

\begin{corollary}\label{c:gin-zero}
Let $[\Sigma, j, \Gamma, \tl u=(a, u)]$ and $[\Sigma', j', \Gamma', \tl v=(b, v)]\in\M(M, \Ha, J)$
be pseudoholomorphic curves,
and assume that no component of $\tl u$ or $\tl v$ lies in an orbit cylinder.
Then the following are equivalent:
\begin{enumerate}
\item The generalized intersection number $\gin{\tl u}{\tl v}=0$.

\item All of the following hold:
     \begin{enumerate}
     \item The map $u$ does not intersect any of the positive asymptotic limits of $v$.
     \item The map $v$ does not intersect any of the negative asymptotic limits of $u$.
     \item If $\gamma$ is a periodic orbit so that at $z\in\Gamma$, $\tl u$ is asymptotic
     to $\gamma^{m_z}$ and at $w\in\Gamma'$, $\tl v$ is asymptotic to $\gamma^{m_w}$, then:
          \begin{enumerate}
          \item If $m_z$ and $m_w$ are both positive then $d_0(\tl u; z)=0$ and
          \[
          \tfrac{\alpha(\gamma^{m_z})}{m_z}\ge\tfrac{\alpha(\gamma^{m_w})}{m_w}.
          \]
          \item If $m_z$ and $m_w$ are both negative then $d_0(\tl v; w)=0$ and
          \[
          \tfrac{\alpha(\gamma^{m_w})}{|m_w|}\ge\tfrac{\alpha(\gamma^{m_z})}{|m_z|}.
          \]
          \item\label{item:gin-zero-pos-neg} If $m_z<0$ and $m_w>0$ then $d_0(\tl u; z)=d_0(\tl v; w)=0$ and
          $\gamma^{m_z}$ and $\gamma^{m_w}$ are both even orbits; or equivalently
          \[
          \tfrac{\winfty(\tl u; z)}{m_{z}}=\tfrac{\winfty(\tl v; w)}{m_{w}}.
          \]
     \end{enumerate}
     \end{enumerate}

\item All of the following hold:
     \begin{enumerate}
     \item The map $u$ does not intersect any of the asymptotic limits of $v$.
     \item The map $v$ does not intersect any of the asymptotic limits of $u$.
     \item If $\gamma$ is a periodic orbit so that at $z\in\Gamma$, $\tl u$ is asymptotic
     to $\gamma^{m_z}$ and at $w\in\Gamma'$, $\tl v$ is asymptotic to $\gamma^{m_w}$, then
     $d_0(\tl u; z)=0$, $d_0(\tl v; w)=0$.  Further
     \begin{enumerate}
          \item if $\gamma$ is elliptic, then $m_z$ and $m_w$ have the same sign, and
          \[
          \tfrac{\alpha(\gamma^{m_z})}{m_z}=\tfrac{\alpha(\gamma^{m_w})}{m_w}.
          \]
          \item if $\gamma$ is odd, hyperbolic then either $m_z$ and $m_w$ are both even, or
          $m_z=m_w$.
     \end{enumerate}
\end{enumerate}

\end{enumerate}
\end{corollary}

\begin{proof}
As discussed preceding the statement of this corollary,
this will follow from \fullref{t:gin-cylindrical} and in particular the fact that
\fullref{t:gin-cylindrical} gives an expression of $\gin{\tl u}{\tl v}$ as a sum of
nonnegative quantities.
Indeed, as discussed following \fullref{c:int-cylinder},
it follows from \eqref{e:int-cylinder} that the
terms involving intersections of the projected maps $u$ and $v$ with periodic orbits are nonnegative.
Moreover, \fullref{l:wind-infinity} establishes that the $d_{0}$ terms are nonnegative,
while
it is clear that each other grouped term on the right hand side
of equation \eqref{e:gin-cylindrical} is nonnegative except possibly for the term in the final
sum of the form
\[
\Delta_1(\gamma, m, n):=-|mn|\bbr{\tfrac{\alpha(\gamma^m)}{|m|}+\tfrac{\alpha(\gamma^n)}{|n|}}
\]
where here $m$ and $n$ have opposite sign.  The nonnegativity of this term will follow from
\fullref{l:alpha-iterate}.  Without loss of generality, assume that
$m>0$ and $n<0$.  So that the expression takes the form
\[
\Delta_1(\gamma, m, n)=mn\bbr{\tfrac{\alpha(\gamma^m)}{m}+\tfrac{\alpha(\gamma^n)}{-n}}
\]
If $\gamma$ is an even orbit then \eqref{e:alpha-even} gives us
\[
\Delta_1(\gamma, m, n)=mn[\alpha(\gamma)-\alpha(\gamma)]=0.
\]
If $\gamma$ is odd, hyperbolic, \eqref{e:alpha-oddhyperbolic} gives us
\[
\Delta_1(\gamma, m, n)=\tfrac{1}{2}[mp(n)-np(m)]
\]
which is nonnegative since $m>0$ and $n<0$, and equals zero if and only if both $m$ and $n$ are even.
Finally if $\gamma$ is elliptic, by \eqref{e:alpha-elliptic} there is an irrational $\theta$ so that
\begin{align*}
\Delta_1(\gamma, m, n)
&=n\fl{m\theta}-m\fl{n\theta} \\
&\ge\ceil{nm\theta}-\fl{nm\theta}=1,
\end{align*}
where the inequality follows from the assumptions that $n<0$, $m>0$ and that $\theta$ is irrational.
To summarize, we have seen that for $m$ and $n$ integers with opposite sign, we have that
\[
-|mn|\bbr{\tfrac{\alpha(\gamma^m)}{|m|}+\tfrac{\alpha(\gamma^n)}{|n|}}\ge 0
\]
and equality occurs if and only if $\gamma^m$ and $\gamma^n$ are both even orbits.

Given this, it is clear from \eqref{e:gin-cylindrical} that
$\gin{\tl u}{\tl v}=0$ if and only if:
\begin{itemize}
\item
$
\inum(\gamma_w, u)=0
$
for all $w\in\Gamma'_+$, 
\item
$
\inum(\gamma_z, v)=0
$
for all $z\in\Gamma_-$,
\item $d_0(\tl u; z)=0$ and
\[
\max\br{\tfrac{\alpha(\gamma_z^{m_z})}{|m_z|}, \tfrac{\alpha(\gamma_w^{m_w})}{|m_w|}}
-\tfrac{\alpha(\gamma_z^{m_z})}{m_z}=0
\]
for all $(z, w)\in\Gamma_+\times\Gamma'_+$ with $\gamma_z=\gamma_w$,
\item $d_0(\tl v; w)=0$ and
\[
\max\br{\tfrac{\alpha(\gamma_z^{m_z})}{|m_z|}, \tfrac{\alpha(\gamma_w^{m_w})}{|m_w|}}
-\tfrac{\alpha(\gamma_w^{m_w})}{m_w}=0
\]
for all $(z, w)\in\Gamma_-\times\Gamma'_-$ with $\gamma_z=\gamma_w$, and
\item
$d_0(\tl u; z)=d_0(\tl v; w)=0$ and
\begin{equation}\label{e:pos-neg-punctures}
-|m_z m_w|\bbr{\tfrac{\alpha(\gamma_z^{m_z})}{|m_z|}+\tfrac{\alpha(\gamma_w^{m_w})}{|m_w|}}=0
\end{equation}
for all $(z, w)\in\Gamma_-\times\Gamma'_+$ with $\gamma_z=\gamma_w$.
\end{itemize}
The discussion of the previous paragraph tells us that \eqref{e:pos-neg-punctures} is true if and only if
$\gamma_z^{m_z}$ and $\gamma_w^{m_w}$ are both even orbits.
Moreover, we have from \fullref{l:wind-infinity} together with $m_{z}<0$, $m_{w}>0$, and the
previously observed nonnegativity of the left-hand-side of \eqref{e:pos-neg-punctures} that
\[
\tfrac{\winfty(\tl v; w)}{m_{w}}\le \tfrac{\alpha(\gamma_{w}^{m_{w}})}{m_{w}}\le
\tfrac{\alpha(\gamma_{z}^{m_{z}})}{m_{z}} \le \tfrac{\winfty(\tl u; z)}{m_{z}}
\]
in general.  This together with the definition \eqref{e:def} of $d_{0}$ tells us
that the last listed condition above, that is $d_0(\tl u; z)=d_0(\tl v; w)=0$ and \eqref{e:pos-neg-punctures},
is true precisely when
\[
\tfrac{\winfty(\tl v; w)}{m_{w}}=\tfrac{\winfty(\tl u; z)}{m_{z}}
\]
so the two conditions given in item \eqref{item:gin-zero-pos-neg} are equivalent.
The fact that
$(1)\iff(2)$ in the statement of the corollary follows immediately then from the above list of conditions.

To see that $(2)\iff(3)$, we use the symmetry of the generalized intersection number with the asymmetry of
statement $(2)$.  Indeed if $(2)$ is true as written,
then $(1)$ is true, and thus $(2)$ will be true with the roles of $\tl u$ and $\tl v$ reversed.
We immediately find that $(2)$ holds precisely when
\begin{itemize}
\item $\inum(\gamma_w, u)=0$ for all $w\in\Gamma'$,
\item $\inum(\gamma_z, v)=0$ for all $z\in\Gamma$,
\item $d_0(\tl u; z)=d_0(\tl v; w)=0$ and
\[
\tfrac{\alpha(\gamma_z^{m_z})}{|m_z|}=\tfrac{\alpha(\gamma_w^{m_w})}{|m_w|}
\]
for all $(z, w)\in\Gamma\times\Gamma'$ with $\gamma_z=\gamma_w$ and $m_z m_w>0$, and
\item $d_0(\tl u; z)=d_0(\tl v; w)=0$ and
\[
-|m_z m_w|\bbr{\tfrac{\alpha(\gamma_z^{m_z})}{|m_z|}+\tfrac{\alpha(\gamma_w^{m_w})}{|m_w|}}=0
\]
for all $(z, w)\in\Gamma\times \Gamma'$ with $\gamma_z=\gamma_w$ and $m_z m_w<0$.
\end{itemize}
Using \eqref{e:alpha-even} we see that
\[
\tfrac{\alpha(\gamma^{m})}{|m|}=\tfrac{\alpha(\gamma^{n})}{|n|}
\]
for any $m$ and $n$ with $mn>0$ if $\gamma$ is even,
and using \eqref{e:alpha-oddhyperbolic} we see that
\[
\tfrac{\alpha(\gamma^{m})}{|m|}=\tfrac{\alpha(\gamma^{n})}{|n|} \iff \tfrac{p(m)}{m}=\tfrac{p(n)}{n}
\]
if $\gamma$ is odd hyperbolic.  This last statement is true
exactly when $m$ and $n$ are both even, or $m$ and $n$ are equal.
These observations with the discussion of the first paragraph shows that the conditions listed above
are equivalent to those listed in $(3)$, so we see that
$(2)\iff(3)$.
\end{proof}

We note that if $\gin{\tl u}{\tl v}=0$ then $\gin{\tl u}{c\cdot\tl v}=0$ for all $c\in\R$, and in particular,
\fullref{t:gin-positivity} implies that the actual algebraic intersection number
$\inum(\tl u, c\cdot\tl v)$ is zero for all $c\in\R$,
provided $\tl u$ and $\tl v$ do not have any components
with images that differ by an $\R$-shift.
By positivity of intersections, this implies that
$\tl u$ is disjoint from every $\R$-translate of $\tl v$ which implies that the projected curves
$u$ and $v$ in the $3$-manifold $M$ either have identical image
or do not intersect.  Thus the preceding corollary gives
sufficient conditions for the projected curves $u$ and $v$ to not intersect.
The converse is not true, as is it possible that projected
pseudoholomorphic maps $u$ and $v$ do not intersect, but
still have $\gin{\tl u}{\tl v}\ne0$ since the asymptotic intersection number of the two curves could be nonzero.
However, since $\R$-shifting a curve changes the asymptotic description from \fullref{t:asymp}
in a predictable way, it is possible to use \fullref{t:gin-cylindrical} to identify necessary and sufficient
conditions for two projected pseudoholomorphic curves to not intersect.

The key element in understanding this problem is the following lemma concerning the behavior of the
total asymptotic intersection number of two curves under $\R$-shifting.

\begin{lemma}\label{l:tain-cylindrical}
Let $[\Sigma, j, \Gamma, \tl u]$ and $[\Sigma', j', \Gamma', \tl v]\in\MM$
be asymptotically cylindrical pseudoholomorphic curves, and let $\tl v_{c}=(b+c, v)$ be
the curve resulting from translating $\tl v$ in the $\R$-coordinate by $c$.
Assume no component of $\tl u$ or $\tl v$ has image contained in an orbit cylinder.
Then the
algebraic intersection number
$\inum(\tl u, \tl v_{c})$
and the total asymptotic intersection index
$\delta_{\infty}(\tl u, \tl v_{c})$
are defined for all but a finite number of values of $c\in\R$.
Moreover
\[
\delta_{\infty}(\tl u, \tl v_{c})
\ge \sum_{\substack{(z, w)\in \Gamma\times\Gamma' \\ \gamma_z=\gamma_w \\ m_{z}m_{w}>0}}
m_{z}m_{w}
\bp{
     \max\br{\tfrac{\alpha(\tl u; z)}{|m_z|}, \tfrac{\alpha(\tl v; w)}{|m_w|}}
     -\max\br{\tfrac{\winfty(\tl u; z)}{|m_z|}, \tfrac{\winfty(\tl v; w)}{|m_w|}}
},
\]
with equality occurring for all but a finite number of the values of $c\in\R$ for which
$\delta_{\infty}(\tl u, \tl v_{c})$ is defined.
Furthermore,
\begin{align*}
\inum\bp{\tl u, \tl v_{c}}
&\le \gin{\tl u}{\tl v} \\
&\hskip.125in
-\sum_{\substack{(z, w)\in \Gamma\times\Gamma' \\ \gamma_z=\gamma_w \\ m_{z}m_{w}>0}}
m_{z}m_{w}
\bp{
     \max\br{\tfrac{\alpha(\tl u; z)}{|m_z|}, \tfrac{\alpha(\tl v; w)}{|m_w|}}
     -\max\br{\tfrac{\winfty(\tl u; z)}{|m_z|}, \tfrac{\winfty(\tl v; w)}{|m_w|}}
}
\end{align*}
and strict inequality occurs for at most a finite set of values of $c\in\R$ for which
$\inum\bp{\tl u, \tl v_{c}}$ is defined.
\end{lemma}

\begin{proof}
There will be values of $c$ for which
 the
algebraic intersection number
$\inum(\tl u, \tl v_{c})$
and the total asymptotic intersection index
$\delta_{\infty}(\tl u, \tl v_{c})$
will not be defined if there are components of $\tl u$ and $\tl v$ which project to the same image
in the three-manifold, but it is clear the lifted maps
$\tl u$ and $\tl v_{c}$ can have components with identical image in $\R\times M$ only for
finite number of values of $c\in\R$, since we assume that no components of the curves
are fixed by the $\R$-action.

We consider a pair of end models $(\Sigma, j, z, \tl u)$ and $(\Sigma', j', w, \tl v=(b, v))$
and assume that there is a $\gamma$ so that $\tl u$ is asymptotic to $\gamma^m$ and $\tl v$ is asymptotic
to $\gamma^n$ with $m$ and $n$ having the same sign.
Assume we have the asymptotic representations
\begin{gather*}
\tl u(\phi(s,t))=(m\tau s, \exp_{\gamma^m(t)}e^{\lambda_1 s}[e_1(t)+r_1(s,t)]) \\
\tl v(\psi(s',t'))=(n\tau s', \exp_{\gamma^n(t')}e^{\lambda_2 s'}[e_2(t')+r_2(s',t')]).
\end{gather*}
with the $\lambda_i$/$e_i$ eigenvalues/eigenvectors of the appropriate asymptotic operators, and
the $r_i$ converging exponentially to $0$ as $s\to\infty$.
Then, according to \fullref{l:rain-compute}, and \fullref{l:wind-rel}, we have that
\begin{equation}\label{e:rain-windbound}
i^\Phi([\tl u; z], [\tl v; w])\ge -mn\max\br{
\tfrac{\winfty^\Phi(\tl u; z)}{|m|}, \tfrac{\winfty^\Phi(\tl v; w)}{|n|}
}
\end{equation}
with strict inequality occurring only if there is an integer $j$ so that
\[
e_1(nt)=e_2(mt+\tfrac{j}{mn})
\]
for all $t\in S^1$.
Considering the shifted end
\[
\tl v_c=(b+c, v)
\]
we have, for $\psi_c(s,t):=\psi(s+\frac{c}{n\tau}, t)$ that
\begin{align*}
\tl v_c(\psi_c(s,t))&=(n\tau s, \exp_{\gamma^n(t)}e^{\lambda_2 (s+\tfrac{c}{n\tau})}[e_2(t)+r_2(s+\tfrac{c}{n\tau},t)] \\
&=(n\tau s, \exp_{\gamma^n(t)}e^{\lambda_2 s}[ e^{c\tfrac{\lambda}{n\tau}} e_2(t)+r_3(s,t)] \\
\end{align*}
with $r_3(s,t)=e^{c\tfrac{\lambda}{n\tau}}r_2(s+\tfrac{c}{n\tau})$
converging exponentially to zero as $s\to\infty$.  Thus, $\R$-shifting an end has the effect
of scaling the eigenvalue appearing in the asymptotic formula (once the formula has been rewritten so that
the first component does not contain an $\R$-shift).
Consequently,  we get either the strict inequality
\begin{align*}
i^\Phi([\tl u; z], [\tl v_c; w])&> -mn\max\br{
\tfrac{\winfty^\Phi(\tl u; z)}{|m|}, \tfrac{\winfty^\Phi(\tl v_c; w)}{|n|}
} \\
&= -mn\max\br{
\tfrac{\winfty^\Phi(\tl u; z)}{|m|}, \tfrac{\winfty^\Phi(\tl v; w)}{|n|}
}
\end{align*}
or that $i^\Phi([\tl u; z], [\tl v_c; w])$ is not defined
for at most one value of $c\in\R$.
For all other values we will have
\[
i^\Phi([\tl u; z], [\tl v_c; w])= -mn\max\br{
\tfrac{\winfty^\Phi(\tl u; z)}{|m|}, \tfrac{\winfty^\Phi(\tl v; w)}{|n|}
}.
\]
Adding $mn\max\br{\tfrac{\alpha(\gamma^m)}{|m|}, \tfrac{\alpha(\gamma^n)}{|n|}}$
to both sides of this, we can conclude that the asymptotic intersection number of the two ends
satisfies
\[
\delta_\infty([\tl u; z], [\tl v; w])\ge mn
\bp{
\max\br{\tfrac{\alpha(\gamma^m)}{|m|}, \tfrac{\alpha(\gamma^n)}{|n|}}
-\max\br{
\tfrac{\winfty^\Phi(\tl u; z)}{|m|}, \tfrac{\winfty^\Phi(\tl v; w)}{|n|}}
}
\]
with the strict inequality occurring for at most one value of $c\in\R$.

The claims of the Lemma now follow from applying the results of the previous paragraph 
pairwise to ends of $\tl u$ and $\tl v_c$ which are asymptotic to coverings of the same orbit with the same sign,
together with \eqref{e:gin-compute}.
\end{proof}

As a corollary of \fullref{t:gin-cylindrical} and \fullref{l:tain-cylindrical} we have the following.

\begin{corollary}\label{c:int-bound-cylindrical}
Let $[\Sigma, j, \Gamma, \tl u=(a,u)]$, $[\Sigma', j', \Gamma', \tl v=(b, v)]\in\MM$
be pseudoholomorphic curves, and let $\tl v_{c}=(b+c, v)$ be
the curve resulting from translating $\tl v$ in the $\R$ coordinate by $c$.
Assume that no component of $\tl u$ or $\tl v$ has image contained in an orbit cylinder,
and let $\gamma_{z}^{m_{z}}$ and $\gamma_{w}^{m_{w}}$ denote the
asymptotic limits of $\tl u$ at $z\in\Gamma$ and $\tl v$ at $w\in\Gamma'$
as in \fullref{t:gin-cylindrical}.
Then the algebraic intersection number $\inum(\tl u, \tl v_{c})$ is defined for all but
a finite number of values of $c\in\R$.
Moreover,
with $\Mx$ as defined in \eqref{e:mxx-define},
we have that
\begin{equation}\label{e:int-projected}
\begin{aligned}
\inum\bp{\tl u, \tl v_{c}}&\le
\sum_{w\in\Gamma'_+}
|m_w|\bp{
\inum(\gamma_w, u)+\sum_{\substack{z\in\Gamma_+ \\ \gamma_z=\gamma_w}}
|m_z|\;\Mxx{\tfrac{\winfty(\tl v; w)}{|m_w|}, \tfrac{\winfty(\tl u; z)}{|m_z|}}
} \\
&\hskip.125in+\sum_{z\in\Gamma_-}
|m_z|\bp{
\inum(\gamma_z, v)+\sum_{\substack{w\in\Gamma'_- \\ \gamma_w=\gamma_z}}
|m_w|\;\Mxx{\tfrac{\winfty(\tl u; z)}{|m_z|}, \tfrac{\winfty(\tl v; w)}{|m_w|}}
} \\
&\hskip.125in+\sum_{\substack{(z, w)\in\Gamma_{-}\times\Gamma'_{+} \\ \gamma_z=\gamma_w}}
|m_w|d_0(\tl u; z)+|m_z|d_0(\tl v; w)-|m_{w}m_{z}|\bbr{\tfrac{\alpha(\gamma_z^{m_z})}{|m_z|}+\tfrac{\alpha(\gamma_w^{m_w})}{|m_w|}}
\end{aligned}
\end{equation}
and strict inequality occurs for at most a finite set of values of $c\in\R$ for which
$\inum\bp{\tl u, \tl v_{c}}$ is defined.
\end{corollary}

\begin{proof}
This follows from subtracting
\[
\sum_{\substack{(z, w)\in \Gamma\times\Gamma' \\ \gamma_z=\gamma_w \\ m_{z}m_{w}>0}}
m_{z}m_{w}
\bp{
     \max\br{\tfrac{\alpha(\tl u; z)}{|m_z|}, \tfrac{\alpha(\tl v; w)}{|m_w|}}
     -\max\br{\tfrac{\winfty(\tl u; z)}{|m_z|}, \tfrac{\winfty(\tl v; w)}{|m_w|}}
}
\]
from \eqref{e:gin-cylindrical} and using \fullref{l:tain-cylindrical}.
\end{proof}

We can now prove \fullref{t:no-isect-intro}, which
gives necessary and sufficient conditions for two projected curves to not intersect.
We restate the result here for the convenience of the reader.

\begin{theorem}[\fullref{t:no-isect-intro}]
Let $[\Sigma, j, \Gamma, \tl u=(a, u)]$ and $[\Sigma', j', \Gamma', \tl v=(b, v)]\in\M(M, \Ha, J)$
be pseudoholomorphic curves,
and assume that no component of $\tl u$ or $\tl v$ lies in in orbit cylinder,
and that the projected curves $u$ and $v$ do not have identical image on any component of their domains.
Then the following are equivalent:
\begin{enumerate}
\item The projected curves $u$ and $v$ do not intersect.

\item All of the following hold:
     \begin{enumerate}
     \item The map $u$ does not intersect any of the positive asymptotic limits of $v$.
     \item The map $v$ does not intersect any of the negative asymptotic limits of $u$.
     \item If $\gamma$ is a periodic orbit so that at $z\in\Gamma$, $\tl u$ is asymptotic
     to $\gamma^{m_z}$ and at $w\in\Gamma'$, $\tl v$ is asymptotic to $\gamma^{m_w}$, then:
          \begin{enumerate}
          \item If $m_z$ and $m_w$ have the same sign then
          \[
          \tfrac{\winfty(\tl u; z)}{m_z}\ge\tfrac{\winfty(\tl v; w)}{m_w}.
          \]
          \item If $m_z<0$ and $m_w>0$ then $d_0(\tl u; z)=d_0(\tl v; w)=0$ and
          $\gamma^{m_z}$ and $\gamma^{m_w}$ are both even orbits; or equivalently
          \[
          \tfrac{\winfty(\tl u; z)}{m_{z}}=\tfrac{\winfty(\tl v; w)}{m_{w}}.
          \]
    \end{enumerate}
     \end{enumerate}

\item All of the following hold:
     \begin{enumerate}
     \item The map $u$ does not intersect any of the asymptotic limits of $v$.
     \item The map $v$ does not intersect any of the asymptotic limits of $u$.
     \item If $\gamma$ is a periodic orbit so that at $z\in\Gamma$, $\tl u$ is asymptotic
     to $\gamma^{m_z}$ and at $w\in\Gamma'$, $\tl v$ is asymptotic to $\gamma^{m_w}$, then
          \[
          \tfrac{\winfty(\tl u; z)}{m_z} = \tfrac{\winfty(\tl v; w)}{m_w}.
          \]
     \end{enumerate}

\end{enumerate}
\end{theorem}

\begin{proof}
The fact that $(1)\iff(2)$ follows from an argument similar to that in the proof of \fullref{c:gin-zero}
using now \eqref{e:int-projected} instead of \eqref{e:gin-cylindrical}.  Indeed, we have that
the projected curves $u$ and $v$ are disjoint if and only if $\tl u$ and $\tl v_c$ are disjoint for
all $c\in\R$, which by positivity of intersections is true precisely when $\inum(\tl u, \tl v_c)=0$
for all $c\in\R$.  Again, by positivity of intersections and \fullref{c:int-bound-cylindrical}
it is necessary and sufficient for the right hand side of \eqref{e:int-projected}
to vanish.
Since the right hand side of \eqref{e:int-projected} decomposes into nonnegative terms, we need
to require each of these terms to vanish.
We can argue as in \fullref{c:gin-zero} that
the vanishing of each of these terms corresponds to the conditions listed
in $(2)$ above.  Therefore $(1)\iff(2)$ as claimed.

The fact that $(2)\iff(3)$, like in \fullref{c:gin-zero}, follows from the symmetry of the intersection number and
the asymmetry of statement $(2)$.  Indeed, if $(2)$ holds for $\tl u$ and $\tl v$ as stated, then $(2)$ holds with
the roles of $\tl u$ and $\tl v$ reversed.
Thus $(2)$ holds precisely when:
\begin{itemize}
\item $\inum(\gamma_w, u)=0$ for all $w\in\Gamma'$,
\item $\inum(\gamma_z, v)=0$ for all $z\in\Gamma$,
\item for every $(z, w)\in\Gamma\times\Gamma'$ with $\gamma_z=\gamma_w$ and $m_z m_w>0$
\[
\tfrac{\winfty(\tl u; z)}{m_z}=\tfrac{\winfty(\tl v; w)}{m_w},
\]
and
\item for every $(z, w)\in\Gamma\times\Gamma'$ with $\gamma_z=\gamma_w$ and $m_z m_w<0$
\[
\tfrac{\winfty(\tl u; z)}{m_z}=\tfrac{\winfty(\tl v; w)}{m_w}.
\]
\end{itemize}
These conditions are equivalent to those in statement $(3)$, and thus $(2)\iff(3)$ as claimed.
\end{proof}

\subsubsection{Positive asymptotic intersection indices and direction of approach to
even orbits}\label{ss:direction-of-approach}
In this section we give a sufficient condition for positivity of the generalized intersection
number of two curves in terms of
direction of approach to shared hyperbolic orbits.

Before stating and proving the result, we will need to establish some terminology.
Let $\gamma$ be a periodic orbit of $\X$, and assume that  either
\begin{itemize}
\item $\gamma$ is a simple even orbit, or
\item $\gamma=\hat\gamma^{2}$ where $\hat\gamma$ is a simple, odd, hyperbolic orbit.
\end{itemize}
In either of these cases, let $\lambda^{-}=\lne(\gamma)\in\sigma(\A_{\gamma, J})$
be the largest negative eigenvalue
of $\A_{\gamma, J}$.

\begin{lemma}\label{l:hyperbolic-direction}
With $\gamma$, $\lambda^{-}$ satisfying the assumptions above, $m\lambda^{-}$
is the largest negative eigenvalue of $\A_{\gamma^{m}, J}$ for any
positive integer $m$.  Moreover, the eigenspace
\[
\ker(\A_{\gamma^{m}, J}-m\lambda^{-})
\] 
is one dimensional for all positive integers $m$, and if $e^{-}\in\ker(\A_{\gamma, J}-\lambda^{-})$
is a basis, then $e^{-}_{m}:=e^{-}(m\cdot)$ is a basis for
$\ker(\A_{\gamma^{m}, J}-m\lambda^{-})$.
\end{lemma}

\begin{proof}
By the assumption that $\gamma$ is an even orbit and the definition of parity \eqref{e:parity},
there exists a $\lambda^{+}>0$
so that eigenvectors with eigenvalue $\lambda^{-}$ and eigenvectors with eigenvalue
$\lambda^{+}$ have the same winding number in any choice of trivialization
of $\gamma^{*}\xh$.  Therefore, by \fullref{l:operator-spectrum}, we know that
\[
\ker(\A_{\gamma, J}-\lambda^{-})
\]
and
\[
\ker(\A_{\gamma, J}-\lambda^{+})
\]
are both $1$-dimensional.

Let $e^{\pm}$ be a basis for $\ker(\A_{\gamma, J}-\lambda^{\pm})$.
By the discussion following \fullref{l:operator-spectrum}, we know that
$e^{\pm}_{m}$ defined by $e^{\pm}_{m}(t)=e^{\pm}(mt)$ are eigenvectors with
eigenvalue $m\lambda^{\pm}$.
Moreover, $e_{m}^{+}$ and $e_{m}^{-}$ will have the same winding in any trivialization,
so it follows from Lemma 3.1 that $e_{m}^{\pm}$ is a basis for
$\ker(\A_{\gamma^{m}, J}-m\lambda^{\pm})$, and since winding is monotonic in
eigenvalue, we know that $m\lambda^{-}$ must be the largest negative eigenvalue
of $\A_{\gamma^{m}, J}$.
\end{proof}

Now, with $\gamma$ as above,
consider a pseudoholomorphic end $[\Sigma, j, z, \tl u]$ which is asymptotic
to $\gamma^{m}$ for some positive integer $m$.
Let $(U, \phi)$ be an asymptotic representation of $u$ near $z$, so that
\begin{equation}\label{e:asymp-rep1}
\tl u(\phi(s,t))=(m\tau s, \exp_{\gamma^m(t)}U(s,t)).
\end{equation}
Then, it follows from a special case of
\fullref{t:asymp} that there exists
a (possibly zero)  $e_{u}\in\ker(\A_{\gamma, J}-\lambda^{-})$ satisfying
\[
e^{-m\lambda^{-} s}[U(s, t)-e_{u}(mt)]\to 0
\]
as $s\to\infty$.
Considering a second end $[\Sigma', j', z', \tl v]$ asymptotic
to $\gamma^{n}$ with $n>0$,
we find an asymptotic representative $(V, \psi)$
and a vector $e_{v}\in\ker(\A_{\gamma, J}-\lambda^{-})$
so that we can write
\begin{equation}\label{e:asymp-rep2}
\tl v(\psi(s,t))=(n\tau s, \exp_{\gamma^{n}(t)}V(s,t))
\end{equation}
with $V$ satisfying
\[
e^{-n\lambda^{-} s}[V(s, t)-e_{v}(nt)]\to 0
\]
as $s\to\infty$.
If there exists a positive real number $c$ so that
\[
e_{u}=ce_{v}
\]
then we say that the ends $[\Sigma, j, z, \tl u]$ and $[\Sigma', j', z', \tl v]$
\emph{approach coverings of $\gamma$ in the same direction}.  Note that in the case
that one of $e_{u}$ or $e_{v}$ is zero, this condition would imply that both of $e_{u}$ and $e_{v}$ are zero.

For simplicity in this discussion, we have only considered ends with positive punctures, but we can analogously
define what it means for ends with negative punctures 
to approach a covering of $\gamma$ in the same direction
by describing the direction of approach
with elements of $\ker(\A_{\gamma_{-}, J}-\hat\lambda^{-})$ where $\gamma_{-}=\gamma^{(-1)}$ is
$\gamma$ traversed backwards, and $\hat\lambda^{-}$ is the largest negative eigenvalue of
$\A_{\gamma_{-}, J}$.\footnote{
	Alternatively, we could use the fact that if $[\Sigma, j, z, (a, u)]$ is an end with a positive puncture
	in $\R\times M$ where $M$ is equipped with the Hamiltonian structure $(\lambda, \omega))$,
	then $[\Sigma, j, z, (-a, u)]$ is an end with a negative puncture in $\R\times M$ where
	$M$ is now equipped with the Hamiltonian structure $(-\lambda, \omega)$.
	This approach is employed in some proofs in \cite{hut02}.
	}

\begin{theorem}
Let $\gamma$ be an even periodic orbit satisfying the assumptions above, let
$[\Sigma, j, z, \tl u]$ and $[\Sigma', j', z', \tl v]$ be distinct (nonintersecting)
pseudoholomorphic ends in $\R\times M$
satisfying the assumptions above, and let
$e_{u}$, $e_{v}\in\ker(\A_{\gamma, J}-\lambda^{-})$ be defined as above.
If $e_{u}=e_{v}$, then the asymptotic intersection index
$\delta_{\infty}([\tl u; z], [\tl v; w])$ of the two ends is positive.
\end{theorem}

\begin{proof}
For simplicity we assume $m$ and $n$ are both positive, but the proof readily adapts to the case
where $m$ and $n$ are both negative.

Letting $e_{u}$ and $e_{v}$ be as above, we first assume that
$e_{u}=e_{v}=0$.
Then, it follows that the eigenvalue of the leading eigenvector
in the formulas for $\tl u$ from \fullref{t:asymp} is strictly less
than $m\lambda^{-}$, and similarly the eigenvalue for the leading term in the formula for
$\tl v$ is strictly less than $n\lambda^{-}$.
Since $\gamma$ is an even orbit,
\fullref{l:operator-spectrum} and
\eqref{e:parity}
imply that
\begin{equation}\label{e:hyp-dir-winding-ineq}
\begin{gathered}
\winfty^{\Phi}(\tl u; z)<\alpha^{\Phi}(\gamma^{m})=m\alpha^{\Phi}(\gamma) \\
\winfty^{\Phi}(\tl v; w)<\alpha^{\Phi}(\gamma^{n})=n\alpha^{\Phi}(\gamma)
\end{gathered}
\end{equation}
where we've used \fullref{l:hyperbolic-direction} to conclude that
$\alpha^{\Phi}(\gamma^{k})=k\alpha^{\Phi}(\gamma)$ for any positive integer $k$.
In the case that $\gamma$ is a simply covered even orbit, we consequently
get from \eqref{e:ain-bound-wind-infty} and \eqref{e:hyp-dir-winding-ineq} that
\begin{align*}
\ain^{\Phi}([\tl u; z], [\tl v; w])
&\ge -mn\max\br{\tfrac{\winfty^{\Phi}(\tl u; z)}{m},\tfrac{\winfty^{\Phi}(\tl v; w)}{n}} \\
&>-mn\alpha^{\Phi}(\gamma) \\
&=-mn\max\br{\tfrac{\alpha^{\Phi}(\gamma^{m})}{m},\tfrac{\alpha^{\Phi}(\gamma^{n})}{n}}
\end{align*}
and so
\[
\delta_{\infty}([\tl u; z], [\tl v; w])
=\ain^{\Phi}([\tl u; z], [\tl v; w])+mn\max\br{\tfrac{\alpha^{\Phi}(\gamma^{m})}{m},\tfrac{\alpha^{\Phi}(\gamma^{n})}{n}}
>0
\]
as claimed.
If $\gamma=\hat\gamma^{2}$ with $\hat\gamma$ a simple, odd, hyperbolic orbit, we have
from \eqref{e:alpha-oddhyperbolic} that
\[
\alpha(\gamma)=\alpha(\hat\gamma^{2})=2\alpha(\hat\gamma)+1.
\]
Consequently we find, again using \eqref{e:ain-bound-wind-infty} and \eqref{e:hyp-dir-winding-ineq},
that
\begin{align*}
\ain^{\Phi}([\tl u; z], [\tl v; w])
&\ge -(2m)(2n)\max\br{\tfrac{\winfty^{\Phi}(\tl u; z)}{2m},\tfrac{\winfty^{\Phi}(\tl v; w)}{2n}} \\
&>-2mn(\alpha^{\Phi}(\gamma)) \\
&=-2mn(2\alpha^{\Phi}(\hat\gamma)+1) \\
&=-(2m)(2n)\max\br{\tfrac{\alpha^{\Phi}(\hat\gamma^{2m})}{2m},\tfrac{\hat\alpha^{\Phi}(\hat\gamma^{2n})}{2n}},
\end{align*}
where in the last line we've used
\[
\alpha^{\Phi}(\hat\gamma^{2k})=k(2\alpha(\hat\gamma)+1)
\]
from \eqref{e:alpha-oddhyperbolic}.
Again we conclude that 
\[
\delta_{\infty}([\tl u; z], [\tl v; w])
=\ain^{\Phi}([\tl u; z], [\tl v; w])+(2m)(2n)
\max\br{\tfrac{\alpha^{\Phi}(\hat\gamma^{2m})}{2m},\tfrac{\alpha^{\Phi}(\hat\gamma^{2n})}{2n}}
>0
\]
as claimed.

Next we assume that $e_{u}=e_{v}\ne 0$.  With
$(U, \phi)$ and $(V, \psi)$ satisfying \eqref{e:asymp-rep1} and \eqref{e:asymp-rep2}
respectively, we can apply \fullref{t:asymp-diff} to write
\[
U(ns, nt)-V(ms, mt)=e^{\mu s}[e(t)+r(s,t)]
\]
with $\mu<0$ an eigenvalue of $\A_{\gamma^{mn}, J}$, $e$ and eigenvector with
eigenvalue $\mu$ and $r$ converging exponentially to zero as $s\to\infty$.
The assumption that $e_{u}=e_{v}$ then implies that $\mu<mn\lambda^{-}$,
which, since $\gamma^{mn}$ is an even orbit, implies that
\[
\wind\Phi^{-1}e<\wind\Phi^{-1}e_{u}(mn\cdot)=mn\wind\Phi^{-1}e_{u}=mn\alpha^{\Phi}(\gamma).
\]
We can therefore conclude that
\[
\wrel^{\Phi}(n\cdot[\tl u; z], m\cdot[\tl v; w])<m^{2}n^{2}\alpha^{\Phi}(\gamma)
\]
which by \fullref{l:rain-compute} implies that
\[
\ain^{\Phi}([\tl u; z], [\tl v; w])>-mn\alpha^{\Phi}(\gamma).
\]
Since $\gamma$ is assumed to be either a simple even orbit, or a double covered odd hyperbolic orbit,
it follows as in the previous paragraph that
\[
\delta_{\infty}([\tl u; z], [\tl v; w])>0
\]
in either case.
\end{proof}

We now prove the main result of this section.

\begin{theorem}[\fullref{t:dir-of-approach-intro}]\label{t:dir-of-approach}
Let $\gamma$ be a periodic orbit satisfying the above assumptions, and
let $[\Sigma, j, \Gamma, \tl u=(a, u)]$ and $[\Sigma', j', \Gamma', \tl v=(b, v)]\in\MM$
be connected pseudoholomorphic curves.
If there are punctures $z\in\Gamma$ and $w\in\Gamma'$ so that
the ends $[\tl u; z]$ and $[\tl v; w]$
approach a cover of $\gamma$ in the same direction,
and so that $[\tl u; z]$ is distinct from the the $\R$-shifted end $[(b+c, v); w]$ for all
$c\in\R$
then
\[
\gin{\tl u}{\tl v}>0.
\]
\end{theorem}

\begin{proof}
For simplicity we assume both punctures $z$ and $w$ are positive, but the argument readily
adapts to the case where they are both negative.

Let $(U, \phi)$ and $(V, \phi)$ satisfy \eqref{e:asymp-rep1} and \eqref{e:asymp-rep2}
and let $e_{u}$ and $e_{v}$ be as defined above.
We first consider the case that $e_{u}=e_{v}=0$.  In this case,
the lemma above shows that the asymptotic intersection index
$\delta_{\infty}([\tl u; z], [\tl v; w])$ is positive, and thus the
total asymptotic intersection index
$\delta_{\infty}(\tl u, \tl v)$ is positive.
Thus, \fullref{t:gin-positivity} gives us
\[
\gin{\tl u}{\tl v}=\inum(\tl u, \tl v)+\delta_{\infty}(\tl u, \tl v)\ge \delta_{\infty}(\tl u, \tl v)>0
\]
as claimed.

In the case that $e_{u}=c e_{v}\ne 0$ for some positive real number $c$, we
apply the observation from the proof of \fullref{l:tain-cylindrical}
that $\R$-shifting a curve has the effect of scaling the eigenvector describing the approach.
Indeed, if we let $\tl v_{c'}$ denote the map
\[
\tl v_{c'}=(b+c', v)
\]
and let $(V_{c'}, \psi_{c'})$ denote an asymptotic representative of $\tl v_{c'}$,
we can choose a $c'$ so that $e_{v_{c'}}$, chosen to satisfy
\[
e^{-n\lambda^{-}s}(V_{c'}(s,t)-e_{v_{c'}}(nt))\to 0 \text{ as } s\to\infty,
\]
is equal to $e_{u}$.
The previous lemma then applies 
to show that $\delta_{\infty}([\tl u; z], [\tl v_{c'}, w])>0$,
and hence, by the homotopy invariance of the generalized intersection number, that
\[
\gin{\tl u}{\tl v}=\gin{\tl u}{\tl v_{c'}}\ge \delta_{\infty}([\tl u; z], [\tl v_{c'}, w])>0
\]
as claimed.
\end{proof}

\subsection{Generalized self-intersection numbers and embeddedness of the projection}
In this section we investigate conditions under which a curve in a cylindrical cobordism has
generalized self-intersection number equal to zero, and we investigate
controls on the embeddedness of the projection of a pseudoholomorphic curve
in $\R\times M$ into the three manifold $M$.  The key observation for this latter question,
originally made by Hofer, Wysocki, and Zehnder in \cite{hwz:prop2},
is that the projected curve is injective if and only if the nonprojected curve
doesn't intersect any of its $\R$-translates.

We first state \fullref{t:gin-cylindrical} in the special case that $\tl u=\tl v$.

\begin{corollary}\label{c:sgin-cylindrical}
Let
$(\Sigma, j, \Gamma, \tl u=(a, u))\in\CM$
be an asymptotically cylindrical map with ends that wind,
and assume that at $z\in\Gamma$, $\tl u$ is asymptotic to $\gamma_{z}^{m_{z}}$,
with $\gamma_{z}\in\Pe^{0}\MH$ a simple periodic orbit and $m_{z}$ a nonzero integer.
Then
\begin{align*}
\gin{\tl u}{\tl u}&=
\sum_{z\in\Gamma}
|m_z|\bp{
\inum(\gamma_z, u)+\sum_{\substack{w\in\Gamma \\ \gamma_w=\gamma_z \\ m_z m_w>0}}
|m_w|\;\Mxx{\tfrac{\alpha(\gamma_z^{m_z})}{|m_z|}, \tfrac{\alpha(\gamma_w^{m_w})}{|m_w|}}
+d_0(\tl u; w)
} \\
&+\sum_{\substack{(z, w)\in\Gamma_{+}\times\Gamma_{-} \\ \gamma_z=\gamma_w}}
|m_w|d_0(\tl u; z)+|m_z|d_0(\tl u; w)-|m_{z}m_{w}|\bbr{\tfrac{\alpha(\gamma_z^{m_z})}{|m_z|}+\tfrac{\alpha(\gamma_w^{m_w})}{|m_w|}}
\end{align*}
with $\Mx$ as defined in \eqref{e:mxx-define}.
\end{corollary}

As a corollary to this and \fullref{c:adj-wind-pi} we
prove the following result
which gives
a set of equivalent conditions to $\gin{\tl u}{\tl u}=0$ for a simple curve.

\begin{corollary}\label{c:self-gin-zero}
Let
$[\Sigma, j, \Gamma, \tl u=(a, u)]\in\MM$
be a simple pseudoholomorphic curve,
and assume that no component of $\tl u$ lies in an orbit cylinder.
Then the following are equivalent:
\begin{enumerate}
\item The generalized self-intersection number of $\tl u$ vanishes, i.e.\ $\gin{\tl u}{\tl u}=0$

\item All of the following hold:
     \begin{enumerate}
     \item The projected map $u$ does not intersect any of its asymptotic limits.
     \item For all $z\in\Gamma$, $d_0(\tl u; z)=0$.
     \item If $\gamma$ is a periodic orbit so that at $z\in\Gamma$, $\tl u$ is asymptotic to $\gamma^{m_z}$
     and at $w\in\Gamma$, $\tl u$ is asymptotic to $\gamma^{m_w}$, then
     \[
     \tfrac{\alpha(\gamma^{m_z})}{m_z}=\tfrac{\alpha(\gamma^{m_w})}{m_w}.
     \] 
     \end{enumerate}

\item All of the following hold:
     \begin{enumerate}
     \item $\ind(\tl u)-\chi(\Sigma)+\#\Gamma_{even}=0$
     \item $\cover(\tl u; z)=1$ for all $z\in\Gamma$.
     \item $\sing(\tl u)=0$ or equivalently, $\tl u$ is embedded and has total asymptotic self-intersection
     index, $\delta_\infty(\tl u)$, equal to zero.
     \end{enumerate}

\item All of the following hold:
     \begin{enumerate}
     \item The map $u$ is an immersion that is everywhere transverse to $\X$.
     \item $d_0(\tl u; z)=0$ for all $z\in\Gamma$.
     \item $\cover(\tl u; z)=1$ for all $z\in\Gamma$.
     \item $\sing(\tl u)=0$ or equivalently, $\tl u$ is embedded and has total asymptotic self-intersection
     index, $\delta_\infty(\tl u)$, equal to zero.
     \end{enumerate}
\end{enumerate}
\end{corollary}

\begin{proof}
The fact that $(1)\iff (2)$ follows from \fullref{c:sgin-cylindrical} and an argument
very similar to that in \fullref{c:gin-zero}.  We omit the details.

The fact that $(1)\iff(4)$ follows from rewriting \eqref{e:adj-wind-pi} as
\begin{align*}
\gin{\tl u}{\tl u}&=2\sing(\tl u) +\wind_\pi(\tl u)+d_0(\tl u) +\cover(\tl u)-\#\Gamma \\
&=2\sing(\tl u)+\wind_\pi(\tl u)+\sum_{z\in\Gamma}\bp{d_0(\tl u; z) +[\cover(\tl u; z)-1]}
\end{align*}
and noting that each term on the right hand side of this is nonnegative.

Finally, the fact that $(4)\iff(3)$ follows directly from \eqref{e:wind-pi-ind}.
Note that \eqref{e:wind-pi-ind} implies that the quantity
$\ind(\tl u)-\chi(\Sigma)+\#\Gamma_{even}$ is nonnegative for curves in
a cylindrical cobordism provided it has no components with image contained in an orbit cylinder.
\end{proof}

We observe that for a connected curve $[\Sigma, j, \Gamma, \tl u]$
satisfying the hypotheses of the previous result,
if $\gin{\tl u}{\tl u}=0$ then the projected curve $u$ is an embedding in the three-manifold.
Indeed, the result shows that $u$ must be an immersion which doesn't intersect any of its asymptotic limits.
Moreover, since for the $\R$-translates $\tl u_c=(a+c, u)$, we have
\[
0\le\inum(\tl u, \tl u_c)\le\gin{\tl u}{\tl u_c}=\gin{\tl u}{\tl u}=0,
\]
it follows from positivity of intersections
that $\tl u$ doesn't intersect any of its $\R$-translates, and hence that the projection $u$ is injective.
As observed in \cite{hwz:prop2}, the asymptotic behavior of $u$ then allows us to conclude that
$u$ is an embedding.
As with the discussion of intersections of curves with distinct projections to the three-manifold,
the converse is not true:
it could well be the case that $u$ is an embedding and $\gin{\tl u}{\tl u}\ne0$ since we could have
that the total asymptotic intersection index $\delta_\infty(\tl u, \tl u_c)$ positive for all $c\in\R\setminus\br{0}$
for which it is defined. 
Again, in this case, since $\R$-shifting the curve changes the asymptotic intersection numbers in a
predictable way, we can
modify the above argument to 
find necessary and sufficient conditions for a curve to have an embedded projection.

Towards this end, we state the following special case of \fullref{c:int-bound-cylindrical}
when $\tl v=\tl u$.

\begin{corollary}\label{c:sint-bound-cylindrical}
Let $[\Sigma, j, \Gamma, \tl u=(a,u)]\in\MM$ be a connected pseudoholomorphic curve
and assume that the image of $\tl u$ is not contained in an orbit cylinder.
Then,
with $\Mx$ as defined in \eqref{e:mxx-define},
we have for any $c\in\R\setminus\br{0}$ that
\begin{align*}
\inum(\tl u, \tl u_{c})
&\le
     \sum_{w\in\Gamma}|m_w|
     \bp{
          \inum(\gamma_w, u)
          +\sum_{\substack{z\in\Gamma \\ \gamma_z=\gamma_w \\ m_z m_w>0 }}
                |m_z|\;\Mxx{\tfrac{\winfty(\tl u; w)}{|m_w|}, \tfrac{\winfty(\tl u; z)}{|m_z|}}
     } \\
&\hskip.325in
     +\sum_{\substack{(z,w)\in\Gamma_{-}\times\Gamma_{+} \\ \gamma_z=\gamma_w }}
          -|m_{z}m_{w}|
          \bp{\tfrac{\winfty(\tl u; z)}{|m_z|}
          +\tfrac{\winfty(\tl u; w)}{|m_w|}}
\end{align*}
with strict inequality occurring for at most a finite number of values of $c\in\R\setminus\br{0}$.
\end{corollary}

Along similar lines, keeping track of how the asymptotic intersection numbers change with $\R$-shifts
allows the following adjustment to \fullref{c:adj-wind-pi}

\begin{lemma}\label{l:sint-bound-wind-pi}
Let $[\Sigma, j, \Gamma, \tl u=(a,u)]\in\MM$ be a connected, simple pseudoholomorphic curve
and assume 
the image of $\tl u$ is not contained in 
 an orbit cylinder.
Then, for any $c\in\R\setminus\br{0}$
we have that
\begin{equation}\label{e:sint-bound-wind-pi}
\begin{aligned}
\inum(\tl u, \tl u_{c})&\le\wind_{\pi}(\tl u)+2\delta(\tl u)
     +\sum_{z_{i}\in\Gamma}[\cov(e_{1}(\tl u; z_{i}))-1+\Dtwo(\tl u; z_{i})] \\
&\hskip.25in
+\sum_{\substack{ z_{i},z_{j}\in\Gamma \\ z_{i}\ne z_{j} \\ \gamma_{z_{i}}=\gamma_{z_{j}} \\ m_{z_{i}}m_{z_{j}}>0 }}
i_{\infty}([\tl u; z_{i}], [\tl u; z_{j}])
+m_{z_{i}}m_{z_{j}}
  \max_{\ell=i, j}\br{\tfrac{\winfty(\tl u; z_{\ell})}{|m_{z_\ell}|}}
\end{aligned}
\end{equation}
where $e_{1}(\tl u; z_{i})$ is as defined in \eqref{e:leading-evect},
$\delta(u)$ is as defined in \eqref{e:self-int-index}, and
$\Dtwo(\tl u; z_{i})$ is as defined in \eqref{e:delta2}.
Moreover equality occurs in \eqref{e:sint-bound-wind-pi} for all but a finite number of values of $c\in\R\setminus\br{0}$.
\end{lemma}

\begin{proof}
By a special case of \fullref{c:int-bound-cylindrical} we have that
\begin{align*}
\inum(\tl u, \tl u_{c})&\le \gin{\tl u}{\tl u}\\
&\hskip.25in
     -\sum_{\substack{ z_{i},z_{j}\in\Gamma \\ \gamma_{z_i} =\gamma_{z_j} \\ m_{z_i}m_{z_j}>0 }}
          m_{z_{i}}m_{z_{j}}
          \bp{
               \max_{\ell=i, j}\br{\tfrac{\alpha(\tl u; z_\ell)}{|m_{z_\ell}|}}
               -\max_{\ell=i, j}\br{\tfrac{\wind(\tl u; z_{\ell})}{|m_{z_\ell}|}}
          } \\
&:= I(\tl u)
\end{align*}
with equality occurring for all but a finite number of values of $c\ne 0$.
Applying the adjunction formula \eqref{e:adj-wind-pi} we have that
the right hand side of this inequality is equal to 
\begin{align*}
I(\tl u)&=\wind_{\pi}(\tl u)+d_{0}(\tl u)
     +[\cover(\tl u)-\#\Gamma]+2\sing(\tl u)\\
&\hskip.25in
     -\sum_{\substack{ z_{i},z_{j}\in\Gamma \\ \gamma_{z_i} =\gamma_{z_j} \\ m_{z_i}m_{z_j}>0 }}
          m_{z_{i}}m_{z_{j}}
          \bp{
               \max_{\ell=i, j}\br{\tfrac{\alpha(\tl u; z_{\ell})}{|m_{z_\ell}|}}
               -\max_{\ell=i, j}\br{\tfrac{\wind(\tl u; z_{\ell})}{|m_{z_{\ell}}|}}
          }.
\end{align*}
We rewrite the final sum in this by grouping the the terms with
$z_i=z_j$, and using the definition \eqref{e:def} of $d_0(\tl u; z)$ and \eqref{e:delta1} to get
\begin{align*}
& \sum_{\substack{ z_{i},z_{j}\in\Gamma \\ \gamma_{z_{i}}=\gamma_{z_{j}} \\ m_{z_{i}}m_{z_{j}}>0}}
          m_{z_{i}}m_{z_{j}}
          \bp{
               \max_{\ell=i, j}\br{\tfrac{\alpha(\tl u; z_{\ell})}{|m_{z_\ell}|}}
               -\max_{\ell=i, j}\br{\tfrac{\wind(\tl u; z_{\ell})}{|m_{z_\ell}|}}
          } \\
&\hskip1in=\sum_{z_{i}\in\Gamma}|m_{z_{i}}|d_{0}(\tl u; z_{i}) \\
&\hskip1.25in
+\sum_{\substack{ z_{i},z_{j}\in\Gamma \\ z_{i}\ne z_{j} \\ \gamma_{z_{i}}=\gamma_{z_{j}} \\ m_{z_{i}}m_{z_{j}}>0 }}
          m_{z_{i}}m_{z_{j}}
          \bp{
               \max_{\ell=i, j}\br{\tfrac{\alpha(\tl u; z_{\ell})}{|m_{z_\ell}|}}
               -\max_{\ell=i, j}\br{\tfrac{\wind(\tl u; z_{\ell})}{|m_{z_\ell}|}}
          } \\
&\hskip1in=d_{0}(\tl u)+\sum_{z_{i}\in\Gamma}\Done(\tl u; z_{i})+\cover(\tl u; z_{i})-\cov(e_{1}(\tl u; z_{i}))  \\
&\hskip1.25in
+\sum_{\substack{ z_{i},z_{j}\in\Gamma \\ z_{i}\ne z_{j} \\ \gamma_{z_{i}}=\gamma_{z_{j}} \\ m_{z_{i}}m_{z_{j}}>0 }}
          m_{z_{i}}m_{z_{j}}
          \bp{
               \max_{\ell=i, j}\br{\tfrac{\alpha(\tl u; z_{\ell})}{|m_{z_\ell}|}}
               -\max_{\ell=i, j}\br{\tfrac{\wind(\tl u; z_{\ell})}{|m_{z_\ell}|}}
          } \\
&\hskip1in=d_{0}(\tl u)+\cover(\tl u)-\#\Gamma+\sum_{z_{i}\in\Gamma}\Done(\tl u; z_{i})+1-\cov(e_{1}(\tl u; z_{i}))  \\
&\hskip1.25in
+\sum_{\substack{ z_{i},z_{j}\in\Gamma \\ z_{i}\ne z_{j} \\ \gamma_{z_{i}}=\gamma_{z_{j}} \\ m_{z_{i}}m_{z_{j}}>0 }}
          m_{z_{i}}m_{z_{j}}
          \bp{
               \max_{\ell=i, j}\br{\tfrac{\alpha(\tl u; z_{\ell})}{|m_{z_\ell}|}}
               -\max_{\ell=i, j}\br{\tfrac{\wind(\tl u; z_{\ell})}{|m_{z_\ell}|}}
          },
\end{align*}
and we therefore have
\begin{align*}
I(\tl u)&=\wind_{\pi}(\tl u)+2\sing(\tl u)
     +\sum_{z_{i}\in\Gamma}[\cov(e_{1}(\tl u; z_{i}))-1-\Done(\tl u; z_{i})] \\
     &\hskip.25in
     -\sum_{\substack{ z_{i},z_{j}\in\Gamma \\ z_{i}\ne z_{j} \\ \gamma_{z_{i}}=\gamma_{z_{j}} \\ m_{z_{i}}m_{z_{j}}>0 }}
          m_{z_{i}}m_{z_{j}}
          \bp{
               \max_{\ell=i, j}\br{\tfrac{\alpha(\tl u; z_{\ell})}{|m_{z_\ell}|}}
               -\max_{\ell=i, j}\br{\tfrac{\wind(\tl u; z_{\ell})}{|m_{z_\ell}|}}
          }.
\end{align*}
Further, rewriting
\begin{align*}
2\sing(\tl u)&=2\delta(\tl u)+2\delta_{\infty}(\tl u) \\
&=2\delta(\tl u)+\sum_{z_{i}\in \Gamma}2\delta_{\infty}(\tl u; z_{i})
+\sum_{\substack{z_{i},z_{j}\in\Gamma \\ z_{i}\ne z_{j} }}\delta_{\infty}([\tl u; z_{i}], [\tl u; z_{j}]) \\
&=2\delta(\tl u)+\sum_{z_{i}\in \Gamma}\Done(\tl u; z_{i})+\Dtwo(\tl u; z_{i}) \\
&\hskip.25in
+\sum_{\substack{ z_{i},z_{j}\in\Gamma \\ z_{i}\ne z_{j} \\ \gamma_{z_{i}}=\gamma_{z_{j}} \\ m_{z_{i}}m_{z_{j}}>0 }}
i_{\infty}([\tl u; z_{i}], [\tl u; z_{j}])
+m_{z_{i}}m_{z_{j}}
  \max_{\ell=i, j}\br{\tfrac{\alpha(\tl u; z_{\ell})}{|m_{z_\ell}|}} \\
&=2\delta(\tl u)+\sum_{z_{i}\in \Gamma}\Done(\tl u; z_{i})+\Dtwo(\tl u; z_{i}) \\
&\hskip.25in
+\sum_{\substack{ z_{i},z_{j}\in\Gamma \\ z_{i}\ne z_{j} \\ \gamma_{z_{i}}=\gamma_{z_{j}} \\ m_{z_{i}}m_{z_{j}}>0 }}
i_{\infty}([\tl u; z_{i}], [\tl u; z_{j}])
+m_{z_{i}}m_{z_{j}}
  \max_{\ell=i, j}\br{\tfrac{\winfty(\tl u; z_{\ell})}{|m_{z_\ell}|}} \\
&\hskip.25in
+\sum_{\substack{ z_{i},z_{j}\in\Gamma \\ z_{i}\ne z_{j} \\ \gamma_{z_{i}}=\gamma_{z_{j}} \\ m_{z_{i}}m_{z_{j}}>0 }}
          m_{z_{i}}m_{z_{j}}
          \bp{
               \max_{\ell=i, j}\br{\tfrac{\alpha(\tl u; z_{\ell})}{|m_{z_\ell}|}}
               -\max_{\ell=i, j}\br{\tfrac{\winfty(\tl u; z_{\ell})}{|m_{z_\ell}|}}
           }
\end{align*}
allows us to write
\begin{align*}
I(\tl u)&=\wind_{\pi}(\tl u)+2\delta(\tl u)
     +\sum_{z_{i}\in\Gamma}[\cov(e_{1}(\tl u; z_{i}))-1+\Dtwo(\tl u; z_{i})] \\
&\hskip.25in
+\sum_{\substack{ z_{i},z_{j}\in\Gamma \\ z_{i}\ne z_{j} \\ \gamma_{z_{i}}=\gamma_{z_{j}} \\ m_{z_{i}}m_{z_{j}}>0 }}
i_{\infty}([\tl u; z_{i}], [\tl u; z_{j}])
+m_{z_{i}}m_{z_{j}}
  \max_{\ell=i, j}\br{\tfrac{\winfty(\tl u; z_{\ell})}{|m_{z_\ell}|}}.
\end{align*}
We therefore conclude that
\begin{align*}
\inum(\tl u, \tl u_{c})&\le\wind_{\pi}(\tl u)+2\delta(\tl u)
     +\sum_{z_{i}\in\Gamma}[\cov(e_{1}(\tl u; z_{i}))-1+\Dtwo(\tl u; z_{i})] \\
&\hskip.25in
+\sum_{\substack{ z_{i},z_{j}\in\Gamma \\ z_{i}\ne z_{j} \\ \gamma_{z_{i}}=\gamma_{z_{j}} \\ m_{z_{i}}m_{z_{j}}>0 }}
i_{\infty}([\tl u; z_{i}], [\tl u; z_{j}])
+m_{z_{i}}m_{z_{j}}
  \max_{\ell=i, j}\br{\tfrac{\winfty(\tl u; z_{\ell})}{|m_{z_\ell}|}}.
\end{align*}
with equality occurring for all but a finite number of values of $c\in\R\setminus\br{0}$.
\end{proof}

We now prove \fullref{t:embedded-projection-intro}
 which gives necessary and sufficient conditions for the projection
of a curve to the $3$-manifold to be embedded.

\begin{theorem}[\fullref{t:embedded-projection-intro}]\label{t:embedded-projection}
Let $[\Sigma, j, \Gamma, \tl u=(a,u)]\in\M(M, \Ha, J)$
be a connected, simple pseudoholomorphic curve, and assume
that $\tl u$ does not have image contained in an orbit cylinder.
Then the following are equivalent:
\begin{enumerate}
\item The projected map $u:\Sigma\setminus\Gamma\to M$ is an embedding.

\item The intersection number $\inum(\tl u, \tl u_{c})$ between $\tl u$ and an $\R$-translate
     $\tl u_{c}=(a+c, u)$ is zero for all $c\in\R\setminus\br{0}$.

\item All of the following hold:
     \begin{enumerate}
          \item $u$ does not intersect any of its asymptotic limits.
          \item If $\gamma$ is a periodic orbit so that $u$ is asymptotic at $z\in\Gamma$ to $\gamma^{m_z}$
          and $u$ is asymptotic at $w\in\Gamma$ to $\gamma^{m_w}$, then
          \[
          \tfrac{\winfty(\tl u; z)}{m_z}=\tfrac{\winfty(\tl u; w)}{m_w}.
          \]
     \end{enumerate}

\item All of the following hold:
     \begin{enumerate}
          \item The map $\tl u$ is an embedding.
          \item The projected map $u$ is an immersion which is everywhere transverse to $\X$
          (equivalently $\wind_\pi(u)=0$).
          \item For each $z\in \Gamma$, we have
          \[
          \gcd(m_z, \winfty(\tl u; m_z))=1
          \]
          (equivalently $\cov(e_1(\tl u; z))=1$).
          \item If $\gamma$ is a periodic orbit so that $u$ is asymptotic at $z$ to $\gamma^{m_z}$ and
          $u$ is asymptotic at $w\ne z$ to $\gamma^{m_w}$ with $m_{z}m_{w}>0$, then the asymptotic intersection
          number of the ends $[\tl u; z]$ and $[\tl u; w]$ achieves the bound from \eqref{e:ain-bound-wind-infty}, i.e.
          \[
          \ain^\Phi([\tl u; z], [\tl u; w])=-m_{z}m_{w}\max\br{\tfrac{\winfty^\Phi(\tl u; z)}{m_z}, \tfrac{\winfty^\Phi(\tl u; w)}{m_w}}.
          \]
     \end{enumerate}
\end{enumerate}
\end{theorem}

\begin{proof}
We first address $(2)\iff(3)$.  This follows from \fullref{c:sint-bound-cylindrical} and positivity of intersections
by an argument analogous to that in \fullref{c:self-gin-zero}. We omit the details.

Next we prove $(2)\iff(4)$.  This is also similar to a part of \fullref{c:self-gin-zero}.  It follows from
positivity of intersections and \fullref{l:sint-bound-wind-pi}, particularly the nonnegativity of
the each of the terms on the right hand side of \eqref{e:sint-bound-wind-pi} .
Indeed, from positivity of intersections
and \eqref{e:sint-bound-wind-pi}
we can conclude
that $\inum(\tl u, \tl u_c)=0$ for all $c\in\R\setminus\br{0}$
if and only if:
\begin{itemize}
\item $\wind_\pi(\tl u)=0$,

\item $\delta(\tl u)=0$ (i.e.\ $\tl u$ is an embedding),

\item $\cov(e_1(\tl u; z))=1$ for all $z\in\Gamma$,

\item $\Dtwo(\tl u; z)=0$ for all $z\in\Gamma$, and

\item the asymptotic intersection numbers satisfies
\[
\ain^\Phi([\tl u; z], [\tl u; w])=-m_z m_w\max\br{\tfrac{\winfty^\Phi(\tl u; z)}{m_z}, \tfrac{\winfty^\Phi(\tl u; w)}{m_w}}
\]
for each pair of distinct punctures $(z, w)\in\Gamma\times\Gamma$ at which $u$ is asymptotic
to coverings of a the same underlying orbit with the same sign.
\end{itemize}
Recalling from \fullref{l:wind2} that $\cov(e_1(\tl u; z))=1$ implies that $\Dtwo(\tl u; z)=0$, and
from \fullref{l:covgcd} and \eqref{e:wind-infinity} that 
\[
\cov(e_1(\tl u; z))=\gcd(\wind(\Phi^{-1}e_1(\tl u; z)), m_z)=\gcd(\winfty^\Phi(\tl u; z), m_z)
\] 
we see that these conditions are equivalent to those listed in $(4)$ above.

Next we observe that $(1)\Rightarrow (2)$.  This has been previously observed
by Hofer, Wysocki, and Zehnder
in \cite{hwz:prop2}, but
we include the argument here for completeness, since it is short and illustrative.
Indeed if there exists a $c\in\R\setminus\br{0}$ so that $\tl u=(a, u)$ intersects $\tl u_c=(a+c, u)$, then there is
a pair of points distinct points $z$, $w\in\Sigma\setminus\Gamma$ so
that $\tl u(z)=\tl u(w)$ which in turn implies that $u(z)=u(w)$.
Thus $(2)$ failing to hold implies that the projected map has a double point and hence can't be an embedding.

Finally we show that  $(2), (3), (4)\Rightarrow (1)$.  Indeed if $(2)$, $(3)$ and $(4)$ are all true, then
the projected curve $u$ is an injective (by $(2)$ and the argument of the previous paragraph)
immersion (by $(3)$), which doesn't intersect any of its asymptotic limits (by $(4)$).
As observed in \cite{hwz:fol}, this with the asymptotic behavior
implies that the map $u$ is an embedding.
\end{proof}

\subsection{Holomorphic open book decompositions}
In this section we examine an application of results from the previous sections to
properties of holomorphic curves in manifolds admitting a holomorphic open book decomposition.
Recall from the introduction
$(M, \Ha, J)$ is said to admit a stable, holomorphic open book decomposition if
there is a link $L\subset M$ made up of
elliptic periodic orbits of $\X$, and
a fibration $\pi:M\setminus L\to S^{1}$
so that for any $\theta\in S^{1}$, $\pi^{-1}(\theta)$ is an embedded
surface bounded by $L$, for which we can write
\[
\pi^{-1}(\theta)=u(S^{2}\setminus\Gamma)
\]
for some punctured $\tildeJ$-holomorphic sphere $[S^{2}, i, \Gamma, (a, u)]\in\M(M, \Ha, J)$
with Fredholm index $\ind(\tl u)=2$ 
having only positive punctures asymptotic to simply covered orbits.

\begin{theorem}[\fullref{t:openbook-intro}]
Assume that $(M,\Ha, J)$ admits a stable, holomorphic open book decomposition.
Let $[\Sigma, j, \Gamma, \tl u=(a, u)]\in\MM$ be a connected
pseudoholomorphic curve, and assume that the image of $u$ is not a page of the open book decomposition,
and that $\tl u$ does not have image contained in an orbit cylinder.
Then at least one of the following is true:
\begin{enumerate}
\item At least one of the positive punctures of $\tl u$ limits to an orbit
that is not a binding of the open book
decomposition.
\item At least one of the positive punctures of $\tl u$ limits to a multiple cover of a binding orbit of
the open book.
\end{enumerate}
\end{theorem}

\begin{proof}
Let $\tl v=(b, v):S^2\setminus\Gamma'\to\R\times M$ be any simple pseudoholomorphic curve which projects
to a page $v(S^2\setminus\Gamma')$ of the open book decomposition.
We apply \eqref{e:gin-cylindrical} with the roles of $\tl u$ and $\tl v$ reversed to find that
\[
\gin{\tl u}{\tl v}=
\sum_{z\in\Gamma_+}
m_{z}\bp{
\inum(\gamma_z, v)
+\sum_{\substack{w\in\Gamma' \\ \gamma_w=\gamma_z}}
\bbr{\max\br{\tfrac{\alpha(\gamma_z^{m_z})}{m_z}, \alpha(\gamma_w)}
-\alpha(\gamma_w)}+d_0(\tl v; w)
}
\]
where we've used the assumption that $\tl v$ only has positive, simply-covered punctures.
Moreover, by the assumption that $\ind(\tl v)=2$, we have from \eqref{e:wind-pi-ind} that
\[
0\le\wind_\pi(v)+d_0(\tl v)\le \ind(\tl v)-\chi(S^2)+\#\Gamma_{even}(\tl v)=2-2+0=0
\]
so we conclude that $d_0(\tl v)=0$
(and that $\wind_\pi(v)=0$, but this already follows from \fullref{t:embedded-projection}
since $v$ parametrizes a page of the open book decomposition and is thus an embedding).
Consequently, our formula for the generalized intersection number of $\tl u$ and $\tl v$
simplifies to
\begin{equation}\label{e:gin-openbook}
\gin{\tl u}{\tl v}=
\sum_{z\in\Gamma_+}
m_{z}\bp{
\inum(\gamma_z, v)
+\sum_{\substack{w\in\Gamma' \\ \gamma_w=\gamma_z}}
\bbr{\max\br{\tfrac{\alpha(\gamma_z^{m_z})}{m_z}, \alpha(\gamma_w)}
-\alpha(\gamma_w)}
}.
\end{equation}
Now, the assumption that $u$ does not have image lying in a page of the open book decomposition, implies that
$u$ intersects some page of the open book decomposition, which in turn, by \fullref{t:gin-positivity},
implies that $\tl u$ has positive generalized intersection number with a pseudoholomorphic curve which
projects to that page of the open book decomposition.  Moreover, by homotopy invariance of the
generalized intersection number, $\tl u$ has positive generalized intersection number with
\emph{every} page of the open book decomposition.
Using this observation with the
the formula \eqref{e:gin-openbook} for the generalized intersection number of $\tl u$ with a page of the open book,
we see that $\gin{\tl u}{\tl v}>0$ implies that either:
\begin{enumerate}
\item there exists a $z\in\Gamma$ so that $\inum(\gamma_z, v)>0$, or

\item there exists a $z\in\Gamma$ and $w\in\Gamma'$ so that $\gamma_z=\gamma_w$ and
\[
\max\br{\tfrac{\alpha(\gamma_z^{m_z})}{m_z}, \alpha(\gamma_w)}- \alpha(\gamma_w)>0,
\]
which in turn implies that $m_z>1$.
\end{enumerate} 
These two conditions are equivalent to the two listed in the statement of the theorem, so this completes the proof.
\end{proof}

\begin{remark}
We remark that in the event that the first alternative of the preceding theorem does not hold, i.e.\
when all the positive punctures of the curve $\tl u$ limit to coverings of binding orbits
of the open book decomposition, the proof can be refined
to give a bound on the covering numbers of the
punctures guaranteed by the second alternative in terms of the iteration formula for the Conley--Zehnder
index of the orbit in question.
Indeed, given an elliptic periodic orbit $\gamma\in\Pe\MH$
and a trivialization $\Phi$ of $\xh|_{\gamma}$, \fullref{l:alpha-iterate}
guarantees that there is an irrational $\theta$ so that
\[
\alpha^{\Phi}(\gamma^{k})=\fl{k\theta}.
\]
While the number $\theta$ here depends on the choice of trivialization, 
the change of trivialization formula from \fullref{l:alpha-triv} allows us to conclude that
the fractional part
\[
\br{\theta}:=\theta-\fl{\theta}\in(0, 1)
\]
of $\theta$ is independent of the choice of trivialization.
We will refer to the number $\br{\theta}$ as the \emph{rotation index} of $\gamma$
and denote it by $\rot{(\gamma)}$.

Now,  assume that
$[S^{2}, i, \Gamma', \tl v=(b, v)]$ is a page of the open book decomposition, and
$[\Sigma, j, \Gamma, \tl u]$ is a pseudoholomorphic curve with every positive puncture
limiting to a covering of a binding orbit of the open book decomposition.
Then the terms of the form $\inum(\gamma_{z}, v)$ in formula \eqref{e:gin-openbook} all vanish.
Writing the remaining terms of the formula using rotation indices
we find that
\[
\gin{\tl u}{\tl v}
=\sum_{\substack{(z, w)\in\Gamma_{+}\times\Gamma' \\ \gamma_{z}=\gamma_{w}}}
\fl{m_{z}\rot(\gamma_{z})}
\]
which the reasoning of the proof allows us to conclude must be positive.
We can thus conclude that for at least one $z\in\Gamma$,
$m_{z}> 1/\rot(\gamma_{z})$.
\end{remark}

\subsection{Generalized holomorphic curves}
In this section we prove the that the generalized intersection number
of two generalized pseudoholomorphic curves with no common components
in nonnegative.

Recall from the introduction that a generalized pseudoholomorphic map in $\R\times M$
is a quintuple $(\Sigma, j, \Gamma, (a, u), \nu)$ satisfying
\begin{gather*}
\pi_{\xh}\circ du\circ j=J\circ \pi_{\xh}\circ du \\
u^*\lambda\circ j=da+\nu \\
d\nu=d(\nu\circ j)=0,
\end{gather*}
and a generalized pseudoholomorphic curve
$[\Sigma, j, \Gamma, (a, u), \nu]$
is an equivalence class of
generalized pseudoholomorphic maps, where two maps are considered equivalent if they
differ by holomorphic reparametrization of the domain.

A key observation from \cite{h2000} is that if
$[\Sigma, j, \Gamma, (a, u), \nu]$ is a generalized pseudoholomorphic curve,
then on any simply connected subset $U$ of $\Sigma$, the $M$-part $u$ of
the map has a pseudoholomorphic lift.  Indeed, if
$h:U\to \R$ satisfies $dh=\nu|_{U}$ then
$(a+h, u):U\to\R\times M$ is $\tildeJ$-holomorphic.  A consequence of this is the following,
again first observed in \cite{h2000}, which states that 
for generalized pseudoholomorphic curves, one still obtains
the first term of the asymptotic formula from \fullref{t:asymp}.

\begin{lemma}\label{l:gen-hol-asymp}
Let $(\Sigma, j, \Gamma, \tl u=(a,u), \nu)$ be a generalized pseudoholomorphic map
in $\R\times M$ with no component having image contained in an orbit cylinder,
and assume at $z\in\Gamma$, $u$ is asymptotic
to $\gamma^{m}$.
Then there is an embedding
\[
\psi:[R,\infty)\times S^{1}\to\Sigma\setminus\br{z}
\]
satisfying $\lim_{s\to\infty}\psi(s,t)=z$
so that
\[
\tl u(\psi(s,t))=\bp{m\tau s, \exp_{\gamma^{m}(t)}e^{\lambda s}[e(t)+r(s,t)]}
\]
where $\lambda<0$ is an eigenvalue of $\A_{\gamma^{m}, J}$, $e$ is an eigenvector
of $\A_{\gamma^{m}, J}$ with eigenvalue $\lambda$, and $r(s,t)\to 0$
exponentially as $s\to\infty$.
\end{lemma}

As a consequence of this Lemma, we observe that the quantities
$\winfty^{\Phi}$ and $d_{0}=\winfty^{\Phi}-\alpha^{\Phi}$, defined in 
\eqref{e:wind-infinity} and \eqref{e:def} respectively are well-defined
for ends of generalized pseudoholomorphic curves, provided the images of the maps are not
orbit cylinders.

We can now prove \fullref{t:gen-hol-pos-intro} which we restate here.

\begin{theorem}[Global positivity of intersections for generalized holomorphic curves]
Let $[\Sigma , j, \Gamma, \tl u=(a, u), \nu]$, $[\Sigma', j', \Gamma', \tl v=(b, v), \nu']\in\Mhar$
be generalized pseudoholomorphic curves, and assume that no
component of $\tl u$ or $\tl v$ is contained in an orbit cylinder.  Then
\[
\gin{\tl u}{\tl v}\ge 0.
\]
\end{theorem}

\begin{proof}
As a consequence of \fullref{l:gen-hol-asymp}, we know that
$\tl u$ and $\tl v$ have ends that wind, so
\fullref{t:gin-cylindrical} tells us that the generalized intersection number of
$\tl u$ and $\tl v$ is given by
\begin{align*}
\gin{\tl u}{\tl v}
&
=\sum_{w\in\Gamma'_+}
|m_w|\bp{
\inum(\gamma_w, u)+\sum_{\substack{z\in\Gamma_+ \\ \gamma_z=\gamma_w}}
|m_z|\;\Mxx{ \tfrac{\alpha(\gamma_w^{m_w})}{|m_w|}, \tfrac{\alpha(\gamma_z^{m_z})}{|m_z|}}
+d_0(\tl u; z)
} \\
&
\hskip.125in+\sum_{z\in\Gamma_-}
|m_z|\bp{
\inum(\gamma_z, v)+\sum_{\substack{w\in\Gamma'_- \\ \gamma_w=\gamma_z}}
|m_w|\; \Mxx{\tfrac{\alpha(\gamma_z^{m_z})}{|m_z|}, \tfrac{\alpha(\gamma_w^{m_w})}{|m_w|}}
+d_0(\tl v; w)
} \\
&
\hskip.125in+\sum_{\substack{(z, w)\in\Gamma_{-}\times\Gamma'_{+} \\ \gamma_z=\gamma_w}}
|m_w|d_0(\tl u; z)+|m_z|d_0(\tl v; w)-|m_{w}m_{z}|\bbr{\tfrac{\alpha(\gamma_z^{m_z})}{|m_z|}+\tfrac{\alpha(\gamma_w^{m_w})}{|m_w|}},
\end{align*}
with $\Mx$ as defined in \eqref{e:mxx-define}.
The theorem would follow immediately if we knew that
each term appearing in this formula were nonnegative.
The nonnegativity of the $d_{0}$ terms
is an immediate consequence of the asymptotic description from \fullref{l:gen-hol-asymp}
and reasoning identical to that in \fullref{l:wind-infinity}.
Every other term has previously been shown to be nonnegative in the
proof of \fullref{c:gin-zero} except for
the terms $\inum(\gamma_{w}, u)$ and $\inum(\gamma_{z}, v)$ since in \fullref{c:gin-zero}
we were assuming that $\tl u$ and $\tl v$ were pseudoholomorphic, while here we only
assume generalized pseudoholomorphic.
It remains true in this case that intersections of a projected generalized pseudoholomorphic curves
with a periodic orbit of $\X$ always occur with positive local intersection index.
Indeed, assume at $z_{*}\in\Sigma$ that $u(z_{*})=\gamma(t)$ where $\gamma$ is a periodic orbit
or $\X$.  Then let $\bar u=(a+h, u):U\to\R\times M$ be a local pseudoholomorphic lift of $u$ on
a neighborhood $U$ of $z_{*}$.
Then $\bar u(z_{*})$ intersects the orbit cylinder $\tl \gamma$.
By the assumption that
$\tl u$ does not have any components with
image contained in an orbit cylinder, it follows that the intersection of
$\bar u$ with $\tl \gamma$ is isolated and has positive local index, and it follows as in
\eqref{e:int-cylinder} that the intersection
of $\gamma$ with $u$ is isolated and has positive local index.
\end{proof}

\appendix

\section{Zeroes of $\pi_{\xh}\circ du$}\label{a:windpi}
Here we will prove that for a 
connected pseudoholomorphic curve $[\Sigma, j, \Gamma, (a, u)]\in\MM$,
the projection of the derivative of $u$ onto the hyperplane distribution
$\xh$ either vanishes identically or has a finite number of isolated zeroes of finite positive order.
This is proved in \cite{hwz:prop2} in the case that the Hamiltonian structure comes from
a contact form, and the proof here is an adaptation of the argument given there.

We start with a local coordinate lemma which is a straightforward modification of the well known version of
Darboux's theorem for presymplectic manifolds.
Because it may be of independent interest,
we prove the result for stable Hamiltonian structures on manifolds of arbitrary odd dimension.
In the definition of stable Hamiltonian structure for 
a $2n+1$-dimensional manifold $M$, the condition \ref{h:vform} needs to be changed to
\[
\lambda\wedge \omega^{n}>0
\]
but otherwise remains the same.  The definitions of $\X$ and $\xh$ are identical, and it remains
true that $\lambda$ and $\omega$ are preserved by the flow of $\X$.
For more details see e.g.\ \cite{behwz,sief08}.
For the following we equip $\R^{2n+1}$ with the coordinates
$\br{(z, x, y)}=\br{(z, x_{i}, y_{i})}\in\R\times\R^{n}\times\R^{n}$.

\begin{lemma}[Darboux's theorem for stable Hamiltonian structures]\label{l:ham-darboux}
Let $(M, \Ha)$ be a closed, $2n+1$-dimensional manifold equipped with a stable
Hamiltonian structure $(\lambda, \omega)$.  
For any $p_{0}\in M$, there exists an $\ep>0$ and an embedding
\[
\phi:(-\ep, \ep)^{2n+1}\to M
\]
with $\phi(0)=p_{0}$ and
\begin{gather*}
\phi^{*}\lambda=dz -\sum_{i=1}^{n}\bp{g_{i}(x,y)dx_{i}+h_{i}(x,y)dy_{i}} \\
\phi^{*}\omega=\omega_{0}=\sum_{i=1}^{n}dx_{i}\wedge dy_{i}
\end{gather*}
where the $g_{i}, h_{i}:(-\ep, \ep)^{2n}\to\R$ are smooth real-valued functions satisfying
$g_{i}(0)=h_{i}(0)=0$.
\end{lemma}

\begin{proof}
Let $\psi:\R\times M\to M$ denote the flow of $\X$, i.e.\ $\dot\psi_{t}(p)=\X\circ\psi_{t}(p)$ for all
$(t,p)\in\R\times M$.  Let $\{e_{i}, f_{i}\}$ be a symplectic basis for $(\xh_{p_{0}}, \omega)$, 
that is, assume that
\[
\omega(e_{i}, e_{j})=\omega(f_{i}, f_{j})=0
\]
for any $i, j$, and that
\[
\omega(e_{i}, f_{j})=\delta_{ij}.
\]
Define a map
$\phi_{0}:\R^{2n+1}\to M$ by
\[
\phi_{0}(z, x_{i}, y_{i})=\psi_{z}(\exp_{p_{0}}(\sum_{i=1}^{n}x_{i}e_{i}+y_{i}f_{i}))
\]
where $\exp$ is the exponential map of any Riemannian metric on $M$.
It follows from the definition that $\phi_{0}$ satisfies:
\begin{itemize}
\item
$d\phi_{0}(z, x, y)\partial_{z}=\X$ for all $(z,x,y)\in\R^{2n+1}$,
\item $d\phi_{0}(0)\partial_{x_{i}}=e_{i}$ and $d\phi_{0}(0)\partial_{y_{i}}=f_{i}$, and thus
\item
$\phi_{0}$ is an embedding on some neighborhood of the origin,
\item $\phi_{0}^{*}(0)\omega=\omega_{0}$, and $\phi_{0}^{*}(0)\lambda=dz$.
\end{itemize}
Since the flow of $\X$ preserves $\lambda$ and $\omega$, we can conclude that $\phi_{0}^{*}\lambda$
and $\phi_{0}^{*}\omega$ are independent of the $z$-variable.
Since $\phi_{0}^{*}\lambda(\partial_{z})=\lambda(\X)=1$, we thus have 
\begin{gather*}
\phi_{0}^{*}\lambda=dz -\sum_{i=1}^{n}\bp{\tl g_{i}(x,y)dx_{i}+\tl h_{i}(x,y)dy_{i}} \\
\phi_{0}^{*}\omega=\omega_{0}+r_{(x,y)}
\end{gather*}
where $\tl g_{i}$, $\tl h_{i}$ are smooth, real-values functions on some neighborhood of the origin in $\R^{2n}$
satisfying  $\tl g_{i}(0)=\tl h_{i}(0)=0$, and where $r_{(x,y)}$ is a two-form on $\R^{2n}$ satisfying
$r_{(0,0)}=0$.
The result then follows from applying a Moser trick in the $x$ and $y$ variables
(see e.g\ Theorem 1, in Section 1.3 of \cite{hz}).
\end{proof}

Next we show that  for
$[\Sigma, j, \Gamma, (a, u)]\in\MM$,
$\pi_{\xh}\circ du$ can only have isolated zeroes of positive order.
It suffices to prove this in ``Darboux coordinates'' provided by the preceding lemma.
In the following $\DD_{\ep}$ will denote the disk of radius $\ep$ in $\C=\br{s+it}$.
We observe this argument readily generalizes to higher dimensions, but is of most use in dimension
$3$ in which case the algebraic count of zeroes of $\pi_{\xh}\circ du$ is topologically controlled.

\begin{lemma}
Let $g, h:\R^{2}\to\R$ be smooth functions, and let
$J\in\J(\R^{3}, \Ha_{0})$ where $\Ha_{0}=(\lambda_{0}, \omega_{0})$ is the Hamiltonian structure defined by
\begin{gather*}
\lambda_{0}=dz-g(x,y)\,dx-h(x,y)\,dy\\
\omega_{0}=dx\wedge dy.
\end{gather*}
Moreover, let $\pi:\R X_{\Ha_{0}}\oplus\xi^{\Ha_{0}}\to\xi^{\Ha_{0}}$ be the projection onto
$\xi^{\Ha_{0}}=\ker\lambda_{0}$ along $X_{\Ha_{0}}$.
If $\tl u=(a, u):(\mathbb{D}_{\ep}, i)\to (\R\times \R^{3}, \tildeJ)$
is a pseudoholomorphic map, then $\pi\circ du$ either vanishes identically, or has isolated zeroes of
finite positive order.
\end{lemma}

\begin{proof}
We will show that in an appropriate basis for $\xi^{\Ha_{0}}$ the section $(\pi\circ du)(\ds)$ satisfies
a perturbed Cauchy--Riemann equation.  The result will then follow from the similarity principle
(see e.g.\ Appendix A.6 in \cite{hz}).

We first note that the vector fields $e=\partial_{x}+g\,\partial_{z}$ and $f=\partial_{y}+h\,\partial_{z}$
form a basis for $\xi^{\Ha_{0}}$
and that $X_{\Ha_{0}}=\partial_{z}$.
Moreover, for any vector field $v=v^{x}\,\partial_{x}+ v^{y}\,\partial_{y}+v^{z}\,\partial_{z}$ on $\R^{3}$ we have that
\begin{align*}
\pi v&=v-\lambda_{0}(v)X_{\Ha_{0}} \\
&=v^{x}\,\partial_{x}+ v^{y}\,\partial_{y}+( v^{x}g+v^{y}h)\,\partial_{z} \\
&=v^{x}e+v^{y}f
\end{align*}
so the coordinates of $\pi v$ in the basis $\br{e,f}$ are given by the
$x$ and $y$ components of $v$ in the standard
basis for $\R^{3}$.

Now, writing $u=(u^{z}, u^{x}, u^{y})\in\R^{3}$ and $v=(u^{x}, u^{y})$, the equation
$\pi\circ du\circ i=J\circ\pi\circ du$ applied to $\ds$ and expressed in the basis
$\br{e,f}$ becomes
\[
v_{s}+\barJ(s,t)v_{t}=0.
\]
Here $\barJ(s,t)$ is $J(u(s,t))$ represented in the basis $\br{e,f}$, and hence satisfies $\barJ^{2}=-I$.
Letting $w=v_{s}$ and differentiating the above equation with respect to $s$ leads to
\[
w_{s}+\barJ(s,t)w_{t}+A(s,t)w=0
\]
with
\[
A(s,t)=\barJ_{s}(s,t)\barJ(s,t).
\]
As explained in the first paragraph, the result is now an easy consequence of the similarity principle.
\end{proof}

\begin{corollary}\label{c:pi-du}
Let $[\Sigma, j, \Gamma, \tl u=(a,u)]\in\MM$
be a  connected pseudoholomorphic curve for some cylindrical $\tildeJ$
associated to a $J\in\J(M,\Ha)$, and let
$\pi_{\xh}:TM=\R\X\oplus\xh\to\xh$ be the projection onto $\xh$ along $\X$.
Then $\pi_{\xh}\circ du$ either vanishes identically or has a finite number of isolated zeroes each of finite
positive order.
\end{corollary}

\begin{proof}
Assume that $\pi_{\xh}\circ du$ does not vanish identically.
By the previous two lemmas, the zeroes of $\pi_{\xh}\circ du$ must be isolated and of finite positive order.
Moreover, according to \fullref{c:pi-du-punctures},
$\pi_{\xh}\circ du$ is nonvanishing in some neighborhood of each puncture since
we assume it doesn't vanish identically on $\Sigma\setminus\Gamma$.
We conclude that $\pi_{\xh}\circ du$ has a finite number of isolates zeroes of finite positive order.
\end{proof}

\section{Local intersections of generalized pseudoholomorphic curves}\label{a:gen-hol}
In this appendix, we construct local examples of generalized pseudoholomorphic curves
exhibiting intersection behavior that can't occur for (genuine) pseudoholomorphic curves.
The author first learned that such examples should exist from C.\ Abbas.

We consider $\R^{3}=\br{(z, x, y)}$ equipped with the stable Hamiltonian structure
$\Ha=(\lambda, \omega)=(dz,dx\wedge dy)$.  Then $\X=\partial_{z}$ and
$\xh=\operatorname{span}\br{\partial_{x}, \partial_{y}}$.
Define $J$ on $\xh$ by $J\partial_{x}=\partial_{y}$.
The projection $\pi_{\xh}:T\R^{3}=\R\X\oplus\xh\to\xh$ is given by
\[
(z, x,y)\mapsto(x,y).
\]
Letting $\DD$ denote the unit disk in $\C=\br{s+it}$,
the generalized holomorphic curve equations \eqref{e:gen-hol} for a map
$\tl u=(a, u_{1}, u_{2}):\DD\to\R\times\R\times\R^{2}$ reduce to
\begin{equation}\label{e:gen-hol-special}
\left\{
\begin{gathered}
d(du_{1}\circ i-da)=0 \\
d(du_{1}+da\circ i)=0 \\
\partial_{s}u_{2}+J\partial_{t}u_{2}=0.
\end{gathered}
\right.
\end{equation}

Consider maps $\tl u$, $\tl v:\DD\to\R\times\R^{3}$ defined by
\[
\tl u(s,t)=(0, 0, s, t)
\]
and
\[
\tl v(s,t)=(cs,t,s,t)
\]
where $c\in\R$ is constant.
Then it is straightforward to check that $\tl u$ and $\tl v$ satisfy \eqref{e:gen-hol-special},
and we moreover observe that $\tl u$ is genuinely holomorphic.
If $c=0$ then $\tl u$ and $\tl v$ intersect along the line
$\tl u(s, 0)=\tl v(s, 0)=(0,0,s,0)$, which can't happen
for two genuinely pseudoholomorphic maps.
If $c=-1$, then $\tl u$ and $\tl v$ have an isolated transverse intersection
at $\tl u(0,0)=\tl v(0,0)=(0,0,0,0)$,
and the sign of the intersection is given by
\[
\det
\begin{bmatrix}
\partial_{s}\tl u & \partial_{t}\tl u & \partial_{s}\tl v & \partial_{t}\tl v
\end{bmatrix}(0,0)
=
\det
\begin{bmatrix}
0 & 0 & -1 & 0 \\
0 & 0 & 0 & 1 \\
1 & 0 & 1 & 0 \\
0 & 1 & 0 & 1
\end{bmatrix}
=-1
\]
so the local intersection number is $-1$.

%
%
%
\bibliographystyle{gtart}
\bibliography{isectbib}

%
\end{document}